[1994/12/01]% LaTeX date must December 1994 or later

\documentclass[twoside,reqno]{crmp-l} %the "native" format for the series
\usepackage{amssymb}

\let\chapter=\part

\numberwithin{equation}{section}

\newtheorem{Cor}{Corollary}
\newtheorem{Thm}{Theorem}
\newtheorem{Prop}{Proposition}
\newtheorem{Lem}{Lemma}

\theoremstyle{definition}
\newtheorem{Def}{Definition}
\newtheorem{Ex-common}{Example}[section]

\theoremstyle{remark}
\newtheorem*{Rem}{Remark}

\newenvironment{Ex}{\begin{Ex-common}}
 {\edef\qedsymbol{$\blacktriangleleft$}\qed
 \end{Ex-common}}

\renewcommand\ge{\geqslant}
\renewcommand\le{\leqslant}
\let\tildeaccent=\~
\renewcommand\~[1]{\widetilde{#1}}

\def\<{\left<}
\def\>{\right>}
\def\({\left(}
\def\){\right)}
\def\f{\phi}
\let\ssm=\smallsetminus

\newcommand\const{\operatorname{const}}

\def\Mat{\operatorname{Mat}}
\def\crit{\operatorname{crit}}
\renewcommand\:{\colon}
\def\R{{\mathbb R}}
\def\C{{\mathbb C}}
\def\e{\varepsilon}
\def\l{\lambda}
\def\bb#1{{\boldsymbol{#1}}}
\def\pd#1#2{\frac{\partial#1}{\partial#2}}

\def\Re{\operatorname{Re}}
\def\Im{\operatorname{Im}}
\def\Arg{\operatorname{Arg}}
\def\dist{\operatorname{dist}}
\def\diam{\operatorname{diam}}

\let\paragraph=\S
\def\secref#1{\paragraph\ref{#1}}
\def\S{\varSigma}
\def\bdelta{\boldsymbol\delta}

\def\setheight{\vrule height15pt depth 10pt width 0pt}
\def\twolines#1#2{$
 \begin{matrix}
 \noalign{\smallskip}\text{#1}
 \\
 \text{#2}
 \end{matrix}
$}
\def\threelines#1#2#3{$
 \begin{matrix}
 \noalign{\smallskip}
 \text{#1}
 \\
 \text{#2}
 \\
 \text{#3}
 \end{matrix}
$}
\makeatletter
\def\invddots{{\mathinner{\mkern1mu%
\raise-1\p@\hbox{$\cdot$} \mkern2mu \raise2\p@\hbox{$\cdot$} \mkern2mu%
\raise5\p@\vbox to 7pt{\vss\kern7\p@\hbox{$\cdot$}} \mkern1mu}}}%
\makeatother

\def\({\ifmmode\left(\else\textup{(}\fi}
\def\){\ifmmode\right)\else\textup{)}\fi}

\begin{document}
\raggedbottom

\title[Quantitative theory of ODE]
{ Quantitative theory of ordinary differential equations \\ and tangential
Hilbert 16$^{\text{th}}$ problem }

\author
 {Sergei Yakovenko}

\thanks{The research was supported by the Israeli Science Foundation grant
no.~18-00/1.}

\address{Department of  Mathematics\\
 Weizmann Institute of Science\\
 P.O.B.~26, Rehovot 76100\\
 Israel}
\email{{\tt yakov@wisdom.weizmann.ac.il}
\endgraf{\it WWW page\/}:
{\tt http://www.wisdom.weizmann.ac.il/\char'176 yakov/index.html}}
\date{April 2001}


\begin{abstract}
These highly informal lecture notes aim at introducing and explaining
several closely related problems on zeros of analytic functions defined by
ordinary differential equations and systems of such equations. The main
incentive for this study was its potential application to the tangential
Hilbert 16th problem on zeros of complete Abelian integrals.

The exposition consists mostly of examples illustrating various phenomena
related to this problem. Sometimes these examples give an insight
concerning the proofs, though the complete exposition of the latter is
mostly relegated to separate expositions.
\end{abstract}

\maketitle

%\setcounter{tocdepth}{1}
%\tableofcontents

\chapter{Hilbert 16th problem: Limit cycles, cyclicity, Abelian
integrals}

In the first lecture we discuss several possible relaxed formulations of
the Hilbert 16th problem on limit cycles of vector fields and related
questions from analytic functions theory.

\section{Zeros of analytic functions}\label{sec:analytic-functions}

The introductory section presents several possible formulations of the
question about the number of zeros of a function of one variable. All
functions below are either real or complex analytic in their domains,
eventually exhibiting singularities on their boundaries. We would like to
stress that only \emph{isolated} zeros of such functions are counted, so
that by definition a function \emph{identically vanishing} on an open set,
\emph{has no isolated zeros} there.

Exposition goes mostly by examples that are separated from  each other  by
the symbol $\blacktriangleleft$. A few demonstrations terminate by the
usual symbol \qedsymbol.

\subsection{Nonaccumulation and individual finiteness}
A function $f(t)$ real analytic on a finite open interval $(a,b)\subset\R$
may have an infinite number of \emph{isolated} zeros on this interval only
if they accumulate to the boundary points $a,b$ of the latter. Thus the
\emph{finiteness problem} of decision whether or not the given function
$f$ has only finitely many zeros in its domain, is reduced to studying the
boundary behavior of $f$.

In particular, if $f$ is analytic also at the boundary points $a,b$, then
accumulation of infinitely many zeros to these points is impossible and
hence $f$ has only finitely many roots on the interval $(a,b)$. However,
this strong condition of analyticity can be relaxed very considerably.

\begin{Ex}\label{ex:real-powers}
Assume that $f(t)$ defined on $(0,1]$ admits a \emph{nontrivial} asymptotic
expansion of the form
\begin{equation}\label{asymp-expansion}
\begin{gathered}
  f(t)\sim\sum_{k=0}^\infty c_k t^{r_k}\qquad\text{ as }t\to0^+,
  \\
  c_k,r_k\in\R, \qquad r_0<r_1<\cdots<r_k<\cdots
  \end{gathered}
\end{equation}
i.e., the difference between $f$ and a partial sum decreases faster than
the next remaining term (the nontriviality means that not all $c_k$ are
zeros). Then in a sufficiently small semi-neighborhood of $t=0$ the
function $f$ has the same sign as the first nonzero coefficient $c_k$
(since the functions $t^r$ are nonvanishing on $\R_+$), and hence roots of
$f$ cannot accumulate to $0$ implying that $f$ has only finitely many
zeros on $(0,1]$.
\end{Ex}

This example can be generalized for functions admitting asymptotic
expansion in any system of \emph{mutually comparable} functions
$f_1,f_2,\dots$,  when $f_{k+1}=o(f_k)$, provided that all of them keep
constant sign near the boundary point(s).

\begin{Ex}\label{ex:dulac-series}
Any function on $(0,a]$, $a>0$, representable as a finite sum of the form
\begin{equation*}
  \sum_{k,r}h_{kr}(t)\,t^r \ln^{k-1}t,\qquad r\in\R,\ k\in\mathbb N,
\end{equation*}
with the functions $h_{kr}(t)$ real analytic on $[0,a]$ (i.e., including
the boundary $t=0$), cannot have infinitely many roots accumulating to
$t=0$. Indeed, the above sum expands in the monomials $t^r\ln^k t$,
$r\in\R$, $k\in\mathbb Z$, that are naturally lexicographically ordered by
their growth rates as $t\to 0^+$. Clearly, a function admitting an
asymptotic expansion involving terms of such form, also possesses the
finiteness property.

In this example we do not exclude the cases when the expansion is trivial.
However, the \emph{convergence} assumption on $h_{kr}$ implies than in
such cases the function is identically zero and has no isolated roots at
all.
\end{Ex}

\subsection{Parametric families of analytic functions, localization and
cycl\-icity} Consider a function $f=f(t,\l)=f(t,\l_1,\dots,\l_n)$ real
analytic in an open domain $U$ of the space $\R\times\R^n$ (this means
that $f$ can be expanded in a converging Taylor series centered around any
point in $U$). The function $f$ can be considered as an analytic family of
functions $f_\l=f(\cdot,\l)$ defined in variable domains
$U_\l=U\cap(\R\times\{\l\})\subset\R$ depending on $\l$.

Our nearest goal is to formulate the parametric finiteness property and
establish simple sufficient conditions for it, similar to the
non-parametric case.

\begin{Def}
Let $A$ be a point set in $\R$ or $\C$. Everywhere below we denote by $\#A$
the number of \emph{isolated} points (finite or not) of $A$.
\end{Def}

\begin{Def}
We say that the analytic family $f=\{f_\l\}\:U\to\R$ possesses the
\emph{uniform finiteness property}, if the number of isolated zeros of all
functions $f_\l$ in their respective domains $U_\l$ is uniformly bounded by
a constant independent of $\l$.
\end{Def}

In the same way as in the non-parametric context, it would be desirable to
derive the uniform finiteness from some local properties of the family $f$.

\begin{Def}
The \emph{cyclicity} of the family $f$ at a point
$(t_*,\l_*)\in\R\times\R^n$ is the upper limit (finite or not)
\begin{equation}\label{counter}
  \mathcal N_f(t_*,\l_*)=
  \limsup_{\e\to 0^+}\sup_{\|\l-\l_*\|<\e}\#\{t\in U_\l:|t-t_*|<\e,
  \ f_\l(t)=0\}.
\end{equation}
The term comes from the bifurcation theory (see below).
\end{Def}

We stress again that the notation $\#\{\cdots\}$ above means the number of
isolated roots of $f_\l$; the upper limit may be infinite even if all
these numbers are finite.

If $f_{\l_*}\not\equiv0$, then cyclicity of any family containing this
function can be majorized in terms of only the function itself.

\begin{Ex}\label{ex:non-identically}
Suppose that $f$ is analytic at an interior point $(0,0)\in U$, and
$f_0\not\equiv0$. Then the cyclicity of $f$ at the origin is finite.
Moreover, if $t=0$ is an isolated root of \emph{multiplicity} $\mu$, that
is, $f(t)=ct^\mu+\cdots$, $c\ne 0$, then the cyclicity is no greater than
$\mu$, that is, $\mathcal N_f(0,0)\le\mu$.

This follows from the fact that the $\mu$th derivative of $f$ is
nonvanishing at the origin and hence at all sufficiently close points. But
a function whose derivative $f^{(\mu)}(t)$ has a constant sign on an
interval, cannot have more than $\mu$ isolated zeros (even counted with
multiplicities), as follows from the iterated Rolle theorem.

This example can be easily generalized for the case of \emph{complex
analytic functions} defined in a domain $U\subset \C\times\C^n$. In the
above assumption $f^{(\mu)}_0(0)\ne0$, the cyclicity of the family will be
\emph{exactly} $\mu$. To prove this, one can choose a small circle around
the origin and apply the Rouch\'e theorem to it. Another possibility would
be to use the Weierstrass preparation theorem.
\end{Ex}

\begin{Ex}\label{ex:one-param-func}
If $\dim\l=1$, then any function analytic in $(t,\l)$ near $t=t_*$,
$\l_*=0$ can be expanded as $f(t,\l)=f_0(t)+\l f_1(t)+\l^2f_2(t)+\cdots$.
If $f(t,\l)\not\equiv0$, then for some finite $k$ necessarily
$f_k(t)\not\equiv0$, and after division by $\l^k$ the question about
\emph{isolated} zeros of $f$ can be reduced to the situation when
$f(t,0)\not\equiv0$, discussed earlier. In this case $\mathcal N_f(t_*,0)$
is no greater than the multiplicity of $f_k$ at $t_*$, where $f_k$ is the
first nonzero term in the expansion.
\end{Ex}

This example illustrates an absolutely general fact about analytic
functions (no matter, real or complex): the cyclicity $\mathcal N_f(t,\l)$
takes finite values at all \emph{interior} points of the domain of
analyticity $U$.

\begin{Thm}\label{thm:fcyc-anal}
If $f$ is analytic at a point $(t,\l)\in\R\times\R^n$, then $\mathcal
N_f(t,\l)<+\infty$.
\end{Thm}

This assertion can be derived from general finiteness properties of
analytic sets, see \cite{loja}. In the present form the theorem was
formulated in connection with bifurcations of limit cycles,
see~\cite{roussarie:h16}.

The proof of finiteness of $\mathcal N_f$ in the general multiparametric
case requires some analytic techniques. One possibility---assuming for
simplicity the point to be at the origin $(0,0)$---is to consider the
expansion $f(t,\l)=\sum_{k\ge 0}a_k(\l)t^k$ and the ideal in the ring
$\mathfrak R$ of analytic germs at $(\R^n,0)$, generated by the
coefficients $a_k(\l)$. This ideal is called the \emph{Bautin ideal}
\cite{roussarie:h16}. Since the ring $\mathfrak R$ is Noetherian, the
Bautin ideal is in fact generated by a finite number of the coefficients
$a_k$. If the germs $a_1,\dots,a_\nu$ generate the Bautin ideal for some
finite $\nu$, then one can show that the function $f(t,\l)$ can be
represented as $\sum_{k=0}^\nu a_k(\l)t^kh_k(t,\l)$ with $h_k$ analytic in
$t,\l$ and $h_k(0,0)\ne0$, see \cite{roussarie:h16}. From this it is
already easy to derive that $\mathcal N_f(0,0)\le \nu<\infty$. Almost no
modification is required to cover the complex analytic case as well: the
only difference is that the usual Rolle cannot be applied to holomorphic
non-real functions. The alternative is to use the complex Rolle theorem
from \cite{jdcs-96a}. An example of such use in a situation very similar to
that discussed above, can be found in \cite{nonlin-00}.

\subsection{Finite cyclicity and uniform
 finiteness}\label{sec:finite-cyclicity}
The main reason for introducing the notion of cyclicity is the following
very simple but basic theorem. In applications to bifurcation of limit
cycles of analytic vector fields, it was observed by R.~Roussarie
\cite{roussarie:cyclicity,roussarie:lacets,roussarie:h16}.

Let $f\:U\to\R$ as before be an analytic family of functions and $\mathcal
N_f\:\overline U\to\mathbb N\cup\{+\infty\}$ its cyclicity function.

\begin{Thm}\label{thm:finite-cyclicity}
If  the closure $\overline U$ is a compact subset of  $\R\times\R^n$, and
cyclicity $\mathcal N_f$ is finite everywhere on $\overline U$, then the
family $f=\{f_\l\}$ admits a uniform upper bound on the number of isolated
roots:
\begin{equation}\label{exist-finiteness}
  \sup_\l\#\{t\in U_\l\:f_\l(t)=0\}<+\infty.
\end{equation}
\end{Thm}

\begin{proof}
By definition of the counting function, any point $(t_*,\l_*)\in\overline
U$ can be covered by a sufficiently small cube $\{|t-t_*|<\e,\
\|\l-\l_*\|<\e\}$ such that the number of isolated roots of $f$ in this
cube is no greater than the number $\mathcal N_f(t_*,\l_*)$ finite by the
assumptions of the theorem. It remains only to choose a finite subcover of
the compact $\overline U$ and add together the corresponding cyclicities.
\end{proof}

Theorem~\ref{thm:fcyc-anal} asserts that the cyclicity is automatically
finite at all \emph{interior} points of the domain $U$, so it is again the
boundary behavior of a given parametric family that determines whether this
family possesses the uniform finiteness or not.

\subsection{Terminology: individual \emph{vs}.~existential
finiteness problems} The arguments proving finiteness of the number of
isolated zeros of an individual analytic function $f(t)$ and the arguments
establishing uniform boundedness of this number for parametric families of
such functions, are both of the same purely existential nature: neither of
them gives any way to \emph{compute} or even \emph{estimate from above}
these numbers. However, the parametric claim is definitely stronger than
the assertion concerning individual functions of this family (one can
easily construct families in which the number of isolated zeros is always
finite but not uniformly bounded).

We shall repeatedly encounter the problem on bounding the number of
isolated zeros for various classes of analytic functions and finite
parameter families of such functions, mostly defined by ordinary
differential equations with polynomial right hand sides. For each class
one can pose several finiteness-type problems in the increasing order of
strength (the gap in the enumeration will be explained below).
\begin{enumerate}
 \item[\textbf{1.}] \emph{Individual finiteness problem}. Prove that each function from the
 family possesses only finitely many isolated zeros.
 \item[\textbf{2.}] \emph{Existential \(uniform\) finiteness problem}. Prove that the
 number of isolated roots is uniformly bounded over all functions from
 this family.
 \item[\textbf{4.}] \emph{Constructive finiteness problem}. Find an explicit upper bound for
 the number of isolated roots or at least find an algorithm for computing
 this bound.
\end{enumerate}

\noindent The adjective ``existential'' stresses the difference between the
last two types of problems, whereas the adjective ``uniform'' will be used
to underscore the difference between the first two assertions.

\begin{Ex}\label{ex:analytic-class}
Consider the class of parametric families of functions real analytic on
$[0,1]$ depending analytically on parameters from $[0,1]^n\subset\R^{n+1}$
(in both cases including the boundary). As was just explained, the
individual finiteness for functions from this class follows from the
uniqueness theorem for analytic functions while the existential (uniform)
finiteness theorem (Theorem~\ref{thm:finite-cyclicity}) follows from
Theorem~\ref{thm:fcyc-anal}.
\end{Ex}

\begin{Ex}
Polynomials in one variable of degree $\le d$ for any finite $d$ form a
finite-parameter family. Constructive finiteness theorem for this ``toy''
class of functions is known as the Principal theorem of Algebra. Less
trivial examples can be found in \secref{sec:quasialg}.
\end{Ex}

\subsection{Constructive finiteness}
In practice the existential finiteness is always derived from finite
cyclicity using Theorem~\ref{thm:finite-cyclicity}. Both steps (finite
cyclicity and its globalization) use arguments of existential nature.
However, at least theoretically it may happen that the function $\mathcal
N_f$ is explicitly bounded at all points of the closure $\overline U$.
This still would not allow to compute explicitly the global uniform bound
on the number of isolated zeros of $f$, but such \emph{local constructive
finiteness} would be clearly a much stronger assertion concerning the
family $f$. The corresponding finiteness problem occupies an intermediate
place between existential and (global) constructive finiteness problems.

\begin{enumerate}
  \item[\textbf{3.}]
  \emph{Constructive finite cyclicity problem}. Find an explicit majorant
   for the counting function $\mathcal N_f$.
\end{enumerate}

\noindent Yet it should be remarked that in order to discuss the
constructive finiteness problems (global or local), the family should be
defined by some algebraic data, otherwise computability does not make
sense. We postpone discussion of these question until
\secref{sec:quasialg}.

\begin{Ex}
Let $F_1(t),\dots,F_n(t)$ be analytic functions (real or complex),
satisfying together a system of polynomial ordinary differential equations
\begin{equation}\label{polysys-intro}
  \frac{dF_i}{dt}=\sum_{k+|\alpha|\le d}c_{ik\alpha}\,t^k F^\alpha,\qquad
  F^\alpha=F_1^{\alpha_1}\cdots
  F_n^{\alpha_n},\quad i=1,\dots,n,
\end{equation}
assuming the degree $d$ and all coefficients $c_{ik\alpha}$ (real or
complex) known. Then for any finite $m\in\mathbb N$ the polynomial
combinations of these functions constitute a finite-parametric family: the
coefficients $\l=\{\l_{k\alpha}\:k+|\alpha|\le m\}$ of the polynomial
combinations $f(t,\l)=\sum_{k+|\alpha|\le m}\l_{k\alpha}\,t^kF^\alpha(t)$
being the parameters. The domain of this family is clearly the Cartesian
product of the linear space of the appropriate dimension and the common
domain $U$ of the functions $F_i$. Though solutions of polynomial systems
may blow up in finite time, we will assume for simplicity that all $F_i$
are analytic on the compact closure $\overline U$ of $U$.

In these assumptions both individual and existential finiteness are
guaranteed by Theorem~\ref{thm:finite-cyclicity} which, as usual, does not
provide any bound on the number of zeros. Yet using algebraicity of the
system defining these functions, one can derive explicit bounds. In
\cite{gabr:simplexp} A.~Gabrielov found an upper bound for the maximal
\emph{order of zero} (\emph{multiplicity of a root}) that a nontrivial (not
identically equal to zero) function from the family $f_m$ may have, in
terms of $n,d$ and $m$. The answer is given by an expression polynomial in
the degrees $d,m$ and exponential in $n$ (the dimension of the system).
Y.~Yomdin in \cite{yomdin:oscillation} derived from this multiplicity bound
the corresponding \emph{cyclicity bound}, using arguments involving Bautin
ideal and its generators (a subtle generalization of the approach outlined
in discussion of Theorem~\ref{thm:fcyc-anal}). This combination of the two
results gives a complete solution for the constructive finite cyclicity
problem for polynomial combinations of functions defined by systems of
polynomial ordinary differential equations.
\end{Ex}

The corresponding global problem was settled by D.~Novikov and the author
\cite{lleida,annalif-99}. Unlike the bound for cyclicity, the bound for
the total number of zeros depends on the size of the domain
$\rho=\max_{t\in U}|t|$ where they are counted and, rather naturally, on
the magnitude of the coefficients $R=\max_{i,k,\alpha}|c_{ik\alpha}|$ of
the system \eqref{polysys-intro}. The answer is polynomial in $\rho,R$ but
as a function of the other (integer) parameters $n,m,d$ is a tower of
height $4$. The detailed explanation of this result is postponed until
\secref{sec:meandering}.

\subsection{Paradigm}
All the above results on elementary properties of analytic functions and
their parametric families, were mentioned in order to set out a paradigm
that will be used when discussing much more complicated parallel problems
on limit cycles of analytic and polynomial vector fields. In order to make
the similarity transparent, all finiteness theorems can be organized in the
form of a table with finiteness types labeling the rows and corresponding
theories represented by different columns. The first column describes
theory of analytic functions as in Example~\ref{ex:analytic-class} above.
Cells of this column list the key theorems ensuring each type of
finiteness types.

The following section \secref{sec:lim-cyc} describes some known results
fitting the second column of Table~1: the strongest form of the Hilbert
16th problem should occupy the place where the question marks appear. The
parallel non-constructive results for Abelian integrals will be discussed
in \secref{sec:AI-finiteness}: their complexity level can be described as
intermediate between relatively elementary theorems on analytic functions
and transcendentally difficult counterparts on limit cycles.

\begin{Rem}
In general, the results in the lower part of the table are more difficult
than those from the upper part. However, this is not true when comparing
existential finiteness with constructive finite cyclicity: while the
former deals with boundary behavior, the latter addresses the issues of
completely different nature that can be (and indeed sometimes are)
simpler. In particular, this is the case with Abelian integrals, where the
bound for multiplicity can be be relatively easily obtained, see
\cite{mardesic:multiplicity} and \secref{sec:mardesic} below.
\end{Rem}

% set \arrayrulewidth1pt to get thicker lines

\begin{table}[hbt]
\begin{tabular}{|c||c|c|c|}
 \hline
  &\multicolumn{3}{c}
 {\setheight \bf Category of objects} \vline
 \\
 \cline{2-4}
 \twolines{\bf Finiteness}{\bf type}  &
 \threelines{\it Functions}{\it analytic}{\it on a compact}&
 \threelines{\it Limit cycles}{\it of polynomial}{\it vector fields} &
 \twolines{\it Abelian}{\it integrals}
 \\
\hline
 \twolines{\bf Individual}{\bf finiteness} &
 \twolines{uniqueness}{theorem} &
 \twolines{Ilyashenko (1991)--}{Ecalle (1992)} &
 \twolines{easy exercise}{see \secref{sec:AI-ind}}
 \\
 \hline
 \twolines{\bf Existential}{\bf finiteness} &
 \threelines{finite cyclicity}{Theorems \ref{thm:fcyc-anal}}
    {and~\ref{thm:finite-cyclicity}, see \secref{sec:finite-cyclicity}} &
 \threelines{in progress for}{quadratic case (121}{polycycles \cite{drr:list})}&
 \threelines{Varchenko--}{Khovanskii}{(1985)}
 \\
 \hline
 \threelines{\bf Constructive}{\bf finite}{\bf cyclicity}&
 \twolines{Gabrielov--}{Yomdin (1998)}&
 \threelines{known for $k$-generic }{elementary and few}{other polycycles}&
 \threelines{Marde\v si\'c}{\cite{mardesic:multiplicity}}{(multiplicity)}
 \\
 \hline
 \twolines{\bf Constructive}{\bf finiteness} &
 \threelines{Novikov--}{Yakovenko}{(1999)}
  &
 \twolines{\bf ???}{(nothing in view)} &
 see Table~2
 \\
 \hline
 \end{tabular}
 \bigskip
 \caption{Various flavors of the Hilbert 16th problem}
% \begin{minipage}{0.9\hsize}
%\begin{small}
%Constructive finiteness problem for analytic functions is considered only
%for special classes of functions defined by ordinary differential equations.
% \end{small}
% \end{minipage}
 \end{table}

\section{Limit cycles of planar vector fields}\label{sec:lim-cyc}

This section briefly surveys a few known general finiteness results for
limit cycles of vector fields.

\subsection{Basic facts}
A polynomial differential equation
\begin{equation}\label{pvf}
  \dot x= P(x,y),\quad \dot y=Q(x,y),\qquad P,Q\in\R[x,y]
\end{equation}
defines a vector field on the real plane $\R^2$, spanning the distribution
(line field) with singularities given by null spaces $\{\omega=0\}$ of the
polynomial Pfaffian form $\Omega=Q\,dx-P\,dy$. If $P$ and $Q$ have a
nontrivial greatest common divisor $R$, then this distribution on the
Zariski open set $\{R\ne 0\}$ coincides with that given by the form
$R^{-1}\Omega$ and hence extends analytically to all but finitely many
points of the curve $\{R=0\}$. For the same reasons the distribution
extends onto the infinite line $\R P^1$ on the projective plane $\R P^2$,
having at most finitely many isolated singularities on the latter after
cancellation of an eventual common factor. In other words, when dealing
with an individual differential equation \eqref{pvf}, one can assume that
the corresponding distribution has only \emph{isolated singularities}. It
is sometimes convenient to talk about \emph{foliations with isolated
singularities} defined by the distribution $\{\Omega=0\}$ on $\R P^2$.

A limit cycle $\gamma$ by definition is an isolated periodic trajectory of
\eqref{pvf}. Limit cycles of a vector field can be tracked by their
intersections with analytic arcs transversal to the field. More precisely,
near any intersection $p=\gamma\cap\sigma$ of a periodic orbit $\gamma$
with a transversal arc $\sigma$, the Poincar\'e first return map
$\Delta_\sigma$, the transport along integral curves until the next
intersection with the transversal, is defined. It has $p$ as a fixed
point, \emph{isolated} if $\gamma$ is a limit cycle. Choosing an analytic
chart $t$ on $\sigma$ allows to describe $p$ as a root of the
\emph{displacement function} $f(t)=\Delta_\sigma(t)-t$. All nearby limit
cycles must intersect $\sigma$ by points that are isolated roots of the
displacement $f$.

From the general theorem on analytic dependence of solutions of
differential equations on initial conditions and parameters, the
displacement function $f$ is analytic near $p$. As a corollary, we conclude
that an infinite number of limit cycles (corresponding to isolated roots
of the displacement function) cannot accumulate to a periodic orbit of the
field. In a similar way, if the vector field depends analytically on
parameters $\l$ and exhibits a periodic orbit $\gamma$ (isolated or not)
for one value $\l_*$ of the parameters, then there exists a finite upper
bound for the number of limit cycles in a small annulus around $\gamma$,
uniform over all values of $\l$ sufficiently close to $\l_*$.

In general, one transversal cannot ``serve'' all limit cycles, and there
is no natural way to define the maximal domain of the Poincar\'e return map
unambiguously and globally. However, the discussion in
\secref{sec:analytic-functions} suggests that it is the boundary behavior
that is important for counting isolated zeros of the displacement.

\subsection{Polycycles and limit periodic sets}
Instead of trying to reduce formally the global investigation of limit
cycles to that of return maps for one or several transversal arcs, it is
better to study the cycles (compact leaves of the foliation) themselves.

For the non-parametric case, the analogue of a boundary point of the
domain of the Poincar\'e return map is a \emph{polycycle}, an invariant
set consisting of one or more \emph{singular points} of the vector field,
and a number of bi-infinite trajectories connecting them in a cyclic order
(repetitions of singular points allowed). More accurately, one can show
that the only sets that can appear as Hausdorff limits of periodic orbits
of a vector field having only isolated singular points, are polycycles.

In the parametric case one cannot assume anymore that singularities of the
foliation are isolated (the polynomials $P$ and $Q$ may have common
factors for some values of the parameters and be mutually prime for the
rest). Still one can show \cite{franc-pugh} that the Hausdorff limit of a
family of limit cycles occurring for converging values of the parameters,
must be a \emph{limit periodic set} (also known as \emph{graphic}), the
object differing from a polycycle in only one instance, namely, it may
contain analytic arcs of non-isolated singularities.

\subsection{Individual finiteness theorem and finite
cyclicity}\label{sec:dulac}
 Following the paradigm set out in \secref{sec:finite-cyclicity}, one can
easily formulate counterparts of the individual and existential finiteness
problems for limit cycles for vector fields. The individual finiteness
theorem (known also as the \emph{Dulac conjecture} after being believed
for some 60 years to be the \emph{Dulac theorem}) asserts that \emph{any
polynomial vector field can have only finitely many limit cycles}. This is
an easy corollary to the highly nontrivial \emph{nonaccumulation theorem}
asserting that limit cycles of an analytic vector field cannot accumulate
to any polycycle. The latter assertion is the most spectacular and the
most general fact established so far in connection with the Hilbert 16th
problem, that was independently proved by Ilyashenko
\cite{ilyashenko:finiteness} and Ecalle \cite{ecalle:book} by totally
different methods. This is a typical example of an assertion on boundary
behavior of the Poincar\'e map.

\begin{Rem}
It would be appropriate to notice here that Dulac reduced the individual
finiteness problem to the situation discussed in
Example~\ref{ex:real-powers}.  However, he did not notice that the
expansion \eqref{asymp-expansion} (only involving logarithms) may well be
\emph{trivial}, hence not allowing for so easy treatment. It took both
Ilyashenko and Ecalle hundreds of pages to prove that even in this
apparently highly degenerate case the nonaccumulation still holds.
\end{Rem}

The existential finiteness problem for polynomial vector fields arises very
naturally, since coefficients of the polynomials $P,Q$  can be treated as
parameters. These parameters can be considered as ranging over the compact
projective space of a suitable dimension, since simultaneous
multiplication of both $P$ and $Q$ by a nonzero constant does not alter
the phase portrait. Existence of a uniform upper bound for the number of
limit cycles would mean that \emph{for any degree $d$ the number of limit
cycles for a polynomial vector field of degree $\le d$ is bounded by some
number $\mathcal H(d)$ depending only on $d$}. According to the Roussarie
localization theorem (an analog of Theorem~\ref{thm:finite-cyclicity} for
limit cycles, see \cite{roussarie:cyclicity}), it would be sufficient to
prove that any limit periodic set has finite cyclicity within the
universal family of polynomial vector fields of degree $\le d$.

However, even for the simplest nontrivial case of quadratic ($d=2$) vector
fields, this is not proved. Dumortier, Roussarie and Rousseau composed in
\cite{drr:list} a list of 121 graphics  occurring for quadratic vector
fields, and reduced the existential finiteness problem to proving that all
these graphics have finite cyclicity. First several cases were studied in
the same paper by applying relatively standard tools of bifurcation theory.
Since then many other considerably more delicate cases were investigated
and their finite cyclicity proved, so that one can hope that the
existential finiteness will be ultimately established in the quadratic
case. At the same time it is clear that such type of case study is even
theoretically impossible for higher degree cases.

On this background the perspective of \emph{computing} or even
\emph{explicitly majorizing} the ``\emph{Hilbert number}'' $\mathcal
H(d)$, which is the strongest (constructive) form of the initial Hilbert's
question ``on the number and position of limit cycles'' \cite{hilbert},
looks very remote.

\subsection{Digression: constructive solutions of localized problems}
As was already noted above, though finiteness of cyclicity for all
graphics is not proved, one can independently work towards obtaining
constructive bounds for cyclicity in cases when its finiteness is already
known. This is possible since cyclicity is very strongly depending on the
types of singular points lying on the polycycle (resp., the limit periodic
set).

For example, if all singular points on the polycycle are elementary
(having non-nilpotent linear parts), then cyclicity of \emph{generic}
$r$-parametric families can be estimated in terms of the number $r$ of
independent parameters. The corresponding algorithm was suggested by
Ilyashenko and the author \cite{elementary} and improved to a very concise
and explicit answer by Kaloshin \cite{kaloshin:elementary}. Some other
cases of polycycles (usually carrying one or at most two singular points)
also have known cyclicity, see \cite{kotova:zoo} and references therein
for a synopsis. These results are partly reflected in Table~1.

Yet if the problem is formulated for polynomial vector fields and the
bounds are required to be given in terms of the degree, the problem
immediately becomes transcendentally difficult. Almost nothing is known,
for example, on the maximal multiplicity of a (nonsingular) limit cycle of
a polynomial vector field. Cyclicity of singular points is a considerably
more ``algebraic'' problem. For example, when the linear part of the
singularity is a nondegenerate rotation, a polynomial algorithm can be
suggested which stops at a step number $N$ if and only if the cyclicity of
a singular point is $N$. However, the running time of this algorithm is
absolutely unknown (see the discussion of a similar problem in
\secref{sec:center-focus}). Recently M.~Briskin, J.-P.~Fran\c coise and
Y.~Yomdin treated in details one problem on Abel equations that may give a
clue to understanding some of the phenomena, see e.g.,
\cite{bfy:bautin-abel,bfy:cras}.

The following section describes one very important case when the local
(with respect to parameters) constructive finiteness can be reduced to an
algebraic context.

\section{Abelian integrals: appearance and basic
properties}\label{sec:AI-intro}

Consider a family of analytic vector fields $X_\l$ analytically depending
on parameters $\l\in(\R^n,0)$ (we consider only a sufficiently small
neighborhood of the origin in the parameter space). Assume that for $\l=0$
the vector field $X_0$ possesses a limit cycle $\gamma_0$ of multiplicity
$\mu<+\infty$. As illustrated by Example~\ref{ex:non-identically}, the
cyclicity of $\gamma_0$ in this case is at most $\mu$, being thus
explicitly bounded.

The problem becomes less trivial if the field $X_0$ possesses an annulus
filled by nonisolated periodic orbits, $\gamma_0$ being one of them.
Vector fields with this property are called integrable or conservative.
Such behavior, very unlikely for an arbitrary analytic family, becomes an
event of finite codimension when \emph{polynomial} vector fields are
considered.

To study bifurcations of limit cycles in perturbations of conservative
systems, an approach described in Example~\ref{ex:one-param-func} can be
used.

\subsection{Perturbation of polynomial vector fields}
Assume that the form $\Omega=Q\,dx-P\,dy$ defining the polynomial
distribution $\{\Omega=0\}$, is \emph{exact}: $\Omega=dH$, where $H=H(x,y)$
is a bivariate polynomial in two variables. The corresponding vector field
is then a Hamiltonian one,
\begin{equation}\label{ham-vf}
  \dot x=-\pd Hy,\quad \dot y=\pd Hx.
\end{equation}
with $H$ being its \emph{Hamiltonian function} or simply the
\emph{Hamiltonian}.

Since $dH$ vanishes on the null spaces of the distribution $\{\Omega=0\}$,
each leaf of the corresponding foliation belongs to a level curve of $H$.
In particular, a periodic orbit must be a compact oval of some level curve
and hence all nearby close leaves must be also closed. Hence a Hamiltonian
vector field cannot have limit cycles: the corresponding Poincar\'e return
map is identity and the displacement identically zero.

Consider the following \emph{one-parameter family} of polynomial
distributions, perturbing the Hamiltonian distribution:
\begin{equation}\label{pert-ham}
  dH+\e\omega=0,\qquad \e\in(\R^1,0),\ \omega=p\,dx+q\,dy,\ p,q\in\R[x,y].
\end{equation}
The polynomial 1-form $\omega$ can be arbitrary.

Consider an analytic segment $\sigma$ transversal to the oval
$\gamma(t_*)\subset\{H=t_*\}$. This transversality implies that the
Hamiltonian $H$ restricted on $\sigma$, gives an analytic chart $t$ on the
latter. On the other hand, for all $t$ sufficiently close to $t_*$, one
can unambiguously choose an oval $\gamma(t)\subset\{H=t\}$ so that
$\gamma(t)$ tends to $\gamma(t_*)$ as $t\to t_*$, e.g., in the sense of
Hausdorff distance (this choice in a more broad context is discussed below,
in \secref{sec:topology}).

The following elementary computation gives the \emph{first variation} of
the displacement function $f(t,\e)=\Delta_\sigma(t)-t$ with respect to the
small parameter $\e$ at $\e=0$.

\begin{Lem}[Poincar\'e--Pontryagin]
\begin{equation}\label{poi-pon}
  f(t,0)\equiv0,\quad
  \left.\frac{d}{d\e}\right|_{\e=0}f(t,\e)=-\oint_{\gamma(t)}\omega.
\end{equation}
\end{Lem}

\begin{proof}
Consider the leaf (analytic curve) of the foliation $\{dH+\e\omega=0\}$
passing through the point $t$ on the transversal, and denote by
$\gamma({t,\e})$  the (oriented) segment of this leaf between the the
initial point and the next intersection with $\sigma$. By definition of
the chart $t$, the displacement function measured in the chart $t$, is the
difference of $H$ between the endpoints of $\gamma(t,\e)$, hence
\begin{equation*}
  f(t,\e)=\int_{\gamma({t,\e})}dH
\end{equation*}
(the equality is exact). The form $dH+\e\omega$ vanishes on
$\gamma({t,\e})$, therefore the integral above is equal to the integral of
$-\e\omega$ along ${\gamma({t,\e})}$ and, since $\gamma({t,\e})$ converges
uniformly to the \emph{closed} curve
$\gamma(t)=\gamma({t,0})\subset\{H=t\}$ as $\e\to 0$, we conclude that
$$
f(t,\e)=-\e\oint_{\gamma(t)}\omega+o(\e)\qquad\text{as }\e\to0,
$$
with $o(\e)$ uniform and analytic in $t$ and $\e$. This yields the formula
\eqref{poi-pon} for the derivative.
\end{proof}

Recall that an \emph{Abelian integral} is the integral of a polynomial
1-form along any algebraic oval (i.e., the closed oval of a real algebraic
curve). Clearly, Abelian integrals can be considered as \emph{functions of
the parameters}, the coefficients of the form and the algebraic equation
defining the curve. However, in applications to the bifurcation theory,
one coefficient is clearly distinguished and plays the role of the
argument, while others are treated as ``true'' parameters. Making a minor
abuse of language, we shall always consider Abelian integrals as analytic
functions in the following sense.

\begin{Def}
A \emph{real} Abelian integral corresponding to a Hamiltonian $H$ and a
polynomial 1-form $\omega=p(x,y)\,dx+q(x,y)\,dy$, $H,p,q\in\R[x,y]$, is a
multivalued function $I(t)$ of the real variable $t$ defined by
integration of $\omega$ over a real oval $\gamma(t)$ of the algebraic curve
$\{H(x,y)=t\}\subset\R^2$.
\end{Def}

The reason for multivaluedness is obvious, since there can be several
ovals lying on the same level curve $\{H=t\}$. In \secref{sec:topology} we
describe the branching of the Abelian integral $I(t)$ after analytic
continuation into the complex domain. This complexification of Abelian
integrals explains connections between different branches of these
functions, see \secref{sec:topology}.

\subsection{Bifurcation of limit cycles from
periodic orbits of Hamiltonian systems} Limit cycles appearing in the
one-parametric system \eqref{pert-ham}, correspond to isolated zeros of
the displacement function $f(t,\e)$, real analytic in $t$ and $\e$, whose
first Taylor term was just computed:
\begin{equation}\label{higher-variations}
  f(t,\e)=I_0(t)+\e I_1(t)+\e^2 I_2(t)+\cdots,
  \quad I_0(t)\equiv0,\ I_1(t)=-\oint_{\gamma(t)}\omega.
\end{equation}
This is exactly the situation treated in Example~\ref{ex:one-param-func}:
if the first variation is nontrivial, $I_1(t)\not\equiv0$, then the number
of isolated roots of $f(t,\e)$ on any interval of analyticity (i.e., not
containing singularities at the endpoints) is no greater than the number of
(necessary isolated) roots of the Abelian integral $I_1$ defined by the
polynomial $H$ and the polynomial 1-form $\omega$. Computing or majorizing
the number of isolated zeros of Abelian integrals is the central theme of
these notes: under the name \emph{Tangential Hilbert problem} it is
discussed below.

The first variation method does not work if $I_1(t)\equiv0$. In this case
higher variations $I_k(\cdot)$, $k=2,3,\dots$, have to be computed and
analyzed.

The computation of higher variations is relatively simple. It follows from
results of Ilyashenko \cite{ilyashenko:1969} and Gavrilov
\cite{gavrilov:petrov-modules} that for almost all Hamiltonians $H$ (and
certainly for Hamiltonians with isolated critical points, pairwise
different critical values and transversal to infinity as defined below in
\secref{sec:trans-infinity}), the condition $I_1\equiv0$ implies that
there exist two polynomials $G,F\in\R[x,y]$ such that
\begin{equation}\label{compensator}
  \omega=G(x,y)\,dH+dF(x,y),\quad\text{hence}\quad d\omega=dG\land dH.
\end{equation}
(the inverse statement is obvious).

\begin{Ex}[see \cite{francoise:iterated}]\label{ex:francoise-quadratic}
Consider the Hamiltonian $H(x,y)=\tfrac12 (x^2+y^2)$. In the complex
coordinates $z=x+iy$, $\bar z=x-iy$ we have $H=\tfrac12 z\bar z$. One can
easily verify that the polynomial 1-form $\omega=A(z,\bar z)\,dz+B(z,\bar
z)\,d\bar z$ has identically vanishing integral over the circles $\{H=t\}$
if and only if the differential $d\omega=(-A_{\bar z}+B_z)\,dz\land d\bar
z$ contains {no monomial terms} of the form $(z\bar z)^k\,dz\land d\bar z$.
Any other monomial can obviously be represented in the form
\eqref{compensator}: $z^i\bar z^j\,dz\land d\bar z=dG\land dH$ with
$G={z^i\bar z^j}/{(i-j)}$.
\end{Ex}

Using the representation \eqref{compensator}, the formula for the second
variation of the displacement map can be expressed as an Abelian integral
again.

\begin{Lem}[\cite{bautin,francoise:iterated}]\label{lem:next-variation}
Assume that for the perturbation \eqref{pert-ham} the first variation
vanishes identically, so that
\begin{equation*}
  I_1(t)=\oint_{\gamma(t)}\omega\equiv0,
\end{equation*}
and $G$ is any polynomial satisfying the condition \eqref{compensator}.
Then the second variation of the displacement map is given by the following
Abelian integral,
\begin{equation}\label{second-variation}
  I_2(t)=\left.\frac {d^2}{d\e^2}\right|_{\e=0}f(t,\e)=
  -\oint_{\gamma(t)}G\omega.
\end{equation}
\end{Lem}

This construction can be further iterated \cite{francoise:iterated} as
long as necessary to obtain a variation that is not identically vanishing.
If, on the other hand, all integrals obtained in this recurrent process,
turn out to be identically zero, then the displacement function itself is
identically zero and hence the perturbation \eqref{pert-ham} consists of
integrable systems for all $\e$.

\begin{Rem}
The exact conditions on the Hamiltonian $H$ for the identity
\eqref{compensator} to hold for any form whose integral is identically
zero, involve connectedness of the complex affine level curves $H^{-1}(t)$
for almost all $t\in\C$ \cite{gavrilov:petrov-modules}. Discussion of the
multidimensional situation can be found in \cite{bonnet-dimca}. As soon as
these conditions fail (even because of coinciding critical values, for
example), then computation of higher variations may involve integration of
non-polynomial forms, see \cite{iliev:nonpolynomial}.
\end{Rem}

\subsection{Open problems: bifurcation of limit cycles from Hamiltonian
polycycles and generalized Poincar\'e center--focus
problem}\label{sec:center-focus}
 It would be wrong to think that the above
approach based on computing consecutive variations, \emph{completely}
reduces the question on limit cycles born by perturbation of Hamiltonian
systems, to investigation of Abelian integrals.

First, even assuming the simplest case $I_1\not\equiv0$, we cannot in
general say anything about limit cycles born from \emph{critical} level
curves (corresponding to polycycles of the unperturbed Hamiltonian
system). In the particular case of \emph{separatrix loops} (level curves
homeomorphic to the circle and carrying only one nondegenerate saddle
point), the problem was settled by R.~Roussarie \cite{roussarie:lacets}. He
estimated how many derivatives of the Abelian integral $I(t)$ should have
zero limits as $t\to t_*$ in order for the perturbation \eqref{pert-ham} to
produce $n$ or more limit cycles close to the separatrix loop on the
critical level curve $\{H=t_*\}$. This result yields an upper bound for the
cyclicity of the separatrix loop in terms of the ``multiplicity of the
root'' of the Abelian integral $I_1$ at the point $t_*$ where $I_1$ in fact
loses its analyticity. The only other type of critical level curve that
can occur for a generic polynomial $H$, if the eight-shaped curve (also
carrying a nondegenerate saddle). The answer in this case is not yet known.

The second, apparently much more difficult problem, appears in connection
with the natural question, \emph{How many consecutive variations $I_k$
should be computed in order to guarantee that the perturbation is
non-conservative?}

The inductive process described in Lemma~\ref{lem:next-variation}, is
algebraic in the following sense. Starting from the perturbation
\eqref{pert-ham} given by the form $\omega=\omega_1$, we construct a
sequence of the polynomial 1-forms $\omega_2,\omega_3,\dots$ which express
the higher variations $I_k(t)$ \eqref{higher-variations}, assuming that
all previous variations vanish identically, $I_1\equiv\cdots\equiv
I_{k-1}\equiv0$. Coefficients of the forms $\omega_k$ are given by
algebraic (polynomial) expressions involving coefficients of the initial
form $\omega$. Vanishing of their integrals implies an infinite number of
polynomial identities between coefficients of the initial form $\omega$.
By the Hilbert basis theorem, all these infinitely many identities are
corollaries to only finitely many of them. Thus only finitely many steps
of the inductive process must be performed (their number $N<+\infty$
depends on the \emph{degree} of the initial form $\omega$ and, naturally,
on the Hamiltonian $H$). If all integrals $\oint\omega_k$, $k=1,\dots,N$
are identically zeros, then all higher variations are necessarily zeros and
the entire family \eqref{pert-ham} is integrable.

However, the problem of finding an upper bound for the number $N$ in terms
of $\omega$ and $H$ is overtly open: even for the most simple case
$H=\tfrac 12 (x^2+y^2)$ considered in Example~\ref{ex:francoise-quadratic}
above, the answer is unknown, moreover, it constitutes the challenging
problem on distinguishing between center and focus, posed by Poincar\'e a
century ago and still open. A similar though apparently more simple
problem was recently studied in \cite{bfy:center1}.

\begin{Rem}
A polynomial vector field can be integrable but not Hamiltonian: it is
sufficient that the corresponding form $\Omega$ possess an integrating
factor. A typical example is that of \emph{Darboux integrable} vector
fields. Let $\l_1,\dots,\l_n\in\R$ be a collection of real numbers and
$H_1,\dots,H_n\in\R[x,y]$ polynomials. The rational 1-form
\begin{equation*}
  \Omega'=\sum_{i=1}^n\l_i\frac{dH_i}{H_i}=\frac{dF}{F},\qquad
  F=H_1^{\l_1}\cdots H_n^{\l_n},
\end{equation*}
determines a conservative (integrable) singular foliation, since its null
spaces $\{\Omega'=0\}$ are tangent to the real level curves
$\{F(x,y)=t\}\subset\R^2$ which are all closed unless beginning and ending
at infinity. The form $\Omega'$ is not polynomial (only rational), but the
form $\Omega=R\Omega'$ already will be, if
$R=\gcd(H_1,\dots,H_n)\in\R[x,y]$ is the common divisor of all $H_i$.

Appearance of limit cycles in polynomial perturbations of the form
\begin{equation*}
  \Omega+\e\omega=0,\qquad \omega=p\,dx+q\,dy,\quad p,q\in\R[x,y],
\end{equation*}
is determined by the same mechanisms as in the Hamiltonian case, in
particular, the first variation of the Poincar\'e displacement map measured
in units of $\ln F=\sum \l_i\ln H_i$, is equal to the integral
\begin{equation*}
  I(t)=\oint_{F=t}\frac 1{R(x,y)}\,\omega
\end{equation*}
of the \emph{rational} 1-form $R^{-1}\omega$ along closed level curves of
the \emph{transcendental} function $F$.

As was noted on several occasions by V.~Arnold, the problem on Abelian
integrals should be posed also for such ``pseudo-Abelian'' integrals.
Nevertheless, there are practically no results of general nature pertinent
to this problem, even in the simplest cases $n=2$ and $n=3$. In
particular, \emph{all} results described below (derivation of the
Picard--Fuchs system, monodromy properties etc) fail for this generalized
class of integrals.
\end{Rem}

\section[Tangential Hilbert problem]{Finiteness problems for Abelian integrals:
tangential Hilbert problem}\label{sec:AI-finiteness}

Despite its ambivalent nature, the connection between limit cycles of
polynomial vector fields and isolated zeros of Abelian integrals justifies
formulation of several finiteness problems for the latter. As suggested
above, we refer to this problem as the {tangential Hilbert problem},
distinguishing between several finiteness types according to the paradigm
laid out before.

\subsection{Individual finiteness}\label{sec:AI-ind}
For any fixed combination of the Hamiltonian $H$ and the 1-form $\omega$
the Abelian integral appearing in \eqref{poi-pon} as a function of the
real variable $t$ can be shown to be real analytic except for finitely
many values of $t$. According to \secref{sec:analytic-functions}, it may
have only finitely many isolated zeros unless they accumulate to one of
these exceptional values. The \emph{individual finiteness problem for
Abelian integrals} is to prove that such accumulation is in fact
impossible. This assertion turns out to be a simple corollary to a general
theorem describing ramification of Abelian integrals after analytic
continuation into the complex domain. The key step in the proof of
individual finiteness for Abelian integrals is the following
representation.

\begin{Lem}\label{lem:AI-expansion}
Any Abelian integral near any exceptional point $t_*$ admits a
\emph{converging} representation of the form
\begin{equation}\label{asymp-expan}
  I(t)=\sum_{r,k}h_{kr}(t)\,(t-t_*)^{r}\ln^{k-1} (t-t_*)
\end{equation}
with finitely many terms, where all exponents $r$ are rational numbers, the
powers $k$ take only finitely many natural values, and all functions
$h_{kr}$ are real analytic at the point $t=t_*$.
\end{Lem}

This Lemma is explained in \secref{sec:topology}. As soon as the
representation \eqref{asymp-expan} is established, the rest of the proof
is as described in Example~\ref{ex:dulac-series}. It remains only to notice
that the convergence of the expansion \eqref{asymp-expan} is crucial: if
all terms in this expansion are zeros, then the integral itself is
identically zero and hence has no isolated roots at all (compare with the
Remark in \secref{sec:dulac}).

\subsection{Existential finiteness}
Abelian integrals depend on the choice of the Hamiltonian and the form; in
order to make the parameter space finite-dimensional, one has to restrict
their degrees. Then the coefficients of the Hamiltonian $H$ and the form
$\omega$ become the natural parameters of the problem.

The parameter space is intrinsically compact: indeed, replacing $H$ and
$\omega$ by $c_1H$ and $c_2\omega$ with $c_1,c_2$ nonzero constants,
clearly does not affect the number of isolated roots of the corresponding
integral. Hence the parameter spaces are in fact projective spaces of
appropriate dimensions.

Fix any two integer numbers $n$ and $d$ and consider all Hamiltonians of
degree $\le n$ and all polynomial forms of degree $\le d$. The existence of
\emph{uniform} upper bounds for the number of isolated zeros of the
Abelian integrals subject to the above restrictions on the degrees, was
proved by Varchenko \cite{varchenko:finiteness} and Khovanskii
\cite{asik:finiteness}.

\begin{Thm}[A.~Varchenko--A.~Khovanskii, 1984]\label{thm:varchenko-khovanskii}
For any $n,d<+\infty$, the number of ovals $\gamma\subset\{H=\const\}$
yielding isolated zeros to the Abelian integral $\oint_\gamma\omega$, is
bounded by a constant $N=N(n,d)$ uniformly over all Hamiltonians of degree
$\le n$ and all polynomial 1-forms of degree $\le d$.
\end{Thm}

One key tool in their proof is again the same Lemma~\ref{lem:AI-expansion}
properly generalized to cover the parametric construction. The second
principal ingredient is the \emph{Pfaffian elimination} technique
\cite{asik:fewnomials}. The latter is a method allowing to reduce the
question on zeros of functions involving real powers, logarithms,
arctangents and other functions that can be defined using Pfaffian
equations with algebraic right hand sides, to the question on zeros of
certain auxiliary systems of equations involving only algebraic functions.
Application of the generalized Lemma~\ref{lem:AI-expansion} allows to
derive the existential finiteness theorem for Abelian integrals, using in
the standard way the compactness and localization arguments, from the
similar existential finiteness assertion for converging multivariate
expressions of the form \eqref{asymp-expan}. The latter assertion can in
turn be derived, using the Pfaffian elimination technique, from the
existential finiteness theorem for analytic families of functions as
introduced in \secref{sec:analytic-functions}. This last case was already
settled by Theorem~\ref{thm:fcyc-anal}, and this completes the proof.

However, since both Theorem~\ref{thm:fcyc-anal} and its corollary,
Theorem~\ref{thm:finite-cyclicity} do not give any information on the
\emph{number} of isolated roots of functions,
Theorem~\ref{thm:varchenko-khovanskii} inherits this quality as a purely
existential statement.

\begin{table}[ht]\label{table:AI}
\begin{tabular}{|c||c|c||c|}
\hline
 {}&\multicolumn{3}{c}
 {\setheight\emph{Degree of the form}}\vline
 \\
 \cline{2-4}
 {\bf Hamiltonian $H$}  &
 \twolines{\bf Low degree}{($\deg \omega\le 2$)}&
 \threelines{\bf Arbitrary}{\bf $d=\deg \omega$}{\bf (asymptotic)} &
 \threelines{\bf Arbitrary}{\bf degree}{\bf (constructive)}
 \\
\hline
 \twolines{\bf Quadratic:}{$H=x^2+y^2$} &
 \multicolumn{3}{c}{
 {Integrals are polynomial} {functions of $t$}
 }\vline
 \\
 \hline
 \twolines{\bf Elliptic}{$H=y^2+x^3-3x$} &
 \twolines{$\le 1$ root}{Petrov \cite{petrov:cubic-real}} &
  &
 \threelines{nonoscillation}{theorem}{Petrov \cite{petrov:cubic-nonosc}}
 \\
 \hline
 \twolines{\bf General}{\bf cubic $H$} &
 \threelines{numerous results}{in particular,}
   {Gavrilov \cite{gavrilov:quadratic}}&
  &
 \threelines{$5(d+2)$}{Horozov--Iliev}{\cite{horozov-iliev:linear}}
 \\
 \hline
 \twolines{\bf Hyperelliptic}{$H=y^2+P_n(x)$}&
 ---&
 Petrov \cite{petrov:trick} &
 \threelines{tower function,}{Novikov--Yako-}{venko \cite{era-99}}
 \\
 \hline
 \threelines{\bf General}{$n$th \bf degree}{\bf polynomial} &
 ---&
 \threelines{$C(n)d+O(1)$,}{Khovanskii}{\cite{asik:unpub}} &
 {\bf ???}
 \\
 \hline
 \end{tabular}

 \bigskip
 \caption{Constructive  tangential Hilbert 16th problem:
 partial synopsis}
\begin{minipage}{0.9\hsize}
\begin{small}
Dashes in the cells mark irrelevant (artificial) problems. Blank cells
indicate that the corresponding problem possesses no specific (more
accurate) solution other than the one implied by stronger versions
appearing to the right or down in the Table.\par
\end{small}
\end{minipage}
\vss
 \end{table}

\subsection{Constructive finiteness}
The simplest (trivial) case when constructive bounds for the tangential
Hilbert problem can be easily produced, is the case of quadratic
Hamiltonians. In this case all Abelian integrals can be explicitly
computed, being actually \emph{polynomial} functions of $t$.

The constructive finite cyclicity problem for Abelian integrals is also a
relatively simple assertion: the maximal multiplicity of an isolated zero
of an Abelian integral admits an upper bound in terms of $\deg \omega$ and
$\deg H$ found by P.~Marde\v si\'c in \cite{mardesic:multiplicity}.

Historically the first nontrivial problem concerning zeros of Abelian
integrals appeared in connection with what later became known as the
Takens--Bogdanov bifurcation. In order to prove that no more than one
limit cycle appears by deformation of the generic cuspidal singular point
on the plane, one has to verify that for any 1-form $\omega=(\alpha+\beta
x)y\,dx$ its integral over the closed ovals of the Hamiltonian
$H(x,y)=y^2+x^3-3x$ has no more than one isolated zero. R.~Bogdanov in
\cite{bogdanov:integral} proved this claim, achieving one of the first
results on Abelian integrals. Later his proof was considerably simplified
by Ilyashenko in the paper \cite{ilyashenko:elliptic} which introduced some
very important tools, among them the idea of complexification of the
Abelian integrals.

There are numerous studies treating other low-degree cases that mostly
appeared in connection with bifurcations of polynomial vector fields of
low degrees. Starting from the paper by Bogdanov, this direction was
pursued, among other, by
     F. Dumortier,
     A. Gasull,
     L. Gavrilov,
     F. Girard,
     E. Horozov,
     I. Iliev,
    Yu.~Ilyashenko,
     A. Jebrane,
     B. Li,
     C. Li,
     J. Llibre,
     P. Marde\v si\'c,
     G. Petrov,
     R. Roussarie,
     C. Rousseau,
     Z. Zhang,
     Y. Zhao,
     H. Zoladek,
to mention only some names and the most recent works. Among these accurate
bounds, the following result is remarkable by its succinct formulation and
difficult proof.

\begin{Thm}[L. Gavrilov--E.Horozov--I.Iliev \cite{gavrilov:quadratic}]
\label{thm:gavrilov-quad} For any cubic Hamiltonian $H$ with four distinct
critical values, and any quadratic 1-form $\omega$, the corresponding
Abelian integral has no more than two isolated roots.

Moreover, in the perturbation \eqref{pert-ham} no more than 2 limit cycles
may appear, including those born from separatrix polygons.
\end{Thm}

This settles the localized version of the Hilbert 16th problem for
quadratic vector fields arbitrary close to Hamiltonian quadratic vector
fields with the specified Hamiltonians (cubic with 4 distinct critical
values). The bound is accurate.

Another remarkable result is due to G.~Petrov \cite{petrov:cubic-nonosc}
who studied completely the \emph{elliptic case} $H(x,y)=y^2+P_3(x)$, where
$P_3\in\R[x]$ is a cubic univariate polynomial: by affine transformations
one can always reduce $P_3$ to the form $P_3(x)=x^3-3x$. Petrov proved
that the Abelian integrals of forms of arbitrary degree form a
\emph{non-oscillating} (in other languages \emph{disconjugate}, or
\emph{Chebyshev}) family: \emph{the number of isolated roots never becomes
equal or exceeds the dimension of the linear space of all such integrals}.

\begin{Rem}
One can easily verify that any linear space spanned by $d$ linear
independent analytic functions $f_1(t),\dots,f_d(t)$ contains a nontrivial
linear combination $c_1 f_1+\cdots+c_d f_d$ exhibiting a root of
multiplicity $d-1$ at any preassigned point.
\end{Rem}

Methods introduced by Petrov were further elaborated and refined.  The
following bound obtained by E.~Horozov and I.~Iliev, though not sharp,
covers the case of Abelian integrals of arbitrary polynomial forms over
level curves of any cubic Hamiltonians.

\begin{Thm}[\cite{horozov-iliev:linear}]\label{thm:horozov-iliev}
For any cubic Hamiltonian $H(x,y)$ the Abelian integral of a form of
degree $d$ cannot have more than $5(d+2)$ isolated zeros.
\end{Thm}

By similar methods the \emph{quartic} Hamiltonians (of degree 4) with
\emph{elliptic} level curves were studied in \cite{girard-jebrane}. The
case of quartic elliptic Hamiltonians perturbed within the Li\'enard
equation was a subject of four recent preprints of C.~Li and F.~Dumortier.

\subsection{Asymptotic bounds}
Results of the different type, asymptotic bounds, take advantage of the
asymmetry of the roles played by 1-forms and Hamiltonians. Unlike the
original Hilbert 16th problem, in which the coefficients of the
polynomials $P,Q$ all enjoy equal rights as parameters, the roles of the
Hamiltonian $H$ and the form $\omega$ in the definition of the Abelian
integrals are fairly different (e.g., the integral depends on $\omega$
linearly, whereas even small variations of $H$ can result in drastic
changes of the domain of definition). Thus it makes complete sense to
separate these two parameters and study first the dependence on $\omega$,
treating $H$ as fixed (``individual'').

The corresponding ``semiconstructive'' problem was addressed in a number
of recent publications. First, very excessive (double exponential in $d$)
upper bounds were obtained in \cite{invmath-95} and almost immediately
improved to simple exponential expression in \cite{annalif-95} by
Ilyashenko, Novikov and the author. The ultimate result in this direction,
a \emph{linear} upper bound of the form $C_1(n)\cdot d+C_2(n)$, was
obtained by Khovanskii \cite{asik:unpub} using some ideas developed
earlier by G.~Petrov \cite{petrov:trick}.  Here the constant $C_1(n)$ is
absolutely explicit (e.g., does not exceed a double exponential of $n$),
while the second constant $C_2(n)$ is purely existential though uniform
over all Hamiltonians of degree $\le n$.

\subsection{Hyperelliptic integrals}
The only particular case covering Ham\-ilton\-ians of arbitrarily high
degrees, for which constructive solution of the tangential Hilbert problem
is known, is that of \emph{hyperelliptic} Hamiltonians,
\begin{equation}\label{hyperel}
  H(x,y)=y^2+P_{n+1}(x), \quad P_{n+1}=x^{n+1}+a_1
  x^{n-1}+\cdots+a_{n-1}x+a_n.
\end{equation}
Singular points of the corresponding Hamiltonian system correspond to
critical points of the polynomial $P$, called \emph{potential}. The reason
why the class of hyperelliptic polynomials is especially simple, one can
vaguely attribute to the fact that from many points of view, hyperelliptic
polynomials behave like their univariate potentials. In particular, this
concerns topology of the bundles defined by complexification of $H$.

Under the addition assumption (believed to be technical, though it occurs
independently and persistently in several related problems) that \emph{all
critical points of the potential $P$ are real}, the tangential Hilbert
problem turns out to be \emph{constructively solvable}. More precisely, as
shown by D.~Novikov and the author \cite{era-99}, there exists an algorithm
defining an elementary function $C(n,d)$ of two integer arguments $n$ and
$d$, such that for any form of degree $\le d$ and any hyperelliptic
Hamiltonian $n$ of degree $\le n$ having only real critical points, the
corresponding hyperelliptic Abelian integral has no more than $C(n,d)$
isolated roots.

The algorithm involves several nested inductive constructions, resulting
in an \emph{extremely excessive} bound: it is given by a \emph{tower
function} (an iterated exponential) of height greater than 5 but probably
smaller than 10.

The proof is based on the fact that Abelian integrals satisfy a system of
first order linear ordinary differential equations with rational
coefficients, called the \emph{Picard--Fuchs} system.

\subsection{Quantitative theory of ordinary differential equations as a
tool for constructive solution of the tangential Hilbert problem} The
method used for constructive solution of the tangential Hilbert problem in
the hyperelliptic case, is fairly general. In \cite{fields} one can find an
introduction to the general theory allowing to investigate zeros of
functions defined by ordinary differential equations with polynomial and
rational coefficients. Basics of this theory were developed in a series of
joint works by D.~Novikov and the author; they are briefly recalled below.

The principal goal of these lecture notes is twofold. First, we show that
the Picard--Fuchs system of differential equations can be
\emph{explicitly} written down for an arbitrary generic Hamiltonian (not
necessarily a hyperelliptic one), at least at the price of certain
\emph{redundancy}. Some additional information may be extracted from the
explicit derivation procedure. In particular, we show that for almost any
Hamiltonian $H$ (chosen for convenience of degree $n+1$):
\begin{enumerate}
\item there exist $\mu$ different Abelian integrals $I_k=\oint\omega_k$ of
monomial 1-forms $\omega_1,\dots,\omega_\mu$, $\mu=n^2$, satisfying
together a \emph{Fuchsian system} of linear ordinary differential
equations with rational coefficients of the form
\begin{equation}\label{fuchsian}
 \begin{gathered}
  \frac{dx}{dt}=A(t)x,\qquad A(t)=\sum_{j=1}^\mu\frac{A_j}{t-t_j},
  \\
  x\in\C^\mu,\quad A_j\in\Mat_{\mu\times\mu}(\C),\quad t\in\C,
  t_1,\dots,t_\mu\in\C,
 \end{gathered}
\end{equation}
such that
\item integral $I(t)$ of any other form $\omega$ can be represented as a linear
combination of the integrals $I_1,\dots,I_\mu$,
\begin{equation}\label{envelope-1}
  I(t)=q_1(t)I_1(t)+\cdots+q_\mu(t)I_\mu(t),
\end{equation}
with polynomial coefficients $q_j\in\C[t]$, $\deg
  q_j\le \deg\omega/\deg H$.
\end{enumerate}
The residue matrices $A_j$ can be sufficiently completely described, in
particular, upper bounds on their norms can be placed.

The second goal is to introduce several results on the number of zeros of
functions defined by systems of ordinary linear equations with polynomial
and rational coefficients, and polynomial combinations of such functions.
We start with the simplest case of one $n$th order linear equation and show
that isolated zeros of its solutions can be described in terms of the
magnitude of coefficients of this equation. In the simplest real case this
is a classical theorem by de la Vall\'ee Poussin \cite{poussin}, which we
generalize for the complex analytic context.  Afterwards we study systems
of polynomial ODE's in the real and complex space of arbitrary dimension;
here for the first time appear explicit bounds given by tower functions.
Finally, we discuss the case of Fuchsian systems of the form
\eqref{fuchsian} and show that under certain natural restrictions on the
monodromy group of this equation, a global upper bound on the number of
zeros can be given in terms of norms of the residue matrices $A_j$.

One can hope that combination of these two techniques ultimately would
allow to construct an explicit bound for the tangential Hilbert problem in
the general case, filling the right bottom corner of the Table~2.

\chapter{Abelian integrals and differential equations}

In this lecture we explain the connection between Abelian integrals and
linear ordinary differential equations.

\section{Complexification of Abelian integrals: topological
approach}\label{sec:topology}

\subsection{General scheme}
We recall here the basic construction of complexification of Abelian
integrals. All details can be found in various textbooks, \cite{odo-2}
being the principal source.

A bivariate polynomial $H\in\C[x,y]$ defines a map from $\C^2$ to $\C^1$
with preimages of points being affine algebraic curves. It turns out that
the map $H$ is a \emph{topological bundle} over a complement to finitely
many points $\S=\{t_1,\dots,t_r\}\subset\C$. This allows to identify in a
canonical way the homology groups of all fibers $H^{-1}(t)$ sufficiently
close to $H^{-1}(t_*)$, which in turn allows to extend integrals over
1-cycles of polynomial 1-forms as complex analytic functions of $t$ in $U$,
ramified over the singular locus $\S$. Geometrically this can be expressed
as introducing a locally flat connexion on the (co)homological bundles over
the punctured sphere.

From the same topological arguments it follows that for any polynomial
1-form $\omega$ the linear space generated by integrals
$\oint_{\delta_i}\omega$ over any family $\delta_1(t),\dots,\delta_\mu(t)$
of cycles forming a basis in the first homology group of the fiber
$X_t=H^{-1}(t)$, is invariant by analytic continuation along the loops
avoiding the locus $\S$. The monodromy group consisting of all
automorphisms of this space occurring as the result of continuation over
all loops, is independent of the choice of the form $\omega$ and depends
only on the Hamiltonian $H$.

These fairly general topological considerations already imply that the
Abelian integrals form a finitely generated module over the ring of
polynomial functions of $t$, and generators of this module satisfy a
system of first order linear ordinary differential equations with rational
(in $t$) coefficients, called \emph{Picard--Fuchs} equations.

However, in order to apply methods described in subsequent sections, it is
necessary to obtain \emph{quantitative} characteristics of these equations,
in particular their dimension, degree and the magnitude of coefficients.
Part of this information (e.g., the dimension) can be easily achieved from
the above construction. To obtain upper bounds on the degrees, some other
more refined considerations are required, but it is practically impossible
to derive bounds for the coefficients using only topological arguments as
above.

The current section contains a brief exposition of well-known facts
leading to derivation of Picard--Fuchs equations and representation of the
space of Abelian integrals as the Picard--Vessiot extension. In the next
section we suggest an alternative approach based on elementary algebraic
consideration, that allows to derive explicitly the Picard--Fuchs system at
the price of certain redundancy.

\subsection{Topological bundles defined by proper maps}
Let $f\:M\to N$ be a smooth map between two manifolds.  Recall that a
point $b\in N$ is a regular value for $f$, if the rank of the differential
$f_*\:T_x M\to T_{b}N$ is maximal (equal to $\dim N$) at all points $x\in
X_{b}$ of the preimage. Complement to the set of regular values consists
of \emph{critical values} and is denoted $\crit f$.

\begin{Lem}\label{lem:morse}
If $b$ is a regular value of a proper map $f\:M\to N$, then there exists a
neighborhood $U\owns b$ such that all preimages $X_y=f^{-1}(y)\subset M$
are diffeomorphic to $X_b=f^{-1}(b)$ so that $f^{-1}(U)\simeq X_b\times U$.
\end{Lem}

\begin{proof}
Consider an arbitrary vector $v_0\in T_{b}N$ and embed it into a vector
field $v$ on $N$. We claim that in a sufficiently small neighborhood of the
preimage $X_{b}$ one can construct a smooth vector field $w$ such that
$f_*w=v$, that is, $w$ and $v$ are $f$-related. Such a field obviously
exists near each point $x\in X_{b}$, since $f_*$ is surjective (and takes
the form of a parallel projection in suitably chosen local coordinates by
virtue of the theorem on rank). Now it remains to choose a finite covering
of $X_{b}$ by these neighborhoods and patch together the corresponding
vector fields, using the appropriate partition of unity.

To conclude the proof, notice that the (local) flows of $v$ and $w$ are
conjugate by $f$ (by construction), hence the local flow of $w$, defined
in some neighborhood of $X_{b}$, takes the latter preimage into preimage
of the corresponding point $y$ on the flow curve of $b$. Since the initial
vector $v_0$ can be chosen pointing to any direction, this proves that all
sufficiently close preimages are diffeomorphic to each other.
\end{proof}

\begin{Rem}
The diffeomorphism between close preimages is not canonically defined, but
any two such diffeomorphisms are homotopically equivalent, since the
vector field $w$ from the proof of the Lemma is homotopically unique.
\end{Rem}

\begin{Cor}
If  $M$ is compact, then $f$ is a topological bundle over $N\ssm\crit
f$.\qed
\end{Cor}

\begin{Rem}
Assertion of this lemma for the case $N=\R$ is the fundamental principle of
the Morse theory.
\end{Rem}

\subsection{Homology bundle and flat connexion on it}
A locally trivial topological bundle $f\:M\to N$ defines in a canonical
way the associated homology bundle over $N$ with the fibers being the
homology groups $H_i(X_y,\mathbb Z)$ (we will be only interested in the
case $i=1$) together with a flat connexion on this bundle. This means that
any 1-cycle on any particular fiber $X_a$ can be transported to any other
fiber $X_b$ along any path $\gamma$ connecting $a$ and $b$ in $N$. Flatness
means that continuation along any sufficiently small loop beginning and
ending at $a$, returns any 1-cycle to its initial position. However,
transport along ``long'' loops (not contractible in $N$) can result in a
nontrivial transformation of the homology, called the \emph{monodromy
transformation}.

Fix a point $a\in N$ and choose a basis $\delta_1,\dots,\delta_\mu$ of
1-cycles in the group $H_1(X_a,\mathbb Z)$, where $\mu$ is the rank of the
homology group. Then the monodromy transformation corresponding to any loop
$\gamma$ from the fundamental group $\pi_1(N,a)$ can be described by the
corresponding \emph{monodromy matrix}
$M=M_\gamma\in\Mat_{\mu\times\mu}(\mathbb Z)$ with integral entries
$m_{ij}$:
\begin{equation}\label{monodromy-cycles}
  \Delta_\gamma\delta_j=\sum_{i=1}^\mu \delta_i m_{ij},\qquad
  i=1,\dots,\mu.
\end{equation}
In the matrix form the monodromy transformation acts on the row vector
$\bdelta=(\delta_1,\dots,\delta_\mu)$ as multiplication by the matrix
$M_\gamma$ from the right, $\Delta_\gamma \bdelta=\bdelta\cdot M_\gamma$.

The integer-valued matrix $M_\gamma$ is obviously invertible. Moreover,
the correspondence $\gamma\mapsto M_\gamma$ is an (anti)representation of
the fundamental group $\pi_1(N, a)$ in $\mathrm{GL}(n,\mathbb Z)$. This
implies, in particular, that $\det M_\gamma=\pm1$ for all loops $\gamma$.

The flat connexion on the homology bundle defines a connexion on the
cohomology bundle which makes it possible to compute the covariant
derivative of 1-forms: any smooth 1-form $\omega$ on $M$ can be restricted
on any fiber $X_a$ and integrated along a continuous (horizontal, locally
constant) family of 1-cycles $\delta(y)$, resulting in a function $I(y)$.
The result of continuation of this function along closed paths in $N$ is
completely determined by the monodromy group and does not depend on the
choice of the form $\omega$. More precisely, integrating the form $\omega$
over each of the 1-cycles $\bdelta(t)=(\delta_1(y),\dots,\delta_\mu(y))$,
we obtain a tuple of continuous functions
$I_j(y)=\oint_{\delta_j(y)}\omega$, $j=1,\dots,\mu$, which after
continuation along a path $\gamma\in\pi_1(N,a)$ undergo the transformation
\begin{equation}\label{monodromy-integrals}
  \Delta_\gamma \boldsymbol I=\boldsymbol I\cdot M_\gamma,\qquad
  \boldsymbol I=(I_1,\dots,I_\mu),\ I_j=I_j(y)
\end{equation}
with the same matrices $M_\gamma$ independently of the form $\omega$.  This
basic fact lies in the core of the topological theory outlined below.

\subsection{Topological bundles defined by polynomial maps} Our goal is
to apply the previous construction to the polynomial map $H\:\C^2\to\C$
considered as a smooth map between smooth manifolds. Since compactness of
the preimages is crucial for these arguments, we need to compactify the
domain (and the range) of $H$.

In contrast with the one-dimensional case, it is in general impossible to
extend $H$ as a map between the natural compactifications $\C P^2$ and $\C
P^1$ respectively, since on the infinite line $\C P^1_\infty\subset\C P^2$
one has several \emph{points of indeterminacy}: they occur at the
intersections between compactified preimages $H^{-1}(t)$ and the infinite
line. Algebraically this can be seen as the indeterminacy of the rational
expression $H(1/z,y/z)=P(y,z)/z^{d}$, $d=\deg H$, at the points where
$P(y,0)=0$ (at all other points of the infinite line $\{z=0\}$ one can
assign the value $H=\infty$ to this ratio).

The problem can be resolved by blowing up these indeterminacy points, in
the same way as blowing up the origin allows to assign values from $\C
P^1=\C\cup\{\infty\}$ to the rational expression $R(x,y)=y/x$ that is
initially indeterminate. After a series of blow-ups at indeterminacy
points one arrives at a \emph{compact} two-dimensional complex manifold
$M$ and a map (still denoted by $H$) from $M$ to $\C P^1$, called
\emph{determination} of the initial polynomial map.

Now one can apply Lemma~\ref{lem:morse}. Since the determination
$H\:M\to\C P^1$ is algebraic, it has only a finite number of critical
values $\S_H=\{t_1,\dots,t_s\}\subset\C P^1$ and we conclude that it
defines a topological bundle over the complement of these \emph{exceptional
values}. Any 1-cycle $\delta(a)\subset H^{-1}(a)$ on a fiber of this
bundle can be embedded into a continuous horizontal family $\delta(t)$ of
1-cycles, ramified over the exceptional locus. A polynomial 1-form
$\omega$ extends as a meromorphic 1-form on $\C P^2$ with the polar
divisor $\C P^1_\infty$ that after the blowing up becomes an algebraic
hypersurface $D$ in $M$. Let $\S$ be the union of $\S_H$ and the critical
values of the projection $H$ restricted on $D$ (including the images of
the non-smooth points of $D$).

The pullback of $\omega$ on $M$ can be integrated along the family
$\delta(t)$: by the Cauchy--Stokes theorem, the result depends only on the
homology class of the cycle. It can be easily shown that the result of this
integration is an analytic function of $t$, ramified over $\S$. Its
monodromy (the result of analytic continuation along closed loops avoiding
the exceptional locus $\S$) is as before determined only by $H$.

This construction proves the following result.

\begin{Prop}[cf.~with \cite{odo-2}]
For any polynomial $H\:\C^2\to\C$, any 1-cycle $\delta$ on a nonsingular
level curve $X_{*}=\{H=t_*\}\subset\C^2$ and any polynomial 1-form
$\omega$ the Abelian integral $I(t)=\oint_\delta\omega$ can be extended as
an analytic multivalued function ramified over a finite number of points
depending only on $H$.\qed
\end{Prop}

Behavior of the integrals near the ramification locus is relatively tame:
it can be seen that any integral can grow no faster than polynomially in
$|t-t_j|^{-1}$ as $t$ tends to some $t_j\in\S$ remaining in any sector
with the vertex at $t_j$.

\subsection{Picard--Fuchs system}
Let as before $\bdelta(t)=(\delta_1(t),\dots,\delta_\mu(t))$ be a
continuous family of 1-cycles forming a basis (frame) of the first homology
group of the respective fibers $H^{-1}(t)$, arranged as a row vector. One
can show that there exist $\mu$ polynomial 1-forms
$\omega_1,\dots,\omega_\mu$ such that the \emph{period matrix}
\begin{equation}\label{period-matrix}
  X(t)=\begin{pmatrix}
    \oint_{\delta_1}\omega_1 & \cdots & \oint_{\delta_\mu}\omega_1 \\
    \vdots & \ddots & \vdots \\
    \oint_{\delta_1}\omega_\mu  &\cdots &\oint_{\delta_\mu}\omega_\mu
  \end{pmatrix}
\end{equation}
is not identically degenerate, $\det X(t)\not\equiv0$. From
\eqref{monodromy-integrals} it follows that
\begin{equation}\label{mono-X}
  \Delta_\gamma X(t)=X(t)M_\gamma,\qquad\forall \gamma\in\pi_1(\C
  P^1\ssm\S,a).
\end{equation}

Differentiating the identity \eqref{mono-X}, we see that the derivative
$\dot X(t)$ has the same monodromy (i.e., $\dot X$ is multiplied by the
same matrix factors $M_\gamma$). Therefore the ``logarithmic derivative''
$A(t)=\dot X(t)\cdot X^{-1}(t)$ is single-valued (invariant by all
monodromy transformations) meromorphic matrix function having poles of
finite order at the points of $\S$ and eventually at the points of
degeneracy of $X(\cdot)$:
\begin{equation*}
  \Delta_\gamma A(t)=
  \dot X(t)M_\gamma\cdot M^{-1}_\gamma X^{-1}(t)=\dot X(t)\cdot X^{-1}(t)
  =A(t)
\end{equation*}
for any loop $\gamma\in\pi_1(\C\ssm\S, t_*)$. From this we conclude that
$A(t)$ is a rational matrix function while the period matrix $X(t)$ is a
fundamental matrix solution to the system of linear ordinary differential
equations with rational coefficients on $\C P^1$,
\begin{equation}\label{pic-fuchs}
  \dot X=A(t)X,\qquad A(\cdot)\in\Mat_{\mu\times \mu}(\C(t)).
\end{equation}
The common name for various such systems satisfied by Abelian integrals,
is the \emph{Picard--Fuchs system} (or Picard--Fuchs equation).

Integrals of any other form $\omega$ can be expressed as linear
combinations of integrals of the framing forms $\omega_i$, $i=1,\dots,\mu$,
with coefficients from the field $\C(t)$ of rational functions. Indeed,
multiplying the row vector $\boldsymbol I(t)=(I_1(t),\dots,I_\mu(t))$,
$I_j=\oint_{\delta_j}\omega$, by $X^{-1}(t)$, we obtain (for the same
reasons as above) a single-valued hence rational row-vector function
$\boldsymbol r(t)=(r_1(t),\dots,r_\mu(t))$, that is,
\begin{equation}\label{envelope}
  \oint_{\delta(t)}\omega=\sum_{j=1}^\mu r_j(t)\oint_{\delta(t)}\omega_j,
  \qquad r_j\in\C(t).
\end{equation}
for any cycle $\delta(t)$ continuously depending on $t$. Note that the
space of functions representable as \eqref{envelope}, is closed by
derivation.

We summarize this as follows. Recall that a \emph{Picard--Vessiot
extension} is a differential field of analytic multivalued functions
obtained by adjoining to the field $\C(t)$ all components of a fundamental
matrix solution of a system of linear ordinary differential equations with
rational coefficients.

\begin{Thm}\label{thm:pic-ves}
Abelian integrals belong to a Picard--Vessiot extension for some system of
linear ordinary differential equations with rational coefficients.
\end{Thm}

Later we discuss this formulation and relevant issues in more details.
However, it is important to stress here that neither entries of the
rational matrix $A(t)$ nor the rational coefficients $r_j(t)$ can be
computed explicitly from the above construction without additional
considerations.

\section{Hamiltonians transversal to infinity}\label{sec:trans-infinity}

\subsection{Definition}
For an arbitrary Hamiltonian $H$, even location of the ramification points
$t_j$ is difficult to describe without effectively resolving all the
indeterminacy points at infinity. However, under some natural and generic
assumptions one may guarantee that no new critical points will appear after
compactification and blowing up the indeterminacy points on the infinite
line $\C P^1_\infty\subset\C P^2$. This will immediately imply that the
set $\S$ of exceptional values must be a subset of $\crit H$, the set of
critical values corresponding to \emph{finite} critical points in $\C^2$
only.

\begin{Def}
A polynomial $H\in\C[x,y]$ is said to be \emph{transversal to infinity}, if
its principal homogeneous part $L=\sum_{i+j=d} a_{ij}x^iy^j$, $d=\deg H$,
factors as a product of pairwise different linear forms.

Equivalent conditions follow.
\begin{enumerate}
 \item The principal homogeneous part has an isolated critical point at the
 origin;
 \item The partial derivatives $\pd Lx$, $\pd Ly$ are mutually prime;
 \item $H$ has exactly $\mu=(\deg H-1)^2$ critical points in $\C^2$ if
 counted with multiplicities;
 \item Each level curve $\{H=t\}$ intersects transversely the infinite
 line $\C P^1_\infty\subset\C P^2$ after projective compactification.
\end{enumerate}
\end{Def}

\subsection{Topology of polynomials transversal to infinity}
We prove now that for $H$ transversal to infinity, $\S=\crit H$.
\begin{Prop}
A polynomial $H\:\C^2\to\C^1$ transversal to infinity is a topological
bundle over the set $\crit H$ of critical values of $H$.
\end{Prop}

Instead of proving this by resolving the indeterminacy points at infinity,
one may modify the proof of Lemma~\ref{lem:morse} sketched above, and
construct the vector field $w$ near an arbitrary infinite point $p\in\C
P^1_\infty$ on $\C P^2$ with the property $H_* w=\pd{}t$.

\begin{proof}
It will be shown below in \secref{sec:division-by-gradient} that if $H$ of
degree $n+1$ is transversal to infinity, then there exist two polynomials
$a,b\in\C[x,y]$ of degree $n-1$ such that
\begin{equation*}
  a\pd Hx+b\pd Hy=x^{2n-1}+\cdots,
\end{equation*}
where the dots stand for a bivariate polynomial of degree $\le 2n-2$. One
can easily check that the rational vector field
$$
 w=\frac {a(x,y)}{x^{2n-1}+\cdots}\pd {}x+
 \frac{b(x,y)}{x^{2n-1}+\cdots}\pd {}y,
$$
$H$-related to the field $\pd {}t$, in the chart $(1/x,y/x)$ is regular
(smooth) \emph{on the infinite line} $1/x=0$ (more precisely, on its
affine part covered by this chart). The other affine part is covered by the
field that is obtained in a similar way from solution of the equation $a\pd
Hx+b\pd Hy=y^{2n-1}+\cdots$.

Thus near each point of the compactified level curve $\overline
X_a\subset\C P^2$ one has a smooth vector field $H$-related to $\pd {}t$
in the finite part $\C^2\subset\C P^2$ (in particular, this implies that
this field vanishes at all points of indeterminacy of $H$ on $\C
P^1_\infty$). The rest of the proof is the same as in
Lemma~\ref{lem:morse}.
\end{proof}

\subsection{Module of Abelian integrals, Gavrilov and Novikov theorems}
For Hamiltonians transversal to infinity, the constructions of
\secref{sec:topology} can be further refined. In particular, the choice of
the forms $\omega_1,\dots,\omega_\mu$ can be made much more explicit.

\begin{Lem}[L. Gavrilov \cite{gavrilov:petrov-modules}]\label{lem:gavrilov}
Let $H$ be a Hamiltonian of degree $n+1$ transversal to infinity, with
distinct critical values $t_1,\dots,t_\mu$, $\mu=n^2$. Then one can choose
$n^2$ monomial 1-forms $\omega_1,\dots,\omega_\mu$ of degrees $\le 2n$ so
that the respective period matrix has the determinant $\det
X(t)=c(t-t_1)\cdots(t-t_\mu)$ with $c\ne 0$.
\end{Lem}

In this assertion we use the convention on degrees of polynomial $k$-forms
formulated in \secref{sec:division-by-gradient}. As a corollary, one can
derive the following result that refines the assertion of
Theorem~\ref{thm:pic-ves}.

\begin{Cor}\label{cor:gavrilov-rep}
Abelian integral of a 1-form $\omega$ of degree $d$ can be represented as
\begin{equation}\label{gavrilov-rep}
\begin{gathered}
  \oint_{\delta(t)}\omega=\sum_{j=1}^\mu p_j(t)\oint_{\delta(t)}\omega_j,
  \\
  p_j(t)\in\C[t],\qquad \deg \omega_j+\deg H\cdot\deg p_j\le\deg\omega.
\end{gathered}
\end{equation}
\end{Cor}
In other words, Abelian integrals constitute a module over the ring
$\C[t]$ that is generated by integrals of the basic forms $\omega_j$.

The constant $c$ from Lemma~\ref{lem:gavrilov} depends on the choice of the
monomial forms and the Hamiltonian. Its value was explicitly computed by
A.~Glutsuk \cite{glutsuk} following some ideas of Yu.~Ilyashenko, and a
simple elementary proof of the inequality $c\ne0$ for an appropriate
choice of the monomial forms was obtained by D.~Novikov
\cite{mit:irredundant}. In the same paper \cite{mit:irredundant} it is
proved, using some of the methods described below, that the period matrix
satisfies a system of linear ordinary differential equations
\begin{equation}\label{irredundant}
  \dot X=\frac1{(t-t_1)\cdots(t-t_\mu)}P(t)X,\qquad P(t)=\sum_{j=0}^\mu
  t^jP_j,
\end{equation}
with a matrix polynomial $P(t)$ of degree $\mu$, in general having
Fuchsian singularities at all points $t_j$ of the ramification locus $\S$,
but a \emph{non-}Fuchsian singularity at $t=\infty$.

As yet another corollary, one can derive a Picard--Fuchs system for the
period matrix. Let $\omega_j$ be as in Corollary~\ref{cor:gavrilov-rep}.
Consider the closed $2$-forms $dH\land \omega_j$ and let $\Omega_j$ be any
polynomial primitives satisfying the conditions
\begin{equation*}
  d\Omega_j=dH\land\omega_j,\qquad j=1,\dots,\mu.
\end{equation*}
Each $\Omega_i$ can be expanded as in Corollary~\ref{cor:gavrilov-rep},
yielding an identity between the integrals,
\begin{equation*}
  \oint\Omega_i=\sum_{j=1}^\mu p_{ij}(t)\oint\omega_j
\end{equation*}
valid for any choice of a continuous family of cycles of integration
$\delta(t)$. On the other hand, as will be shown in
\secref{sec:gelfand-leray}, the derivative of each integral
$\oint\Omega_i$ is exactly the integral $\oint\omega_i$. Differentiating
the above identities, we arrive to the matrix differential equation
\begin{equation*}
  X=\dot PX+P\dot X,\qquad
  P=P(t)=\|p_{ij}(t)\|\in\Mat_{\mu\times\mu}(\C[t])
\end{equation*}
with a matrix polynomial $P(t)$ of some known degree.  If required, this
identity can be resolved to the form \eqref{pic-fuchs}.

\subsection{Commentaries}
The proof of Lemma~\ref{lem:gavrilov} and Corollary~\ref{cor:gavrilov-rep}
is based on a finer than before analysis of topology of the bundle $H$ for
polynomials transversal to infinity. In particular, in assumptions of the
Lemma, one can choose a special framing of the homology bundle by
\emph{vanishing cycles} $\delta_j(t)$, represented by loops on the
preimage $X_t$ that shrinks to a point when $t\to t_j$ (a special
precaution is required to avoid problems with multivaluedness). For such
choice of the cycles, the period matrix $X(t)$ must have a vanishing
column at each of the points $t_j$, for any collection of the framing
forms. Next, in this case the determinant $\det X$ is a single-valued
function that therefore must be a polynomial (being locally bounded
everywhere on $\C$). Its growth as $t\to\infty$ depends on the degrees of
the framing forms, since the ``size'' of the cycles $\delta_j(t)$ grows in
a known way (depending only on $H$).

All this implies that the determinant of any period matrix $\det X(t)$ is
a polynomial in $t$ divisible by $\prod_{j=1}^\mu (t-t_j)$. The
coefficients $p_j(t)$ of representation \eqref{gavrilov-rep} can be found
by solving a system of linear algebraic equations. When solved by the
Cramer rule, this system yields $p_j(t)$ as a ratio of two determinants,
the determinant of the period matrix for the basic forms
$\omega_1,\dots,\omega_\mu$ in the denominator, and that of a similar
period matrix for the collection of 1-forms with $\omega_j$ replaced by
$\omega$ in the numerator. By Lemma~\ref{lem:gavrilov}, these ratios are
polynomials in $t$ and their degrees can be easily majorized in terms of
$d=\deg \omega$.

\subsection{Multiplicity of the roots of Abelian integrals: constructive
finite cyclicity}\label{sec:mardesic}
 A similar construction allows to majorize the order of a zero of any Abelian
integral, at least under the assumption that $H$ is transversal to
infinity and has only Morse critical points \cite{mardesic:multiplicity}.
Instead of the period matrix $X$, consider the matrix function $J(t)$ whose
entry $J_{ij}(t)$ is the $(i-1)$-st derivative of the integral
$I_j(t)=\oint_{\delta_j(t)}\omega$.

For the same reasons as before, its determinant $w(t)=\det J(t)$, the
Wronskian of the integrals $I_1,\dots,I_\mu$, is a single-valued hence
rational function of $t$. Its poles may occur only at the points $t_j$ and
$t=\infty$. The assumption on finite singular points implies that all
integrals $I_j(t)$ have at worst logarithmic growth near each $t_j$, and
this growth rate allows for differentiation so that $I_j^{(i-1)}(t)$ grows
no faster than $|t-t_j|^{1-j}$ as $t\to t_j$ without spiraling. These
estimates imply an upper bound on the total order of all poles of $w$ at
all finite points. The growth rate of $w$ as $t\to\infty$ depends on $\deg
\omega$ and can be easily estimated. This gives an upper bound on the
degree $\nu=\deg w(t)$ of the rational function (the total number of its
poles on $\C P^1$ including those at infinity). This degree is obviously an
upper bound for the order of any nontrivial zero of $w$. Since the order
of the Wronskian $w$ is equal to $\mu$, the number $\nu+\mu-1$ is an upper
bound for the order of any root of any integral $I_j(t)$ at any point
$t\ne t_j$.

\subsection{Reservations}
Despite more detailed constructions and more accurate considerations, the
approach based only on topological ideas cannot provide many important
data. For example, even the ``constant'' $c$ from Lemma~\ref{lem:gavrilov}
depends in a rather nontrivial way on both $H$ and the choice of the
framing forms $\omega_j$, see \cite{glutsuk} and \cite{mit:irredundant}.
Among other things, this means practical impossibility of majorizing the
polynomial coefficients $p_j$ in \eqref{gavrilov-rep}. The same refers to
the derivation of the Picard--Fuchs system: after an accurate computation,
it can be reduced to the form determined by three constant matrices
$P_1,P_2,P_3$,
\begin{equation*}
  (P_0+tP_1)\dot X=P_2 X,\qquad P_i\in\Mat_{\mu\times\mu}(\C)
\end{equation*}
but no bounds on the norms of these matrices (or their inverses) can be
obtained except for some especially simple cases.

In short, it is combinatorial parameters like degrees, dimensions, ranks
and so on that can be more or less easily derived from even the most
explicit topological constructions. In contrast, all magnitude-like
parameters (norms, absolute values, diameters of point sets etc.), require
additional arguments. Some of them can be obtained using rather elementary
algebraic considerations.

\section{Elementary derivation of the Picard--Fuchs system}

In this section we derive yet another Picard--Fuchs system by explicit
linear algebraic considerations. The advantage of this approach (besides
its transparency) is that it allows to bound explicitly the magnitude of
coefficients of the system. In addition, the Picard--Fuchs system obtained
this way possesses a nice \emph{hypergeometric} form, exhibiting only
Fuchsian singularities (though this fact was not yet fully exploited).
However, this transparency and explicitness is achieved at the price of a
certain redundancy: the dimension of the system obtained this way, is two
times bigger than the minimal possible one.

\subsection{Division by the gradient ideal}\label{sec:division-by-gradient}
Division with remainder by an ideal can be expressed in the language of
polynomial differential forms. In what follows we consider $k$-forms
$\Lambda^k$, $k=0,1,2$, with polynomial coefficients on the plane $\C^2$.
The degree of a $k$-form \emph{by definition} is the maximum of degrees of
its coefficients, plus $k$. Under such convention, $\deg
(\xi\land\eta)\le\deg\xi+\deg\eta$ for all admissible combination of ranks
of $\xi$ and $\eta$ between $0$ and $2$, and also
\begin{equation*}
  \deg d\omega\le\deg\omega\qquad\forall \omega\in\Lambda^k,
\end{equation*}
for any rank $k=0,1$. The linear space of forms of rank $k$ and degree $d$
will be denoted by $\Lambda_d^k$, and we denote $\Lambda^k_{\le
d}=\Lambda^k_0+\Lambda^k_1+\cdots+\Lambda^k_d$.

Consider $\omega=a\,dx+b\,dy\in\Lambda^1$ with only isolated singularities.
This implies that the ideal $(a,b)$ generated by the polynomials
$a,b\in\C[x,y]\simeq\Lambda^0$ has a finite codimension $\mu$ in $\C[x,y]$,
that is, there exist $\mu$  polynomials $r_1,\dots,r_\mu\in\C[x,y]$ such
that any polynomial $q$ from this ring can be represented as
\begin{equation*}
  q=av-bu+\sum_1^\mu c_i r_i,\qquad u,v\in\C[x,y],\ c_1,\dots,c_\mu\in\C.
\end{equation*}

Introducing 2-forms $\Omega=q\,dx\land dy$, $R_i=r_i\,dx\land dy$,
$i=1,\dots,\mu$, and 1-form $\eta=u\,dx+v\,dy$ the above identity can be
rewritten as \emph{division with remainder},
\begin{equation}\label{div-rem}
  \forall\Omega\in\Lambda^2\quad \exists\eta\in\Lambda^1:\qquad
\Omega=\omega\land\eta+R,\quad R=\sum_{i=1}^\mu c_i R_i\in\Lambda^2.
\end{equation}

In particular, any 2-form $\Omega$ can be divided with remainder by the
differential $dH$ of any polynomial $H\in\Lambda^0$. Note that the space
of remainders may be arbitrarily enlarged if necessary: the uniqueness of
the division \eqref{div-rem} will be lost then, but in exchange one may
get better norms of the ratio and remainder.

In general the procedure of division can be very delicate. However, if $H$
is transversal to infinity, then one can easily describe the outcome,
explicitly majorizing the \emph{degrees} of the remainder $R\in\Lambda^2$
and the incomplete ratio $\eta\in\Lambda^1$.

\begin{Lem}[see \cite{redundant}]\label{lem:div-forms}
If the polynomial $H\in\Lambda^0$ of degree $n+1$ is transversal to
infinity, then any 2-form $\Omega$ can be divided by $dH$ with the
incomplete ratio $\eta$ of degree $\le\deg\Omega-\deg H$ and the remainder
$R$ of degree $\le 2n$.
\end{Lem}

\begin{proof}
If $L\in\Lambda^0$ is a homogeneous polynomial of degree $n+1$ without
multiple linear factors, then the map between subspaces of homogeneous
forms,
\begin{equation}\label{sylvester}
 \mathfrak J=\mathfrak J_{L}\:\Lambda^1_{n}\to\Lambda^2_{2n+1},
 \qquad \eta\mapsto dL\land\eta,
\end{equation}
is an isomorphism. Indeed, in the bases consisting of all monomial forms
of the given degrees, the matrix of $\mathfrak J$ is the \emph{Sylvester
matrix} whose determinant is the \emph{resultant} of the two partial
derivatives $L_x$ and $L_y$. The assumption on $H$ implies that this
resultant is nonzero.

Therefore any homogeneous form $\Omega$ of degree exactly $2n+1$ is
divisible by $dL$. Any monomial form of degree greater than $2n+1$ can be
represented as a monomial 1-form of degree $2n+1$ times a monomial
function and hence is divisible by $dL$ \emph{without remainder} with the
same relation between the degrees,
$$
\deg\Omega\ge 2n+1\implies\Omega=dL\land\eta,\qquad
\deg\eta\le\deg\Omega-\deg L.
$$
Applying this observation to all homogeneous components of a 2-form
$\Omega=\Omega_0+\cdots+\Omega_{2n}+\Omega_{2n+1}+\cdots+\Omega_d=
R+\Omega_{2n+1}+\cdots+\Omega_d$, we prove the assertion of the Lemma for
homogeneous polynomials.

To divide a 2-form $\Omega$ of degree $\ge 2n+1$ by a nonhomogeneous
differential $dH=dL+\xi$,  where $L$ is the principal homogeneous part of
$H$, $\deg\xi<\deg L$, we divide it by $dL$ first, and then transform the
result as follows,
$$
\begin{gathered}
\Omega=dL\land\eta+R=(dH-\xi)\land \eta+R=dH\land \eta+\Omega',
 \\
\Omega'=R-\xi\land\eta, \qquad \deg\Omega'\le\max(\deg\Omega-\deg \xi+\deg
L,\deg R)<\deg\Omega,
\end{gathered}
$$
reducing division of $\Omega$ by $dH$ to division of another form $\Omega'$
of strictly inferior degree. Iterating this step, we prove the Lemma in the
general case. Notice that this is essentially the algorithm of division
with remainder of univariate polynomials.
\end{proof}

\begin{Rem}
Lemma~\ref{lem:div-forms} is an example of the \emph{redundant} division.
For a Hamiltonian $H$ of degree $n+1$ transversal to infinity, the gradient
ideal (ideal of 2-forms divisible by $dH$) has codimension $n^2$. Indeed,
the codimension is equal to the number $\mu$ of critical points of $H$ in
$\C^2$, counted with their multiplicities. This latter number is exactly
$n^2$ by virtue of B\'ezout theorem, since no critical points are allowed
to ``escape to infinity'' by the assumption on $H$.

On the other hand, the linear space of bivariate monomials of degree $\le
2n-2$ (the space of 2-forms $\Lambda^2_{\le 2n}$) is
$\nu=2n(2n-1)/2\approx2n^2$, roughly two times greater than $\mu$.

The irredundant analog of this theorem can be easily restored if
necessary. However, the choice of monomial 2-forms generating the
remainder, will depend on the principal homogeneous part $L$ of $H$.
\end{Rem}

\subsection{Gelfand--Leray residue and derivative}\label{sec:gelfand-leray}

\begin{Lem}
Let $\omega,\eta\in\Lambda^1$ be two polynomial 1-forms such that
$$
d\omega=dH\land \eta.
$$
Then for any continuous family $\delta(t)$ of 1-cycles on the level curves
$H^{-1}(t)$,
$$
\frac d{dt}\oint_{\delta(t)}\omega=\oint_{\delta(t)}\eta.
$$
\end{Lem}

The proof in the real case (assuming only smoothness of the forms) can be
achieved by integration of $d\omega$ over the annulus between $\delta(t)$
and $\delta(t+\Delta t)$ and passing to limit as $\Delta t\to0$.

This formula allows to differentiate explicitly Abelian integrals of a
form $\omega$, expressing the result as an Abelian integral once again if
$d\omega$ is divisible by $dH$. In fact, $\eta$ can be a \emph{rational}
1-form having all zero residues after restriction on each level curve
$H^{-1}(t)$.

\begin{Ex}
Let $H(x,y)=x^2+y^2$ and $\omega=y\,dx$. Then $\oint_{H=t}\omega=-\pi t$
(consider only the real values of $t$ and use the Stokes formula for the
circle positively oriented). Clearly, the form $\eta=\tfrac12y^{-1}\,dx$
satisfies the assumption of the Lemma, and indeed
$$
\oint_{H=t}\eta=\int_0^{2\pi}\frac12\cdot\frac{d\cos s}{\sin s}= -\pi.
$$
This example helps to memorize the order of the wedge multiplication in
the Gelfand--Leray formula.
\end{Ex}

\subsection{Derivation of the redundant Picard--Fuchs system}\label{sec:redundant}
The linear space $\Lambda^2_{\le 2n}$ of possible remainders occurring in
the division \eqref{div-rem}, is spanned by monomial forms
$x^ry^s\,dx\land dy$, $r+s\le 2n-2$. Denote its dimension $(2n-1)n$ by
$\nu$ and choose any monomial primitives $\omega_i$, $i=1,\dots,\nu$, so
that $d\omega_i$ span the quotient space $\Lambda^2_{\le
2n}/d\Lambda^1_{\le 2n}$ (modulo exact forms). Below we refer to
$\omega_i$ as the \emph{basic} forms.

Consider the 2-forms $H\,d\omega_i\in\Lambda^2_{\le 3n+1}$ and divide them
with remainder by $dH$:
\begin{equation}\label{div-bydH}
  H\,d\omega_i=dH\land\eta_i+R_i,\qquad i=1,\dots,\nu.
\end{equation}
By the assertion on the degrees, $\deg\eta_i\le\deg\omega_i+\deg H-\deg
dH\le \deg\omega_i\le 2n$, therefore each of the forms can be represented
(modulo an exact polynomial form) as a linear combination of the basic
forms,
$$
\eta_i=\sum_{j=1}^\nu b_{ij}\omega_j+dF_i,\qquad b_{ij}\in\C,\
F_i\in\Lambda^0_{\le 2n}.
$$
Similarly, being all of degree $\le 2n$, the remainders
$R_i\in\Lambda^2_{\le 2n}$ can be represented as linear combinations of the
forms $d\omega_i$:
$$
R_i=\sum_{j=1}^\nu a_{ij}\,d\omega_j,\qquad a_{ij}\in\C.
$$
Let $\delta(t)$ be any continuous family of cycles. Then for any
$t\notin\S$ the forms
$$
H\,d\omega_i-\sum_{j=1}^\nu a_{ij}d\omega_j,\qquad \forall i=1,\dots,\nu,
$$
are all divisible by $dH$ with the ratios being cohomologous to
$\sum_{j=1}^\nu b_{ij}\omega_j$. Denote
$$
 X_i(t)=\oint_{\delta(t)}\omega_i,\qquad i=1,\dots,\nu.
$$
Note that integration of a form $H\omega$ over any cycle
$\delta\subset\{H=t\}$ yields $t\oint_\delta\omega$, since $H$ is constant
on the cycle. Integrating both sides over the oval $\delta(t)$ and using
the Gelfand--Leray formula, we arrive to the identities
\begin{equation}
  t\dot X_i(t)-\sum_{j=1}^\nu a_{ij}\dot X_j(t)=\sum_{j=1}^\nu
  b_{ij}X_j(t),
\end{equation}
which means that the \emph{column} vector $(X_1,\dots,X_\nu)$ of Abelian
integrals satisfies the system of linear ordinary differential equations
\begin{equation}\label{rpf}
  (tE-A)\dot X=BX,\qquad X\in\C^\nu,\quad A,B\in\Mat_{\nu\times\nu}(\C),
\end{equation}
with the constant matrices $A=\|a_{ij}\|$, $B=\|b_{ij}\|$ as parameters
($E$ is the identity matrix). Writing
\begin{equation}\label{adjugate}
 (tE-A)^{-1}=\frac1{\chi(t)}\cdot P(t),\qquad P=\sum_{k=0}^{\nu-1}P_k t^k,
 \ \chi(t)=\det (tE-A),
\end{equation}
where $P(t)$ is a matrix polynomial of degree $\nu-1$, the adjugate matrix
for $tE-A$, one sees immediately that \eqref{rpf} is a system of linear
ordinary differential equations with rational coefficients.

\subsection{Hypergeometric systems}
The form \eqref{rpf} is rather specific: for instance, all singular points
of this system are Fuchsian (simple poles), including the point at
$t=\infty$. This is obvious if the spectrum of $A$ is simple, but holds
true in the general case as well, as follows from the explicit formula for
$(tE-A)^{-1}$ for $A$ in the Jordan normal form.

It would be appropriate to remark here that the residues $A_j$ of the
Fuchsian system \eqref{hyperel} are invariant by any conformal change of
the independent variable $t$. In the case of hypergeometric systems
\eqref{rpf} the point $t=\infty$ is distinguished: the residues $A_j$ at
all finite points $t_j\in\operatorname{Spec}A$ have rank $1$ for a generic
matrix $A$, whereas the rank of the residue $A_\infty=-\sum_1^\mu A_j$ at
infinity is generically full. Thus the natural symmetry group of
hypergeometric systems is not the full group of conformal automorphisms of
$\C P^1$, but rather the affine group of transformations $t\mapsto at+b$,
$a,b\in\C$ fixing the point $t=\infty$. Making an affine transformation
transforms the system \eqref{rpf} into the system $(tE+A')\dot X=B'X$ with
the same matrix $B'=B$ and $A'=a^{-1}(A-bE)$.

In the subsequent sections it will be shown that in order to estimate the
number of isolated zeros of solutions to Fuchsian systems, it is
sufficient to know the norms of the residue matrices. As follows from the
explicit inversion formula \eqref{adjugate}, norms of the residues can be
bounded if the norm $\|A\|$ is bounded from above and pairwise distances
between the critical points $t_j$ are bounded \emph{from below}. Of
course, choosing a suitable affine transformation as above, one can change
the norm $\|A\|$, but at the same rate the distances between the singular
points (eigenvalues of $A$) will be affected. In other words, the norm of
the matrix $\|A\|$ should be majorized relative to the spread of its
eigenvalues.

\subsection{Explicitness and bounds}
The above derivation does not involve any existential assertion: all
constructions are completely transparent and allow for explicit bounds,
say, on the norms of the matrices $A,B$ from \eqref{rpf}. To do this, we
introduce the norms on the ring of polynomials, letting $\|p\|$ being the
sum of absolute values of all its coefficients. This norm is
multiplicative, $\|pq\|\le\|p\|\,\|q\|$, and extends on polynomial
$k$-forms, remaining multiplicative with respect to the wedge product. The
exterior derivative is a bounded operator on forms of bounded degrees.

The entries of the matrices $A,B$ appear as coefficients of linear
expansion of a known $1$-form in the chosen basis. But since the basic
forms $\omega_i$ are \emph{monomial} with coefficients equal to $1$, to
majorize these entries it is sufficient to majorize the norms of the
2-forms $R_i$ (the remainders) and 1-forms $\eta_i$ (incomplete ratios).
In other words, one has to control only the division step, since
multiplication by $H$ is an operator whose norm is no greater than $\|H\|$.

The division step is also rather transparent, its well-posedness being
determined by the norm of the inverse Sylvester matrix $\mathfrak J_L$
from \eqref{sylvester} and the norm of the non-principal terms $\|H-L\|$,
where $L$ is the principal homogeneous part of $H$. The group of affine
transformations of the complex plane $\C^2$ naturally acts on the space of
all Hamiltonians of degree $n+1$ not affecting the critical values of $H$.
By an appropriate transformation of this group, one can always achieve the
normalizing condition $\|\mathfrak J_L^{-1}\|=1$ that is a condition on
the principal homogeneous part $L$. The problem on bounding the norm
$\|A\|$ is reduced therefore to studying how the magnitude of the
non-principal coefficients of $H-L$ may affect configuration of the
critical values of $H$ subject to the above normalizing condition on the
principal part $L$. The problem can be explicitly solved for the
univariate polynomials, implying an answer in the hyperelliptic case as
well \cite{redundant}. Moreover, one can show that if all critical values
of a bivariate polynomial transversal to infinity coincide, then
necessarily the polynomial $H$ must coincide with its principal part,
being thus homogeneous, eventually after a suitable parallel translation
in the $(x,y)$-plane. The inequality between the non-homogeneity of a
bivariate polynomial and the spread of its critical values is still
unknown, see \cite{redundant} for partial results.

\subsection{Preliminary conclusion}
The tangential Hilbert problem for gen\-eric Hamiltonians, gets reduced to
the question about the number of isolated zeros of linear combinations of
functions satisfying together a system of linear ordinary differential
equations with rational coefficients.

The procedure of derivation of this system is very transparent. In
particular, it can be written in the hypergeometric form \eqref{rpf} with
explicit bounds on the norms of the corresponding matrices $A$ and $B$.
These bounds in turn imply that when reduced to the Fuchsian form
\eqref{fuchsian}, the system will have the residue matrices $A_j$ bounded
(as usual, in the sense of the norms) in terms of the inverse distance
$\max_{i\ne j}\left\{|t_i-t_j|^{-1},|t_i|\right\}$ between singular points
of the system.

When this inverse distance tends to infinity (which corresponds to
confluence of singular points), the resulting bounds on the norms
$\|A_j\|$ of the residues in \eqref{fuchsian} explode. However, this
explosion is of a very specific nature: the monodromy group of the system
with one or several confluent singularities, remains the same. In
particular, the spectral data of the residues remain bounded.

In the subsequent lectures we will find out how far away are these
conditions from sufficient conditions allowing for an explicit solution of
the tangential Hilbert problem.

\chapter{Quasialgebraicity of function fields}

Starting from this moment, we will pursue the same path towards the
tangential Hilbert problem, but this time in the opposite direction.
Namely, we will establish conditions on Fuchsian systems guaranteeing that
their solutions are similar to algebraic functions, in particular, admit
explicit bounds for the number of isolated zeros.

\section{Functional fields and their quasialgebraicity}\label{sec:quasialg}

The main objects of study in this section are functional fields obtained
by adjoining one or several analytic (in general, multivalued) functions
to the field $\C$ of complex numbers (or slightly more generally to that of
rational functions $\C(t)$). Such fields admit filtration (grading) by
degrees.

The goal is to obtain conditions on the field (in terms of properties of
the adjoined functions) guaranteeing that the question on the \emph{global}
number of isolated zeros of functions from this field can be
\emph{algorithmically} (effectively) solved. An accurate definition will be
given at an appropriate moment, after explaining all technicalities
pertinent to the problem. We begin by examples illustrating the goals.

\subsection{Algebraic functions}
The field of rational functions $\C(t)$ in one variable $t$ possesses the
following obvious but nevertheless remarkable property:
\begin{enumerate}
 \item any element $f(t)=p(t)/q(t)$, $p,q\in\C[t]$ from this field has a
well defined degree $\deg f=\max(\deg p,\deg q)$ (assuming that the
representation is irreducible), and
 \item the number of isolated zeros of $f$ on the whole projective line
$\C P^1=\C\cup\{\infty\}$ is no greater than $\deg f$ (actually, equal to
it if counted with multiplicities).
\end{enumerate}
In other words, there is a direct relationship between the combinatorial
complexity of representation of $f$ in the field (i.e., the number of
field operations necessary to produce $f$ from constants and the
independent variable $t$), and its analytic complexity measured by the
number of isolated zeros.

This example can be easily generalized by considering fields generated by
one or several algebraic functions.

Let $f_i(t)$, $i=1,\dots,n$ be algebraic functions of one variable, defined
implicitly by the polynomial equations
\begin{equation*}
  P_i(t,x_i)=0,\qquad P_i\in\C[t,x_i],\quad \deg P_i=d_i,\quad i=1,\dots,n,
\end{equation*}
with respect to $x_i$. Consider the ring $\C[f_1,\dots,f_n]$ formed by
polynomial combinations of the functions $f_i$, and the corresponding
field of fractions $\C(f_1,\dots,f_n)$. Both consist of analytic
multivalued functions ramified over a finite point set on the projective
line $\C P^1$, though the number of distinct branches of every function is
finite (no more than $d_1\cdots d_n$).

One can define unambiguously the degree of functions in this ring and the
respective field as the degree $d$ of the \emph{minimal} representation
\begin{equation*}
  f=\sum_{|\alpha|\le d}c_\alpha f^\alpha,\qquad
  \alpha=(\alpha_1,\dots,\alpha_n)\in\mathbb Z^n_+,\ c_\alpha\in\C,\
  f^\alpha=\prod_1^n f_i^{\alpha_i}.
\end{equation*}

\begin{Prop}\label{prop:bezout}
The total number of isolated zeros of a polynomial combination
$f(t)=P(f_1(t),\dots,f_n(t))$, $f\in\C[f_1,\dots,f_n]$ of degree $\le d$
on all branches of this function, does not exceed $d\cdot d_1\cdots d_n$.
\end{Prop}

\begin{proof}
This is an immediate corollary to the B\'ezout theorem applied to the
system of algebraic equations
\begin{equation*}
\left\{
\begin{aligned}
  P_1(t,x_1)&=0,
  \\
  \vdots\qquad&
  \\
  P_n(t,x_n)&=0,
  \\
  P(x_1,\dots,x_n)&=0.
  \end{aligned}\right.
\end{equation*}
The same proof actually works for a more general case of functions $f_i$
defined by a system of algebraic equations $P_i(t,x_1,\dots,x_n)=0$,
$i=1,\dots,n$, involving all functions \emph{simultaneously}.
\end{proof}

\begin{Rem}
The function $f_0(t)=t$ is clearly algebraic, and if required, we can
always assume it being among the collection of the functions $f_i$, thus
avoiding particular cases and awkward notation. This agreement will allow
us to assume that all functional fields are extensions of the field
$\C(t)$.
\end{Rem}

\subsection{Existential quasialgebraicity}
Some parts of the above construction can be reproduced in a completely
general context. Let $U\subset\C$ be an open domain (for simplicity assume
it to be bounded, connected and simply connected).

Consider an arbitrary collection of $n$ functions
$F=\{f_1(t),\dots,f_n(t)\}$ analytic in $U$. They define the ring $\C[F]$
and the corresponding field $\C(F)$ of functions meromorphic in $U$. As
before, for any function $f\in\C[F]$ one can define its degree as the
minimal possible degree of the polynomial $P\in\C[x_1,\dots,x_n]$ realizing
the given function $f(t)=P(F(t))=P(f_1(t),\dots,f_n(t))$ (there can be
algebraically dependent functions among the generators). This grading
extends naturally for rational combinations from $\C(F)$.

\begin{Prop}\label{prop:exist-quasialg}
For any compact $K\Subset U$ there exists a counting function
$C=C_K\:\mathbb N\to\mathbb N$, taking only finite values $C_K(d)<+\infty$
for any finite $d$, such that the number of isolated zeros of any function
$f$ in $K$ can be at most $C_K(d)$:
\begin{equation*}
  \deg f\le d\implies \#\{t\in K\: f(t)=0\}\le C_K(d).
\end{equation*}
\end{Prop}

\begin{proof}
First we show that under the assumptions of the Theorem, the number of
isolated zeros of any \emph{linear} combination $f_c(t)=\sum c_i f_i(t)$
with complex constant coefficients $c_1,\dots,c_n\in\C$, is bounded in any
compact $K$ uniformly over all such linear combinations. Indeed, without
loss of generality we may assume that the functions $f_i$ are linear
independent---this does not affect the supply of all linear combinations.
Next, it is sufficient to consider only combinations with coefficients on
the unit sphere, satisfying the equality $\sum_j |c_j|^2=1$. The functions
$f_c$ for such $c$ are all different from identical zero, hence each of
them has only a finite number of isolated zeros in the compact $K$
(accumulation of roots to the boundary of $K$ is forbidden). Now the
standard semicontinuity arguments using compactness of the unit sphere,
prove that the number of zeros of all $f_c$ is uniformly bounded.

To deal with arbitrarily polynomial combinations, we can treat them as
linear combinations of \emph{monomials} $f^\alpha(t)$, $|\alpha|\le d$,
reducing the general case to the already studied one.
\end{proof}

\begin{Rem}
One can easily recognize in this demonstration some minor variations on
the theme already exposed in \secref{sec:analytic-functions}. Of course,
Proposition~\ref{prop:exist-quasialg} follows from
Theorem~\ref{thm:finite-cyclicity}, since the parameters $c_\alpha$ can be
considered varying over the compact sphere. Here we could explicitly avoid
dealing with functions vanishing identically for some values of the
parameters, simplifying considerably the proof.
\end{Rem}

\subsection{Comparison}
Two above finiteness assertions, Propositions~\ref{prop:bezout} and
\ref{prop:exist-quasialg}, differ in two important instances:
\begin{enumerate}
 \item the bounds on roots of algebraic functions are \emph{global},
i.e., valid on the maximal domain of definition of the functions from the
field $\C(F)$, whereas the bounds on the roots of arbitrary analytic
functions in general depend on the choice of the compact $K$ and can in
many cases blow up as the compact approaches the boundary of the maximal
domain $U$;
 \item the bounds on roots of algebraic functions are given by an explicit
formula involving some basic parameters defining the field, whereas the
function $C_K(d)$ is totally existential (see below).
\end{enumerate}

The nature of the counting function $C_K(d)$ from
Proposition~\ref{prop:exist-quasialg} remains totally non-effective. One
can easily construct examples of functions (even entire functions) such
that the growth of $C_K(d)$ will be arbitrarily fast \cite{jde-96}.

\subsection{In search of quasialgebraicity: reappearance of Picard--Vessiot
 extensions}\label{sec:insearch}
Our goal is to provide sufficient conditions on the functions $f_i$
guaranteeing that the corresponding field will be similar to the field
obtained by adjoining algebraic functions. This condition, still
understood informally until made precise in \secref{sec:quasialgebraicity},
will be referred to as \emph{quasialgebraicity} of the function field. The
accurate definition is postponed since it  involves some technical details.

However, even prior to giving any accurate formulation, the class of
function fields among which one could hope to find nontrivial cases of
quasialgebraicity, can be substantially restricted.

First, \emph{the generating functions $f_1,\dots,f_n$ must be multivalued}
(ramified). Indeed, a single-valued function having at most polar
singularities on the projective line $\C P^1$ (recall that we are looking
for \emph{global} bounds, hence the functions should be defined globally),
must be rational. Any field generated by such functions, is a subfield of
$\C(t)$ and hence we get nothing new.

On the other hand, if one of the functions $f_i$ has an essential
singularity on $\C P^1$, then by classical theorems of complex analysis
this function near such point must take infinitely many times almost all
values, hence one can easily construct an \emph{individual} polynomial
combination having infinitely many roots accumulating to the essential
singularity. This precludes quasialgebraicity whatever exact meaning it
may be assigned.

Thus any field (or ring, what is almost the same for our purposes)
exhibiting nontrivial quasialgebraicity, must consist of functions
ramified over some finite set $\S\subset\C P^1$. As above, globality means
that functions should be analytically continuable along any path avoiding
the ramification locus $\S$.

The possibility of analytic continuation along paths (and loops)
introduces an additional structure, the monodromy group action. Choose
arbitrarily a nonsingular point $a\notin\S$. Then any element from the ring
$\C[F]$ can be identified with the full analytic continuation of its germ
at $a$. Denote as in \secref{sec:topology} the monodromy operator
associated with a loop $\gamma\in\pi_1(\C\ssm\S,a)$ by $\Delta_\gamma$.
Then it would be natural to assume that $\C[F]$ is closed (invariant) by
analytic continuations, that is,
\begin{equation}\label{monodromy-invariance}
  \Delta_\gamma f_i\in\C[f_1,\dots,f_\mu]
  \qquad\forall\gamma\in\pi_1(\C P^1\ssm\S,a).
\end{equation}
Moreover, since the filtration of the ring $\C[F]$ by degrees should be
well-defined, it must be preserved by analytic continuations, that is,
continuation of a polynomial $f=P(f_1,\dots,f_\mu)$ (of minimal degree
among all polynomials representing the given $f$) along any loop should
again be a polynomial of the same degree. In particular, \emph{the
$\C$-linear span of the germs $f_1,\dots,f_\mu$ in the space of all
analytic germs must be invariant by all monodromy operators}. In other
words, there should exist invertible square matrices $M_\gamma$ such that
\begin{equation*}
  \Delta_\gamma(f_1,\dots,f_n)=(f_1,\dots,f_n)\cdot M_\gamma.
\end{equation*}

In the same way as in \secref{sec:topology}, this implies that the
functions $f_1,\dots,f_\mu$ must be solutions of a linear ordinary
differential equation with single-valued coefficients on $\C P^1\ssm\S$.
As before, allowing essential singularities of the coefficients would
immediately relinquish any control over nonaccumulation of roots of
solutions, hence we arrive to the following important conclusion: \emph{To
be quasialgebraic, the field $\C(F)=\C(f_1,\dots,f_\mu)$ must be a
Picard--Vessiot extension of $\C$ or $\C(t)$ obtained by adjoining
solutions of a linear ordinary differential equation with rational
coefficients}.

This fact, in particular, implies that adding to the generating tuple of
functions $(f_1,\dots,f_n)$ their derivatives of orders $\le n-1$ will
make the field $\C(t,f_1,\dots,f_n)$ a \emph{differential field} (closed by
differentiation). Without loss of generality we can assume that such
completion was already done and we deal with the differential field
$\C(X)$ (and the corresponding ring $\C[X]$) obtained by adjoining to $\C$
all entries of a fundamental matrix solution $X(t)$ for a system of first
order linear ordinary differential equations with rational coefficient
matrix,
\begin{equation}\label{linsys}
  \dot X(t)=A(t)X(t),\qquad A(t)\in\Mat_{\mu\times\mu}(\C(t)).
\end{equation}
In the particular case we are discussing, one can take $X$ to be the
Wronskian matrix of the collection $\{f_i\}$, $X_{ij}(t)=f_j^{(i-1)}(t)$,
$i,j=1,\dots,n$.

\subsection{Quasialgebraicity of Fuchsian systems: examples and
counterexamples} As in the case of single-valued functions, a precondition
for quasialgebraicity is nonaccumulation of roots of solutions to singular
points. Singularities of linear systems \eqref{linsys} can be easily
described: they occur only at the poles of the coefficients matrix $A(t)$
or rather at the poles of the matrix-valued differential 1-form $A(t)\,dt$
on the Riemann sphere $\C P^1$. The corresponding classical theory
\cite{hartman} distinguishes between two types of singularities:
\begin{itemize}
 \item those exhibiting at most polynomial growth of entries of the
fundamental matrix $X(t)$ and its inverse $X^{-1}(t)$ and called
\emph{regular singularities}, and
 \item those exhibiting abnormally fast (faster than polynomial) growth
of solutions, called \emph{irregular singularities}.
\end{itemize}

\begin{Rem}
In order to measure growth rate of multivalued functions near a
ramification point, they should be restricted on a sector bounded by two
rectilinear rays with the vertex at this point. Otherwise one can
construct a curve approaching the singular point while spiraling around it
in such a way that the growth in terms of the distance to the singularity
will be arbitrarily fast even for the most innocent multivalued function
$\ln t$.
\end{Rem}

The dichotomy between regular and irregular singularities is closely
related to dichotomy between poles and essential singularities for
single-valued functions. Consider a neighborhood of a singular point
$t_*\in\S$, assuming for simplicity that $t_*=0$. Let $M$ be the monodromy
operator associated with a small loop around the origin, so that
\begin{equation}\label{local-monodromy}
  \Delta_\gamma X(t)=X(t)M,\qquad \det M\ne0.
\end{equation}
Let $A\in\Mat_{\mu\times \mu}(\C)$ be any matrix logarithm of $M$, so that
\begin{equation*}
  \exp 2\pi i A=M.
\end{equation*}
Then the multivalued matrix function $t^A=\exp (A\ln t)$ also has the same
monodromy,
\begin{equation*}
   \Delta_\gamma \,t^A=\exp(A\Delta_\gamma\ln t)=\exp(A(\ln t+2\pi
  i))=t^A M=Mt^A.
\end{equation*}
Therefore the matrix ratio $H(t)=X(t)t^{-A}$ is single-valued in a small
punctured neighborhood of $0$. If the origin is a regular singularity, then
$H(t)$ has at most a pole and hence by choosing a different valuation of
the logarithm and replacing $t^A$ by $t^{kE+A}$ for sufficiently large
natural $k$, one can make $H(t)$ holomorphic at the origin. By explicit
elementary formulas for matrix exponents, one can derive from this a local
representation of entries of the matrix $X(t)$ in the form
\begin{equation*}
  \sum_{k,\l} h_{k\l}(t)\,t^\l\ln^{k-1} t
\end{equation*}
with the exponents $\l$ ranging over the spectrum of $A$, the natural $k$
being no greater than the maximal size of Jordan blocks of $A$ and
$h_{k\l}(t)$ holomorphic at $t=0$. Notice the remarkable coincidence with
\eqref{asymp-expan} and Exercise~\ref{ex:dulac-series}: the latter implies
that  at least in the situation when all functions $h_{k\l}$ and all
eigenvalues of $A$ are real, isolated zeros cannot accumulate to the
regular singularity at the origin.

If, on the other hand, the singularity is irregular, then at least some
entries of the matrix $H(t)$ must exhibit essential singularity at the
origin. In the same way as with single-valued function, one can in this
case construct functions from $\C(X)$ that would have infinitely many
roots accumulating to the origin.

\begin{Ex}\label{ex:essential-sing}
The linear system
\begin{equation}\label{essent-sing}
 \left\{
  \begin{aligned}
  \dot x_1&=t^{-2}{x_1},
  \\
  \dot x_2&=0
  \end{aligned}
  \right.
\end{equation}
generates the field $\C(X)$ containing a function $f(t)=\exp(-1/t)-1$
whose zeros at the points $t_m=(2\pi i m)^{-1}$, $m=\pm1,\pm2,\dots$
accumulate to two essentially singular points at $t=0$ and $t=\infty$
along the imaginary axis.
\end{Ex}

Thus occurrence of irregular singularities destroys any hope to achieve
quasialgebraicity, and we are left with the class of Picard--Vessiot
extensions exhibiting only regular singularities.

There is a simple sufficient condition guaranteeing that a singular point
of the linear system \eqref{linsys} is regular. By the Fuchs theorem
\cite{hartman}, if $A(t)$ has a \emph{simple pole} (of the first order) at
a point $t_*$, then this point is a regular singularity. Such
singularities are called \emph{Fuchsian}.

The inverse to the Fuchs theorem is in general not true: there exist
regular non-Fuchsian singularities. But a system exhibiting regular
singularity at $t_*$, by a \emph{meromorphic} (locally near $t_*$) linear
transformation $X(t)\mapsto R(t)X(t)$ can be reduced to the system having a
simple pole at $t_*$. Thus on the level of \emph{local meromorphic
equivalence} there is no difference between Fuchsian and regular
singularities.

For a globally defined system exhibiting several regular singularities at
finitely many points $\S=\{t_1,\dots,t_d\}$ on the projective line, one can
ask whether there exists a globally meromorphic (hence \emph{rational})
linear transformation simultaneously taking all regular singularities into
Fuchsian ones. The problem (that constitutes a part of the so called
Riemann--Hilbert, or 21st Hilbert problem) turns out to be very delicate,
the result depending essentially on the structure of the monodromy group
of the system, and not always admitting solution, as shown recently by
A.~Bolibruch \cite{bolibruch:icm-94,bolibr:umn}. However, from the
classical result by Plemelj \cite{plemelj,forster} it follows that if an
additional singular point is allowed to be created anywhere, then the
answer is always positive and a rational matrix function $R(t)$ can be
found such that $Y(t)=R(t)X(t)$ satisfies a system of linear ordinary
differential equations having only Fuchsian singularities on the whole
projective line $\C P^1$. Clearly, if $t\in\C(X)$ (which we may always
assume without loss of generality, as noted above), then $\C(X)=\C(Y)$ and
hence when discussing quasialgebraicity, one can deal with Fuchsian
systems only.

A Fuchsian system with $d$ finite singular points $t_1,\dots,t_d$ can
always be written in the form
\begin{equation}\label{fuchsian2}
 \begin{gathered}
  \frac{dx}{dt}=A(t)x,\qquad A(t)=\sum_{j=1}^d\frac{A_j}{t-t_j},
  \\
  x\in\C^n,\quad A_j\in\Mat_{n\times n}(\C),\quad t\in\C,\quad
  t_1,\dots,t_d\in\C,
 \end{gathered}
\end{equation}
explicitly indicating the corresponding residue matrices $A_j$ (we return
to the initial notation $X$ for the dependent variables). Thus the natural
problem arises, \emph{When the Picard--Vessiot extension $\C(X)$
constructed by adjoining all components of the fundamental matrix solution
of a Fuchsian system \eqref{fuchsian2}, is quasialgebraic?} Note that this
field depends not on the choice of the matrix solution $X(t)$ but rather
on its ``logarithmic derivative'' $A(t)=\dot X(t)X^{-1}(t)$ which is a
rational matrix function. The Fuchsian system \eqref{fuchsian2} is
determined by its dimension and the collection of algebraic data
$\{A_i,t_i,\, i=1,\dots,d\}$. The bound for the number of zeros should be
given in terms of these algebraic data.

\begin{Ex}
The simplest class of Fuchsian systems is that having only two
singularities (one simple pole of the matrix $A(t)$ on the whole line $\C
P^1$ is impossible since the sum of all the residues, including the one at
infinity, must be zero). By a conformal transformation of the independent
variable the two points can be placed at $t=0$ and $t=\infty$. The
corresponding system will then take the \emph{Euler form},
\begin{equation}\label{euler}
  \dot X=\frac{A}{t}\cdot X,\qquad A\in\Mat_{\mu\times\mu}(\C).
\end{equation}

The Euler system \eqref{euler} can be immediately integrated: $X(t)=t^A$.
The associated Picard--Vessiot extension has the form that is already
familiar:
\begin{equation}\label{euler-field}
  \C(X)=\C(\{t^\l\ln^{k-1}t\}_{\l,k}),\qquad \l\in\operatorname{Spec}A,
  \quad k\le m(\l),
\end{equation}
where $m(\l)\ge 1$ is the maximal size of Jordan cells corresponding to
the eigenvalue $\l$.

Roots of \emph{real} functions from the field $\C(X)$ cannot accumulate to
the origin, if all eigenvalues of $A$ are real
(Example~\ref{ex:dulac-series}). We will see that the ``Euler field'' for a
system with only real spectrum is indeed quasialgebraic.

On the other hand, the field $\C(t^{\sqrt{-1}},t^{-\sqrt{-1}})$ associated
with the $2\times2$-Euler system with non-real spectrum
$\{\pm\sqrt{-1}\}$, contains the function $\cos\ln t$ that obviously has
infinite number of roots accumulating to $t=0$.

This suggests that quasialgebraicity of the fields $\C(X)$ should be
somehow related to spectral properties of the residue matrices.
\end{Ex}

\subsection{Counting zeros of multivalued functions globally}
Since functions from the Picard--Vessiot extensions are in general
multivalued, counting their zeros ``on the whole Riemann sphere'' should
be supplied with a precise meaning.

\begin{Ex}
The function $f(t)=t^\l-1$ for a real irrational $\l$ has infinitely many
different roots $t_k=q^k$, $q=\exp \frac{2\pi i}{\l}$, $k\in\mathbb Z$,
lying on different branches of this function.

In order to single out a branch, it is required to choose a simply
connected domain $U\subset\C\ssm\S$. Usually this choice is done by
slitting the complex plane along segments or rays with endpoints in the
singular locus $\S$. In the above example, slitting $\C$ along the negative
semiaxis yields a simply connected domain $U$ with the following property:
\emph{any branch of the function $t^\l-1$ in this domain has at most
$\lfloor \l^{-1}\rfloor$ roots in it} (the same bound holds also for
$-U$). Here and below $\lfloor a \rfloor$ stands for the integer part of a
real number $a$. This explicit bound is global in the sense that $U$ and
$-U$ together cover $\C\ssm\S$.
\end{Ex}

However, the choice of the simply connected domain $U$ can affect even
finiteness of bounds on the number of zeros.

\begin{Ex}
If $\l$ is a non-real number, then the roots fall on the logarithmic spiral
(since $|q|\ne 1$ in this case). If $U$ is chosen by slitting $\C$ along
positive or negative semiaxis, then the number of isolated roots on any
branch will still be bounded by $\lfloor (\Re \l)^{-1}\rfloor<+\infty$ for
$\l\notin i\R$.

On the other hand, one can choose a simply connected spiral domain
containing the origin on the boundary, so that it would simultaneously
contain an infinite number of the points $t_k=q^k$.

Though the spiral slits are not algebraic (even not analytic) curves, this
example can be easily modified to construct \emph{polygonal} simply
connected domains containing as many roots of the function $t^\l-1$ as
necessary, provided that the number of sides of the polygon can be
arbitrary.
\end{Ex}

These examples motivate the following definition.

\begin{Def}
An analytic multivalued function $f\:\C P^1\to\C$ ramified over a finite
set $\S\subset\C P^1$ is said to admit a \emph{global bound on the number
of isolated zeros}, if there exists a natural number $N<+\infty$ such that
the number of isolated zeros of any branch of $f$ in any \emph{open
rectilinear triangle} $T\subset\C\ssm\S$ has no more than $N$ isolated
zeros, the bound being uniform over all such triangles.
\end{Def}

\begin{Rem}
Having this property, one can place an explicit upper bound on the number
of isolated roots of $f$ in any simply connected domain $U$ bounded by
finitely many algebraic curves of known degrees. The number of curves,
their degrees and the number of ramification points will explicitly enter
the answer together with the number $N$.

Indeed, the domain $\C\ssm\S$ can be triangulated into finitely many
triangles as above. The number of simply connected components of any
intersection $U\cap T$ with each triangle of the triangulation can be
easily majorized by the B\'ezout theorem, in terms of the above data. The
total number of such components should be multiplied by $N$ to construct
the required upper bound.
\end{Rem}

\begin{Ex}
The function $t^\l-1$ admits a global upper bound for the number of
isolated roots for all $\l$ with $\Re\l\ne0$. For $\l=\sqrt{-1}$ any real
interval $(0,a)$, $a>0$, contains infinitely many roots of this function.
\end{Ex}

\begin{Rem}
In practice, however, we will always choose a convenient system of
rectilinear slits as in the textbooks on complex variables: if, for
example, an upper bound for the number of zeros is known in both
$U=\C\ssm\R_+$ and $-U=\C\ssm\R_-$, then the number of zeros in any
triangle not containing the origin, does not exceed the maximum of these
two bounds, since such triangle cannot intersect both the positive and
negative semiaxes simultaneously.
\end{Rem}

\section{Digression on computability}

\subsection{Hierarchy of computability}
The discussion in the previous section addressed the issue of
\emph{globality} of bounds on the number of zeros of multivalued
functions. Now we pass to a brief discussion of \emph{computability} of
these bounds. This would lead us again to functions defined by
differential equations with algebraic right hand side parts.

The bounds for functions of the field $\C(f_1,\dots,f_\mu)$, besides being
dependent on the degree $d$ of these functions, should obviously depend on
the field itself, that is, ultimately, on the choice of the generators
$f_i$. While computability of functions of natural arguments is a well
developed area (several notions of computability can be distinguished, see
below), dependence on functional parameters is definitely outside the
scope of any theory (e.g., the number of roots of a function is not a
computable characteristic \emph{per~se}).

An intermediate position occupy ``computable'' functions of one or several
real or complex variables. Here we take the most restrictive attitude,
legalizing only \emph{polynomial functions}, the absolute value
$|\,\cdot\,|$ of a real or complex variable, and the integer part
$\lfloor\cdot\rfloor$ of a real variable.

Returning to functions of one or several natural arguments, one can
classify them in an ``increasing order of computability'' into:
\begin{enumerate}
 \item (general) recursive functions, whose values can be computed
 for any combination of arguments in a finite time by some deterministic
 algorithm, Turing machine etc.;
 \item primitive recursive functions, that can be defined by one or several
 iterated inductive rules of the form
\begin{equation}\label{prim-rec}
  \f(k+1,m)=\varPhi(k,m,\f(k,m)),\qquad k\in\mathbb N,\ m\in\mathbb N^{n},
\end{equation}
 assuming that the functions $\varPhi$ and $\f_1=\f(1,\cdot)$
 are already known;
 \item elementary functions, given by finite compositions of powers,
 exponents, algebraic operations \emph{etc}.
\end{enumerate}

\begin{Ex}[tower functions]
Consider the function $\tau(k,n)$ of two integer arguments, defined by the
recursive rule
\begin{equation}\label{tower-rule}
  \tau(0,k)=k,\qquad \tau(n+1,k)=2^{\tau(n,k)}.
\end{equation}
These rules define \emph{tower functions} (iterated exponents):
$\tau(n,\cdot)$ is a tower of height $n$ and for each particular $n$ is an
elementary function of the second argument $k$:
\begin{equation*}
  \tau(1,k)=2^k,\quad \tau(2,k)=2^{2^k},\quad
  \tau(3,k)=2^{2^{2^k}},\quad\dots
\end{equation*}
However, as a function of the first argument, it is not elementary.
\end{Ex}

\begin{Ex}[Ackermann generalized exponential]
Consider the function $A(z,x,y)$ of three natural arguments, defined by the
recursive rules and initial conditions as follows,
\begin{equation}\label{ackermann}
\begin{gathered}
  A(z+1,x,y+1)=A(z,x,A(z+1,x,y)),
  \\
  A(0,x,y)=y+1,\ A(1,x,0)=x,\ A(2,x,0)=0,
  \\
  A(z,x,0)=0,\qquad \forall z\ge 2.
  \end{gathered}
\end{equation}
These rules define each function $\f=A(z+1,\cdot,\cdot\,)$ for any
particular $z$ unambiguously provided that the function
$\varPhi=A({z},\cdot,\cdot\,)$ is already defined, by the \emph{simple}
recursive rule $\f(x,y+1)=\varPhi(x,\f(x,y))$. In other words, the set of
conditions defines $A$ for all nonnegative combinations of arguments.

One can easily check that
\begin{equation*}
  \begin{aligned}
A(1,x,y)&=x+1+\cdots+1=x+y,
\\
A(2,x,y)&=x+x+\cdots+x=xy,
\\
A(3,x,y)&=xx\cdots x=x^y,
\end{aligned}
\qquad
\begin{gathered}
A(4,x,y)=\underset{y\text{ times}}{\underbrace{x^{x^{\invddots^x}}}}
\\
\text{etc.}
\end{gathered}
\end{equation*}
This suggests that as a function of the first argument, $A$ is not
primitive recursive. The reason is that, unlike in the simple recursive
rule \eqref{prim-rec}, the function $A$ is applied to itself (this does
not prove that $A$ is not primitive recursive, but the fact remains true
and can be rigorously demonstrated). It can be also shown that as a
function of the first argument, the Ackermann generalized exponential grows
faster than any primitive recursive function.
\end{Ex}

\subsection{Transcendental functions defined by algebraic data}
The above brief discussion suggests that in order to speak about computable
bounds depending on several analytic functions as parameters, these
analytic functions must themselves be defined in terms of finitely many
integer, real or complex parameters.

The simplest example of such functions are polynomials (their coefficients
play the role of the parameters) and, slightly more generally, algebraic
functions. However, this example does not allow to produce nontrivial
quasialgebraic fields.

Differentiation of algebraic functions leaves them algebraic. On the
contrary, integration (taking primitives) in general results in
transcendental functions. Another possibility of constructing
transcendental functions from algebraic ones is exponentiation. However,
both primitives and exponentials are only particular cases of
\emph{functions defined by ordinary differential equations with algebraic
coefficients} (e.g., $x(t)=\exp\int f(t)\,dt$ is a solution to the
differential equation $\dot x=f(t) x$). This looks as a most general
mechanism transforming algebraic input data into transcendental output
functions.

Thus we arrive to the following natural conclusion: in order to speak
about quasialgebraicity of the function fields $\C(f_1,\dots,f_\mu)$, the
generating functions $f_i(t)$ must be defined by polynomial ordinary
differential equations or systems of such equations. Then by computability
of any bound we would always assume that this bound can be expressed in
terms of the (real or complex) coefficients of these equations and/or
integer parameters (dimensions, degrees, etc), so that as functions of real
or complex parameters these bounds are polynomial, while being elementary
(or at worst primitive recursive) functions of the remaining integer
variables.

\begin{Ex}[principal]
Suppose that the functions $f_i(t)$, $i=1,\dots,\mu$, together satisfy a
system of polynomial ordinary differential equations of the form
\begin{equation*}
  \dot x_i=\sum_{|\alpha|\le m}c_{i\alpha} x^\alpha,\qquad
  i=1,\dots,\mu,\quad
  \alpha=(\alpha_1,\dots,\alpha_\mu)\in\mathbb Z_+^\mu,
\end{equation*}
with complex coefficients $c_{i\alpha}\in\C$. Then the expressions
\begin{equation*}
  C_1=(\max_{i,\alpha}|c_{i\alpha}|)^{N(\mu,m)},
  \qquad C_2=N(\mu,m)\cdot\sum_{i,\alpha}|c_{i\alpha}|
\end{equation*}
are both computable characteristics of the collection $\{f_1,\dots,f_\mu\}$
provided that $N$ is an elementary or at worst primitive recursive
function of two integer arguments $m$ (the degree) and $\mu$ (the
dimension of the system).
\end{Ex}

\subsection{Restricted computable bounds}
Explicit computability does not imply (neither does it assume) globality of
the bounds. In the same way as the parameters defining the functions $f_i$
may enter the answer, sometimes numeric characteristics of the domain (in
our settings, always a triangle $T\subset\C\ssm\S$) may enter into the
expressions.

\begin{Def}
Let $U\subset\C$ be a domain (usually polygonal) eventually containing
some of the singular points inside. We say that the field $\C(X)$ is
\emph{quasialgebraic in $U$}, if the number of zeros of any function from
this field can be bounded in any triangle $T\subset U\ssm\S$, uniformly on
all such triangles, but the bound may depend on $U$.
\end{Def}

The most commonly occurring form of this dependence is through the distance
between $\partial U$ and the singular locus $\S$.

\begin{Rem}
Since the independent variable ranges over the Riemann sphere $\C P^1$ and
the point $t=\infty$ may well belong to the singular locus, the distance
from $\partial U$ to $\S$ should be defined in such cases as minimum of the
above distance to the finite part of $\S$ and the number $\inf_{t\in
\partial U}|t^{-1}|$ measuring the ``distance from $\partial U$ to infinity''.
\end{Rem}

\section{Quasialgebraicity and uniform quasialgebraicity}
\label{sec:quasialgebraicity}

\subsection{Quasialgebraicity of Picard--Vessiot extensions: accurate
formulation of the problem} Let $X(t)$ be a fundamental matrix solution of
a Fuchsian system of $n$ linear ordinary differential equations
\eqref{fuchsian2}, and $\C(X)$ (resp., $\C[X]$) the field obtained by
adjoining all entries of this matrix to the field $\C$ (resp., the ring of
all polynomial combinations of these entries). As was already noted, we
always assume that $\C(X)$ contains the subfield of rational functions
$\C(t)$.

\begin{Def}
The field $\C(X)$ is called \emph{quasialgebraic}, if the number of
isolated roots of any function $f\in \C(X)$ of degree $k$ in this field, in
any triangle $T$ free from singular points $t_j$ of the Fuchsian system
\eqref{fuchsian2}, is bounded by a number depending only on:
\begin{enumerate}
  \item the degree $k=\deg_{\C(X)}f$;
  \item the dimension $n$ and the number $d$ of (finite) singular points
  (as an elementary or at worst primitive recursive function);
  \item the entries of the residue matrices $A_1,\dots,A_d$ (in an
  algebraic way);
  \item on the coordinates of the singularities $t_1,\dots,t_d$ (also in an
  algebraic way).
\end{enumerate}
\end{Def}

In fact, in all cases when quasialgebraicity of the field will be
established, the bounds would depend on the complex (matrix) parameters in
a very simple way, via the \emph{residual norm} of the rational matrix
function $A(t)$,
\begin{equation}\label{the-norm}
  R(A(\cdot))=\max_{j=1,\dots,d}\|A_j\|.
\end{equation}
In a similar manner, dependence of the bounds on the position of the
singularities $t_j$ will be expressed via the (inverse) \emph{spread} of
these points,
\begin{equation}\label{spread}
  \rho(\S)=\max_{i\ne j}\{|t_i-t_j|^{-1},\ |t_i|\}.
\end{equation}
(this number is large only when some of the singular points approach each
other).

\subsection{Uniform quasialgebraicity}
According to the above definition, the field is quasialgebraic if the
simple algebraic data (parameters) defining it can be used to produce an
explicit bound on the number of zeros of any function from this field.

Some of these parameters are always relevant. There is no question why the
dimension and the degree must necessarily enter any bound on zeros. It may
require some efforts to see that other parameters can also affect the
answer.

\begin{Ex}\label{ex:eigenvalues-do-matter}
The Euler system
\begin{equation*}
  \dot x_1=a\,t^{-1} x_1,\quad \dot x_2=0,\qquad a\in\mathbb N,
\end{equation*}
defines a field containing the function $f(t)=t^a-1$ that is of degree
\emph{one} in this field (though, as a polynomial, it has degree $a$). The
number of isolated zeros of the function $f$ in a triangle can be as large
as $\lfloor \tfrac a2\rfloor$.
\end{Ex}

This example shows that the magnitude of eigenvalues of the residue
matrices clearly affects the number of zeros.

It is not very difficult (see the subsequent sections) to construct upper
bounds (global or not) that would involve \emph{rational} expressions of
the parameters (in this case, entries of the residue matrices) that have
poles, making the bounds exploding for certain combinations of parameters.
On the other hand, there are no visible reasons for appearance of the
infinite number of zeros for these values of the parameters. It requires
considerable efforts to show that in terms of the residue matrices, the
bounds can be given by \emph{polynomial} expressions, which is equivalent
to expressing them in terms of the residual norm \eqref{the-norm} as
above. Finally, it is the most difficult part to prove that under certain
additional but rather natural assumptions, the bounds can be given
\emph{uniformly over all configurations of singular points}.

\subsection{Notes, remarks}
We conclude this highly informal discussion by several remarks, also of a
very general nature.

As stated in \secref{sec:AI-finiteness}, the ultimate goal of the theory
is to establish a constructive bound for the number of zeros of Abelian
integrals (tangential Hilbert problem). These integrals were shown to
belong to certain Picard--Vessiot extensions. Application of the methods
explained in subsequent sections requires repeated algebraic manipulations
with the integrals, therefore bringing into play the whole field of
functions. Moreover, on the final stage of the construction not one but
rather several Picard--Vessiot fields are considered simultaneously, while
carrying out induction in the number of ramification points. This explains
why the existence of upper bounds on the number of zeros
(quasialgebraicity) is defined as a property of the corresponding
functional fields (rather than elements constituting these fields).

There was also a special reason for choosing the definition of
quasialgebraicity without attempting to specify explicitly the ``counting
function'' measuring the number of isolated zeros. This was done primarily
because the bounds that one can obtain on this way are \emph{enormously}
excessive. Several nested inductive constructions immediately produce
tower-like bounds even from the modest exponential contributions on each
inductive step.

Finally, we would like to note that the Fuchsian representation is not the
only possible. Actually, the Picard--Fuchs system of differential equations
for Abelian integrals is written in the hypergeometric form
\begin{equation*}
  (tE+A)\dot X(t)=BX(t)
\end{equation*}
determined by two constant matrices $A,B\in\Mat_{n\times n}(\C)$, and it
would be natural to require that quasialgebraicity were expressed in terms
of the norms $\|A\|$, $\|B\|$ rather than in terms of the respective
residues after transforming to the Fuchsian form. This type of bounds is
not proved (if it were, this would imply constructive solution of the
tangential Hilbert problem, as explained in \secref{sec:redundant}).

\chapter{Quantitative theory of differential equations}

\section[Bounded meandering principle]
{Bounded meandering principle: explicit bounds on zeros of
functions defined by ordinary differential equations}\label{sec:meandering}

The preceding section contained motivations for introducing the notion of
quasialgebraicity. In this section we survey some results of constructive
(though not always global) nature, bounding the number of isolated zeros
of functions defined by differential equations.

\subsection{Linear $n^{\text{th}}$ order equations with bounded
analytic coefficients}\label{sec:alode-bounds} The basis for all other
considerations is a classical theorem by de la Vall\'ee Poussin concerning
solutions of linear ordinary differential equations with bounded
coefficients, not necessarily polynomial or even analytic. It gives a
sufficient condition for a linear ordinary differential equation of order
$n$ with \emph{real bounded coefficients} guaranteeing absence of
solutions with more than $n-1$ isolated roots on a given real interval
$I\subset\R$. Such equations are called \emph{disconjugate} on $I$. In
order to stress the difference with the complex case, we denote the
independent variable by $s$, and consider a linear equation
\begin{equation}\label{lode}
 y^{(n)}+a_1(s)\,y^{(n-1)}+ \cdots+ a_{n-2}(s)\,y''+a_{n-1}(s)y'+a_n(s)y=0
\end{equation}
on the real interval $I=[s_0,s_1]$ of length $r$ with real bounded
coefficients:  $|a_k(s)|<c_k<\infty$ for all $s\in I$.

\begin{Lem}[\cite{poussin}]\label{lem:lode-nonosc}
If
\begin{equation}\label{smallc}
\sum_{k=1}^n \frac{c_k\,r^k}{k!}<1,
\end{equation}
then any  $C^n$-smooth function $f(s)$ satisfying a linear equation
\eqref{lode} may have at most $n-1$ isolated roots on $I$, counted with
multiplicities.
\end{Lem}
This result can be seen as a generalization of the Sturm nonoscillation
theorem for equations of order greater than $2$. The inequality can be
slightly improved, see \cite{levin:non-osc}.

This (simple) statement implies a number of corollaries. First,
subdividing any interval into sufficiently short segments satisfying
\eqref{smallc}, and adding together the bounds for each segment, one can
obtain an explicit bound on zeros valid on any real interval where the
coefficients $a_i(s)$ are explicitly bounded. Next, one can consider
equations of the same form \eqref{lode} with \emph{complex-valued}
coefficients and solutions. There is no sense to count zeros of such
solutions (since in general there are no roots), instead an interesting
question is to estimate their \emph{topological index}, the variation of
argument $\Arg f(s)$ between the endpoints of the real interval. By a
simple modification of the method used in the proof of Lemma
\ref{lem:lode-nonosc} one can obtain an upper bound for this topological
index, again in terms of the magnitude (uniform upper bounds for the
absolute value) of coefficients of the equation \eqref{lode}, see
\cite{fields}, as follows.

\begin{Lem}\label{lem:index-sturm}
Variation of argument of any solution $f(t)$ of the equation \eqref{lode}
with complex-valued coefficients all bounded by $C$ in the absolute value,
along the segment $I$ of length $|I|=s_1-s_0$ is no greater than
\begin{equation*}
  \pi(n+1)(1+3C|I|).
\end{equation*}
\end{Lem}

If the equation has \emph{analytic} coefficients $a_i\in\mathcal O(U)$
defined in some domain $U\subset\C$, then the number of complex isolated
zeros inside any polygonal domain can be majorized using the above
described bounds for the topological index and the classical argument
principle. As a result, the following inequality for the number of
isolated zeros can be obtained (we return to the initial notation $t$ for
the complex independent variable).

\begin{Thm}[\cite{fields}]\label{thm:lin-eq}
Suppose $f_1(t),\dots,f_n(t)$ solve
\begin{equation}\label{lin-eq}
\begin{gathered}
  y^{(n)}+a_1(t)y^{(n-1)}+\cdots+a_{n-1}(t)y'+a_n(t)y=0,\qquad t\in U,
  \\
  a_i\in\mathcal O(U),\quad |a_i(t)|\le R,\ i=1,\dots,n.
\end{gathered}
\end{equation}
Then for any triangle $T\subset U$ of perimeter $\ell$ any linear
combination $f=\sum_1^n c_i f_i$ with complex constant coefficients $c_i$,
has no more than
\begin{equation}\label{eq-bound}
\begin{gathered}
\tfrac32(n+1)(1+\ell R)
\end{gathered}
\end{equation}
isolated zeros in $T$.
\end{Thm}

\begin{Rem}
Analyticity of solutions of the equation \eqref{lin-eq} follows from
analyticity of its coefficients. However, the theorem apparently can be
modified to cover also equations with meromorphic coefficients, provided
that the solutions are analytic in $U$ (the so called apparent
singularities).
\end{Rem}

\subsection{Systems of linear equations: Novikov's counterexample}
Motivated by potential applications to Abelian integrals and
quasialgebraic fields, we would like to find generalizations of
Lemma~\ref{lem:lode-nonosc} or Theorem~\ref{thm:lin-eq} for systems of
linear equations, in hope to majorize the number of isolated roots of any
linear combination of coordinates $c_1 x_1(s)+\cdots+c_n x_n(s)$ for a
system
\begin{equation*}
  \dot x(s)=A(s)x(s),\qquad s\in I\subset\R,\qquad\forall s\in
  I\quad \|A(s)\|\le R,
\end{equation*}
in terms of $\ell=|I|$ and $R$. The following counterexample due to
D.~Novikov \cite{mit:counterexample}, shows that it is impossible, neither
for real nor for complex systems.

\begin{Ex}
Let $t_1,...,t_d$ be collection of different numbers (real or complex)
from a real interval $I$ or a triangle $T$ (consider only the last case
for simplicity). Denote $a(t)=\e(t-t_1)...(t-t_d)$ and let $\e$ be so
small that $|a(t)|+|\dot a(t)+a^2(t)|<1$ for any $t\in T$.

Solution $\phi_1=\exp(\int a(t) dt)$ of the linear differential equation of
the first order $\dot x_1=a(t)x_1$, has no zeroes at all. However, the
derivative $\phi_2=\dot\phi_1=a(t)\phi_1$ has the same zeroes as $a(t)$
and satisfies the equation $\dot \phi_2=(\dot a+a^2)\phi_1$.

Together the pair $(\phi_1,\phi_2)$ satisfies the linear $2\times2$-system
\begin{equation}\label{mit-system}
  \dot x_1=ax_1,\quad\dot x_2=(\dot a+a^2)x_1,\qquad a(t)=\e\prod_1^d
  (t-t_j),
\end{equation}
whose coefficients are bounded by $1$ everywhere in $T$ by the choice of
$\e$. However, the second component has the specified number $d$ of
isolated zeroes there, where $d$ can be arbitrarily large.
\end{Ex}

This example suggests that when dealing with systems of linear equations,
instead of general analytic (variable) bounded coefficients, only rational
or even polynomial coefficients were allowed, so that
\begin{equation}\label{matr-poly}
  A(t)=\sum_{k=0}^d A_k t^k,\qquad A_k\in\Mat_{n\times n}(\C),
\end{equation}
and the magnitude of the matrix coefficients, e.g., the total norm
$R=\sum_k \|A_k\|$ used to produce upper bounds.

\subsection{Reduction of a system to an equation:
obstructions}\label{sec:linear-derivation} The above example is not very
surprising since the standard procedure of reducing a linear $n\times
n$-system
\begin{equation}\label{n-equations}
  \dot x=A(t)x,\qquad x=(x_1,\dots,x_n)\in\C^n
\end{equation}
to one $n$th order linear equation is discontinuous with respect to
parameters. Indeed, without loss of generality one can assume that it is
the number of zeros of the first component $y=x_1$ in a domain
$U\subset\C$ that interests us. Differentiating $x_1$ by virtue of the
system, we see that the derivatives are linear combinations:
\begin{equation}\label{derivatives}
  y^{(k)}(t)=\bb q_k(t)\cdot x(t),\qquad k=1,\dots,n,\dots,
\end{equation}
where the  analytic (co)vector functions $\bb
q_k(t)=(q_{k,1}(t),\dots,q_{k,n}(t))\in{\C^n}^*$ are determined by the
recurrent rule
\begin{equation}\label{iterations}
  \bb q_0(t)=(1,0,\dots,0),\quad
  \bb q_{k+1}(t)=\dot {\bb q}_k(t)+\bb q_k(t)\cdot A(t)
\end{equation}
and do not depend on the choice of the trajectory $x(t)$. (We use the
right multiplication to stress that $\bb q_k(t)$ are row vector functions).
Consider the field $\Bbbk=\mathcal M(U)$ of functions meromorphic in $U$
and the linear $n$-space $\Bbbk^n$ over this field.

The functions $\bb q_k(\cdot)$ are vectors in this space. Consider the
linear subspaces
\begin{equation}\label{chain-of-spaces}
  \{0\}\subset L_0\subseteq L_1\subseteq\cdots\subseteq
  L_{n-1}\subseteq L_n\subseteq\cdots\subseteq\Bbbk^n,
\end{equation}
each $L_k$ being spanned by $\bb q_0,\bb q_1,\dots,\bb q_k$ over $\Bbbk$.
If all inclusions are strict, the dimensions of $L_k$ over $\Bbbk$ must
strictly increase, while being bounded by $n=\dim_\Bbbk \Bbbk^n$. Hence the
ascending chain of subspaces \eqref{chain-of-spaces} must stabilize no
later than after $n$ steps, i.e., for some $\mu\le n$ the inclusion must be
non-strict, $\bb q_\mu\in L_{\mu-1}$, which implies a linear identity over
$\Bbbk=\mathcal M(U)$,
\begin{equation}\label{lin-comb}
  \bb q_\mu+a_{1}(t)\bb q_{\mu-1}+\cdots+a_{\mu-1}(t)\bb q_{1}+a_\mu(t)\bb q_0=0,
  \quad a_i\in\Bbbk=\mathcal M(U).
\end{equation}
This implies that $y=\bb q_0\cdot x$ satisfies the equation \eqref{lin-eq}.

Though the recurrent formulas generating the vectors $\bb q_1,\bb
q_2,\dots$ are explicit, it is impossible to control the magnitude of the
coefficients $a_i(\cdot)$, whatever that may mean, from this construction.
This can be seen in a simpler settings when the field $\Bbbk$ is replaced
by $\C$.

\begin{Ex}\label{ex:instable}
Consider a sequence of vectors $ q_0, q_1,\dots$, in the $n$-dimen\-sional
space over $\R$ or $\C$, that grows at most exponentially in the sense of
the norm: $\| q_k\|\le c^k$, $c>0$.

For the same dimensionality reasons as above, no later that at the $n$th
step a linear dependence must occur, allowing to express some $ q_k$ as a
linear combination of the preceding vectors $ q_0,\dots, q_{k-1}$,
\begin{equation}\label{lin-comb-C}
  - q_k+a_1 q_{k-1}+\cdots+a_k  q_0=0,\qquad a_i\in\C,
\end{equation}
similarly to \eqref{lin-comb}. However, coefficients of this dependence
are out of control and can be arbitrarily large. To see this, consider the
situation when the angle between $ q_{k-1}$ and $L_{k-2}$ is very small
but nonzero so that $L_{k-2}\subsetneq L_{k-1}$, while $ q_k$ belongs to
$L_{k-1}$ and is orthogonal to $L_{k-2}$ in $L_{k-1}$. Even assuming some
\emph{lower} bounds on $\| q_k\|$ will not help to improve the situation
in this case.
\end{Ex}

The situation becomes completely different under any of the following two
additional assumptions.

\begin{Ex}\label{ex:lattice}
Assume that all the vectors $ q_k$ actually belong to a lattice $\mathbb
Z^n\subset\C^n$ and their norms are bounded from above as in the preceding
example. Then the coefficients of the linear combination \eqref{lin-comb}
can be explicitly bounded from above: indeed, if $ q_k$ is a linear
combination of $ q_0,\dots, q_{k-1}$, then the coefficients $a_i$ in
\eqref{lin-comb-C} can be found by solving a system of linear
nonhomogeneous algebraic equations with integer matrix of coefficients and
integer free terms. By the Cramer rule, solutions of this system can be
obtained (after elimination of redundant equations and assuming that
solutions indeed exist) as ratios of appropriate minors (determinants of
some square submatrices). These minors are integer numbers, explicitly
bounded from above by virtue of the assumptions on the norms $\| q_i\|$.
Hence in each ratio the numerator is bounded from above, whereas the
denominator is no smaller than $1$ in the absolute value. This clearly
implies an upper bound on all the coefficients $a_i$.
\end{Ex}

\begin{Ex}\label{ex:iterated-linear-map}
Assume that the vectors $q_k\in \C^n$ are obtained by iterations of a
linear map $P\:\C^n\to\C^n$. Then the bounds on the norms will be
automatically satisfied with $c=\|P\|$. The \emph{first} linear
combination between the vectors still can have very large coefficients,
exactly as in Example~\ref{ex:instable}. However, instead of looking for
the first combination, we may continue until the step number $n$. Since
any operator $P$ is a matrix root of its characteristic polynomial
$\chi(z)=z^n+a_1 z^{n-1}+\cdots+a_{n-1}z+a_n$, the linear combination of
the form \eqref{lin-comb-C} in this case can be obtained via the identity
$\chi(P)q_0=0$. Coefficients $a_i$ of the characteristic polynomial admit
an upper bound in terms of the eigenvalues of $P$, each of which is no
greater than the norm $\|P\|$. Finally, we conclude that a linear
combination \eqref{lin-comb} in this case (iterations of a linear map $P$
in the finite-dimensional space $\C^n$)) can be constructed with
coefficients $a_i$ explicitly bounded in terms of $\|P\|$ and $n$.
\end{Ex}

\subsection{Reduction of a system to an equation: commutative algebra
versus linear algebra} The operator $\bb q(t)\mapsto \dot {\bb q}(t)+\bb
q(t)\cdot A(t)$ from the $n$-dimen\-sional $\C(t)$-linear space into
itself, describing the iterations \eqref{iterations}, is \emph{not} linear
over the field $\C(t)$, since the derivative $d/dt$ is not $\C(t)$-linear
operator. On the other hand, considered over the field $\C$, the space of
polynomial vector-functions is \emph{not} finite-dimensional unless $\deg
A(t)=0$ (in which case the degrees are uniformly bounded). Thus neither of
the above methods can work.

However, a solution can be found in terms of the commutative algebra. It
allows to treat in a similar way both linear and nonlinear systems. Before
proceeding further, we give a useful technical definition.

\begin{Def}
The \emph{norm} of a polynomial $p\in\C[x_1,\dots,x_n]$ is the sum of
absolute values of its coefficients:
\begin{equation}\label{def-norm}
  \left\|\sum\nolimits_\alpha c_\alpha x^\alpha\right\|=\sum\nolimits_\alpha|c_\alpha|.
\end{equation}
One can easily verify that in addition to the usual triangle inequality,
this norm is multiplicative, $\|pq\|\le\|p\|\,\|q\|$ for any two
polynomials. The norm of a derivative $\|\partial_i p\|$ (in any variable)
can be easily bounded, $\|\partial_i p\|\le \deg p\,\|p\|$.
\end{Def}

Generalization of the construction from \secref{sec:linear-derivation} is
very simple.

\begin{Ex}[basic]
Consider the linear system \eqref{n-equations} with a \emph{polynomial}
matrix of coefficients $A(t)$ as in \eqref{matr-poly}. The rule
\eqref{iterations} defines a sequence of \emph{polynomials} $q_k(t,x)=\bb
q_k(t)\cdot x\in\C[t,x]$, that are always linear in the variables
$x=(x_1,\dots,x_n)$: the recursive rules allow to estimate the degrees
$\deg_t q_k(t,x)$ and, if necessary, their norms.

Instead of the linear subspaces $L_k$ over the field $\Bbbk=\C(t)$ in this
case, consider the \emph{polynomial ideals}
$I_k=(q_0,\dots,q_k)\subset\C[t,x]$. The ascending chain
\begin{equation}\label{chainofideals}
  \{0\}\subset I_0\subseteq I_1\subseteq\cdots\subseteq
  I_n\subseteq\cdots\subseteq\C[t,x]
\end{equation}
of these ideals must eventually stabilize in the sense that some inclusion
$I_{\ell-1}\subset I_{\ell}$ becomes non-strict (an equality). This follows
from the fundamental fact that the polynomial ring $\C[t,x]$ is Noetherian
and \emph{any} ascending polynomial chain in it eventually stabilizes.

The above stabilization condition means that $\pm q_\ell\in I_{\ell-1}$,
hence for appropriate polynomial coefficients $h_1,\dots,h_\ell\in\C[t,x]$
\begin{equation}\label{inideal-general}
  q_{\ell}+\sum_{i=1}^{\ell}h_{i}q_{\ell-i}=0.
\end{equation}

\emph{A priori}, the polynomial coefficients $h_i$ can depend on $x$ in a
nontrivial way. However, since all polynomials $q_k$ are linear in $x$
(homogeneous), one can truncate the identity \eqref{inideal-general}
retaining only constant terms of $h_i$ (of degree $0$ in $x$) and
construct a new identity of exactly the same form \eqref{inideal-general}
but with coefficients $a_i(t)=h_i(t,0)$ from the \emph{univariate}
polynomial ring $\C[t]$,
\begin{equation}\label{inideal-linear}
  q_{\ell}+\sum_{i=1}^{\ell}a_{i}q_{\ell-i}=0,\qquad a_i\in\C[t],\
  i=1,\dots,\ell.
\end{equation}

The identity \eqref{inideal-linear} obviously means that the function
$q_0(t,x)=x_1$ satisfies the linear ordinary differential equation
\begin{equation}\label{derived-lode}
  y^{(\ell)}+\sum_{i=1}^{\ell}a_{i}(t)\,y^{(\ell-i)}=0,\qquad
  a_i(t)=h_i(t,0)\in\C[t],\ i=1,\dots,\ell.
\end{equation}

In order to apply the results form \secref{sec:alode-bounds}, one has to
estimate the absolute values $|a_i(t)|$ from above. In any disk of known
radius $\{|t|< R\}$ this can be done if the decomposition
\eqref{inideal-general} is explicitly known, in particular, if the
following parameters are explicitly bounded from above:
\begin{itemize}
  \item the length $\ell$ of the chain (equal to the order of the resulting
  differential equation);
  \item the degrees of the coefficients $h_i$ in $t$;
  \item the norms $\|h_i\|$, $i=1,\dots,\ell$.
\end{itemize}
\end{Ex}

\begin{Ex}\label{ex:chains}
The construction involving chains of ideals, is not very degree-specific.
Consider a system of \emph{polynomial} ordinary differential equations,
\begin{equation}\label{poly-sys}
  \dot x_i=P_i(t,x),\qquad i=1,\dots,n,\ P_i\in\C[t,x_1,\dots,x_n].
\end{equation}
Let $f$ be a polynomial combination $f(t)=Q(t,x_1(t),\dots,x_n(t))$ for
some $Q=Q(t,x)\in\C[t,x]$. Consider the infinite sequence of polynomials
$q_0$, $q_1$, \dots, $q_n,\dots\in\C[t,x]$ formed by iterations of the Lie
derivative, i.e., the recursive rule
\begin{equation}\label{lie-der}
  q_0=Q,\quad q_{k+1}=\pd{q_k}{t}+\sum_{i=1}^n\pd{q_k}{x_i}\cdot P_i.
\end{equation}

Let $I_k=(q_0,q_1,\dots,q_k)\subset\C[t,x]$ be the polynomial ideals
generated by the first $k+1$ polynomials from this sequence in the
polynomial ring $\C[t,x]$. They obviously form an ascending chain,
$I_k\subset I_{k+1}$. Since the ring $\C[t,x]$ is Noetherian, this chain
must stabilize and hence for some natural $\ell$,
\begin{equation}\label{inideal}
  q_\ell+h_1q_{\ell-1}+\cdots+h_{\ell-1}q_1+h_\ell q_0=0,\qquad
  h_i\in\C[t,x].
\end{equation}
As in the linear case, this identity implies a polynomial relationship
between the unknown function $y=Q(t,x(t))$ and its derivatives up to order
$\ell$. Unlike the linear case, this time the coefficients $h_i$ may
depend on $x$ explicitly, so that in addition to the data described in the
previous example, one needs upper bounds for $|x_i(t)|$ in $U$, which is a
non-algebraic piece of information. However, in many cases this
information can be easily achieved (or even \emph{a priori} known) and the
problem reduces to getting the same information as in the linear case,
namely: the length $\ell$ of the ascending chain of polynomial ideals, the
degrees and the norms of the polynomials $h_i$ appearing in the
representation \eqref{inideal}.
\end{Ex}

\subsection{Generalizations and improvements}
It is not clear from the very beginning, how replacing ascending chains of
linear subspaces by ascending chains of polynomial ideals may resolve the
problems related to unboundedness of the coefficients of the
decompositions \eqref{lin-comb} and \eqref{inideal} respectively, see
Exercise~\ref{ex:instable}. We explain it now.

The first advantage of the suggested approach allows to treat systems
depending polynomially on additional parameters,  without the risk of
producing bounds that blow up for certain values of the parameters.
Actually, the difference between parameters and phase variables disappears
almost completely.

\begin{Ex}\label{ex:parametric}
Assume that the linear system \eqref{linsys} with a polynomial matrix
$A(t)$ \eqref{matr-poly} depends on additional parameters
$\l=(\l_1,\dots,\l_p)$ is a polynomial way: $A\in\Mat_{n\times
n}(\C[t,\l])$. Then one can consider the chains of ideals in the bigger
ring $\C[t,\l]$ and only minor notation changes are necessary to construct
a linear ordinary differential equation \eqref{derived-lode} whose
coefficients will be in fact polynomial in $t$ \emph{and} $\l$ (note that
the resulting differential equation is always monic: its leading
coefficient before the principal derivative is $1$).

This polynomiality eliminates the danger that for some values of the
parameters the coefficients of the derived equation \eqref{derived-lode}
will blow up (which was earlier the case).
\end{Ex}

Another advantage appears as a generalization of Example~\ref{ex:lattice}.

\begin{Ex}\label{ex:lattice-ideal}
Assume that the coefficients matrix $A(t)$ of the system
\eqref{n-equations} is polynomial as in \eqref{matr-poly}, and in addition
all matrix coefficients $A_k$ have only integer entries:
$A_k\in\Mat_{n\times n}(\mathbb Z)$. Then, since the polynomial
$q_0(t,x)=x_1$ also belongs to the subring $\mathbb Z[t,x]$, all
subsequent polynomials $q_k$ will also have integer coefficients, and
their degrees are growing not faster than linearly in $k$. The growth of
the norms $\|q_k\|$ can also be easily controlled.

Suppose that the length $\ell$ of the ascending chain of the corresponding
ideals \eqref{chainofideals} is already known (in the univariate case it is
relatively simple, see \cite{annalif-99}). Then one can explicitly compute
an upper bound $r$ for the degrees of polynomial coefficients $h_i$ in
\eqref{inideal-general} in terms of $\ell$, the degree $d=\deg A(t)$ and
$n=\dim x$. To find the polynomials $h_i$ themselves (or rather the
univariate polynomials $a_i=h_i(\cdot,0)$), it is possible now to use the
method of indeterminate coefficients: writing each $a_i$ as $\sum_{j=0}^r
c_{ij}t^j$ and substituting them into the identity \eqref{inideal-linear},
we obtain a system of non-homogeneous linear algebraic equations for the
unknown variables $\{c_{ij}\:i=1,\dots,\ell, j=0,\dots,r\}$, with integral
coefficients matrix and the free terms column, all explicitly bounded from
below. For the same reason as in Example~\ref{ex:lattice}, in this
``lattice'' case explicit bounds for $\|a_i\|$ can be immediately produced
in terms of $n$, $d$ and $r$.
\end{Ex}

The lattice polynomial system may look artificial, but in fact their
appearance is natural. The explanation is given in the following principal
example.

\begin{Ex}[universal system]\label{ex:universal}
All entries of the matrix coefficients $A_k$ of the polynomial matrix
function \eqref{matr-poly} can be considered as parameters and denoted by
$\l_i$, $i=1,\dots,(d+1)n^2$. With respect to these variables, the
``universal matrix polynomial'' $\sum_{k=0}^d A_k t^k$ is in fact in the
``lattice'' $\Mat_{n\times n}(\mathbb Z[t,\l])$, of known degree $d+1$:
moreover, all coefficients of $A(t)$ over the ring $\mathbb Z[t,\l]$ are
zero-one matrices.

In a similar way one can start iterations from the ``general linear form''
$q_0(t,x)=\sum_{i=1}^n \beta_i x_i$ which becomes a polynomial of degree
$n+1$ with zero-one coefficients in $\C[x,\beta_1,\dots,\beta_n]$ and add
the string of the coefficients $\{\beta_i\}$ to the parameter list. This
would allow to treat simultaneously isolated zeros of all nontrivial linear
combinations.
\end{Ex}

\subsection{Lengths of ascending chains}
Example~\ref{ex:chains} in fact proves that for any dimension $n$ and
degree $d$ of the matrix polynomial $A(t)$ from \eqref{matr-poly}, one can
majorize explicitly the number of isolated zeros of any linear combination
$\sum \beta_i x_i(t)$ in terms of $n,d$, and $R=\sum_{0}^d\|A_k\|$
provided that the length of ascending chain of polynomial ideals
\eqref{chainofideals} is explicitly known. We claim that this length is a
``computable'' function.

Notice that for each combination of $n,d$ we have a uniquely defined
parameter space $\{\l\}=\C^{(d+1)n^2+n}$ (coefficients of the matrix
polynomial $A(t)$ and the initial linear form $q_0$), and the uniquely
defined chain of ideals in the polynomial ring $\C[t,x,\l]$ generated by
polynomials $q_k$ that in fact belong to the lattice $\mathbb Z[t,x,\l]$.

The construction of the generators $q_k$ is absolutely explicit. The
problem of verifying whether the next polynomial $q_k$ belongs to the
ideal $I_{k-1}$ generated by the previous polynomials, is constructive:
there exists an algorithmic procedure allowing to get a positive or
negative answer in a finite number of steps for each $k$. Thus the length
$\ell$ is a ``computable'' (in the weakest sense) function of $n$ and $d$:
the above description can be transformed into an algorithm computing the
value $\ell(n,d)$ for any given combination of $n$ and $d$.

Coupled with the preceding discussion, this computability means that the
problem on the number of zeros of functions defined by polynomial ODE's,
is algorithmically solvable. However, the complexity of this algorithm
turns out to be very high in the sense that the bound on the growth rate
of the function $\ell(n,d)$ is tremendous and higher than any elementary
function. In the next section we address this question in details and show
what other modifications are necessary in order to produce a
``theoretically feasible'' upper bound.

\section{Lengths of chains}\label{sec:chains}

\subsection{Descending chains of algebraic varieties}
The problem on ascending chains of polynomial ideals, belongs entirely to
the realm of commutative algebra. However, there is a parallel geometric
problem that admits a considerably more transparent solution. Moreover,
the construction used in the ``geometric'' proof, can be adjusted (after
introducing appropriate technical tools) to the algebraic case.

Consider a \emph{strictly} decreasing chain of complex algebraic varieties,
\begin{equation}\label{chainofvarieties}
  \C^n=X_0\supset X_1\supset X_2\supset\cdots\supset X_k\supset\cdots
\end{equation}
where each variety $X_k$ is given by a finite number of polynomial
equations in $\C^n$ of degree no greater than $d_k$ (no restrictions are
imposed on the number or the structure of these equations). Without loss
of generality we may assume that the sequence $d_1,d_2,\dots$ is
non-decreasing.

By the same Noetherian property, the chain \eqref{chainofvarieties} must
terminate after finitely many steps at some variety $X_\ell$. The problem
is to compute the length $\ell$ from the dimension $n$ and knowing the
bounds for the degrees $d_k$.

If the degrees $d_k$ are all bounded, then the problem belongs to linear
algebra, since the polynomial equations defining the varieties $X_k$ will
in fact constitute a finite-dimensional linear space, and its dimension
will be the natural bound for the length of the chain
\eqref{chainofvarieties}. In order to avoid technical troubles when talking
about computable dependence on \emph{infinite} input data
$\{d_1,d_2,\dots,\}$, we will consider \emph{finite-parameter} examples.
The most important are the linear growth case, when
\begin{equation}\label{d+k}
  d_k=d+k,\qquad k=1,2,\dots
\end{equation}
or the exponential growth case
\begin{equation}\label{d-to-k}
  d_k=d^k,\qquad k=1,2,\dots.
\end{equation}
In both cases the natural number $d$ is a parameter, and in both cases
$\ell$ should be majorized in terms of $n$ and $d$.

\subsection{Lexicographically decreasing sequences of words}
The advantage of the ``geometric'' problem on chains of complex algebraic
varieties, is that it can be reduced to a purely combinatorial problem on
lexicographically decreasing sequences of words.

Recall that any algebraic variety can be uniquely represented as the union
of irreducible algebraic subvarieties of different dimensions varying from
$0$ (isolated points) to $n-1$ (hypersurfaces). Thus with any
$X\subset\C^n$ one can associate a vector $\nu(X)\in\mathbb Z_+^n$ with
$n$ integer nonnegative coordinates $(\nu^{n-1},\dots,\nu^1,\nu^0)$, where
$\nu^i=\nu^i(X)$ stands for the \emph{number} of irreducible
$i$-dimensional components of $X$.

Denote by $\prec$ the lexicographic order on $\mathbb Z_+^n$, letting
\begin{equation*}
  (\nu^{n-1},\dots,\nu^0)\prec(\bar\nu^{n-1},\dots,\bar\nu^0)
\end{equation*}
if and only if for some $k=1,\dots,n$
\begin{equation*}
  \nu^{n-1}=\bar\nu^{n-1},\ \dots,\ \nu^k=\bar\nu^k,\ \text{but}\
  \nu^{k-1}<\bar\nu^{k-1}.
\end{equation*}

The following elementary observation is crucial.

\begin{Lem}
If $X\subsetneq Y\subset\C^n$, then $\nu(X)\prec\nu(Y)$.
\end{Lem}

\begin{proof}
This holds since:
\begin{enumerate}
 \item any irreducible component of $X$ should belong to an irreducible
 component of $Y$, and
 \item if $A\subset B$ is a pair of \emph{irreducible} varieties, then
 $\dim A\le\dim B$ and in the case of equal dimensions necessarily $A=B$.
\end{enumerate}
In other words, when passing from $Y$ to $X$ each irreducible component
either completely survives, or is split into a number of other components
of lower dimensions.
\end{proof}

The fact that a descending chain of algebraic varieties must stabilize,
follows now from the following purely combinatorial claim.

\begin{Prop}
A lexicographically strictly decreasing chain
\begin{equation}\label{lexic-chain}
  \nu_1\succ\nu_2\succ\cdots\succ\nu_k\succ\cdots,\qquad\nu_k\in\mathbb
  Z_+^n,
\end{equation}
must be finite.
\end{Prop}

\begin{proof}
For $n=1$ the claim is obvious, since any decreasing sequence of
nonnegative integers must be finite. For an arbitrary $n$, the first
``letters'' $\nu^{n-1}_k\in\mathbb Z_+$, $k=1,2,\dots$, form a
non-increasing sequence (it must \emph{not} necessarily be strictly
decreasing). However, no more than a finite number of values is taken.
Along any interval of constancy of the first letter, the tails
$(\nu^{n-2}_k,\dots,\nu^0_k)$ also form a lexicographically strictly
decreasing sequence in $\mathbb Z^{n-1}_+$. Hence the length of each such
segment is finite by the induction assumption, and the length of the whole
chain is finite as the sum of finitely many finite numbers.
\end{proof}

This proof, being extremely simple, can be supplied by quantitative
estimates of the lengths under additional assumptions on the \emph{norms}
of the words. Denote $\|\nu\|=\nu^{n-1}+\cdots+\nu^0$ the norm on $\mathbb
Z^n_+$.

\begin{Ex}\label{ex:appearance-Ack}
Assume that the norms of the words $\nu_k$ forming the chain
\eqref{lexic-chain} are bounded:
\begin{equation}\label{lexic-norm-bounds}
  \|\nu_k\|\le d+k,\qquad\forall k=1,2,\dots \ \text{($d$ a natural
  parameter)}.
\end{equation}
Then the length of the chain as a function of $n$ and $d$ is bounded by a
general recursive (but \emph{not} primitive recursive) function. The
explanation is as follows.

Let $f(n,d,i)$ be the maximal length of the chain under the additional
constraint that the first ``letter'' of the first word is no greater than
$i$,  $\nu_1^{n-1}\le i$. Then the restricted function $f(n+1,\cdot,i+1)$
can be expressed via $f(n,\cdot,\cdot)$ and $f(n+1,\cdots,i)$.

Indeed, the length of the \emph{initial} segment of the chain, on which
the first letter maintains its initial value $i$, can be at most
$f(n,d,d)$, since the tails start with a word whose first letter can be at
most $d$. After the first segment is exhausted, the remaining part of the
chain begins with a word of length $n$, whose first letter is at most
$i-1$, and with the restriction on the norms of the words as follows,
\begin{equation*}
  \|\nu_k\|\le d+(k+f(n,d,d)),
\end{equation*}
if the words of the remaining chain are numbered after the drop of the
first letter. This is tantamount to replacing $d$ by $f(n,d,d)$, therefore
we obtain the recurrent inequality
\begin{equation}\label{ackermann-for-words}
  f(n+1,d,i)\le f(n,d,d)+f(n+1,d+f(n,d,d),i-1).
\end{equation}
Coupled with the boundary conditions
\begin{equation*}
  f(n,d,0)\le f(n-1,d,d),\qquad f(1,d,i)\le i,
\end{equation*}
this determines the upper bound for $f$ completely. Notice the remarkable
similarity between the recurrent rule \eqref{ackermann-for-words} for $f$
and that for the Ackermann generalized exponential \eqref{ackermann}.
\end{Ex}

The growth rate determined by the above recursive function, is enormous:
in terms of $n$, the dimension of the ambient space, the function $f$
grows faster than any primitive function. On the other hand, this bound is
essentially sharp: for lexicographically ordered chains it is rather
obvious, for chains of polynomial ideals it was proved by G.~Moreno
\cite{moreno,moreno:conference}, for chains of algebraic varieties the
claim is apparently also true though not written anywhere. In any case the
bounds based on such estimate of the length of chains should not be
considered as constructive, even theoretically.

However, under rather mild additional assumptions the bounds can be
improved very considerably, in fact to become elementary functions of
$n,d$.

\subsection{Dynamically generated chains}
The chain of ideals \eqref{chainofideals}  is generated by a dynamical
system. More precisely, the rule \eqref{lie-der} for the general case of a
polynomial vector field \eqref{poly-sys} is a Lie derivative of the ring
$\C[t,x]$, and the ideals $I_k$ are generated by iterated Lie derivatives
of the initial polynomial $q_0\in\C[t,x]$.

A parallel construction for chains of algebraic varieties also involves a
dynamical ingredient \cite{lleida,annalif-99}. Namely, instead of the
general chains \eqref{chainofvarieties} with the only restriction on the
degrees of the varieties $X_k$, one should consider chains generated by a
discrete time dynamical system in $\C^n$.

Let $F\:\C^n\to\C^n$ be a polynomial map of some known degree $d$ and
$X\subset\C^n$ an algebraic subvariety. Consider the decreasing chain
$\{X_k\}$ defined by the recursive rule
\begin{equation}\label{dynamic-chain-X}
  X_0=X,\qquad X_{k+1}=X\cap F^{-1}(X_k),\quad k=0,1,2,\dots.
\end{equation}
In other words,
\begin{equation*}
  X_k=X\cap F^{-1}(X)\cap F^{-2}(X)\cap\cdots\cap F^{-k}(X),
\end{equation*}
where $F^k$ stands for the $k$ times iterated map $F\circ\cdots\circ F$
and $F^{-k}(\cdot)=(F^k)^{-1}(\cdot)$ denotes the corresponding preimage.

Dynamically $X_k$ can be described as the set of points $a\in\C^n$ such
that the point $a$ itself, together with its $F$-orbit $F(a),\dots,F^k(a)$
of length $k$, belong to $X$. This immediately implies that
\begin{equation}\label{convexity-X}
  F(X_k\ssm X_{k+1})\subseteq X_{k-1}\ssm X_k,\qquad k=1,2,\dots
\end{equation}
(the differences occurring above consist of initial conditions of orbits
that jump off the variety $X$ \emph{exactly} after the specified number of
steps).

The condition \eqref{convexity-X} means that the difference $X_k\ssm
X_{k+1}$ cannot be too large compared with the preceding difference
$X_{k-1}\ssm X_k$. For an arbitrary polynomial map this is not true, but
under the additional assumption that $F$ preserves dimensions of
semialgebraic sets (i.e., takes curves into curves and not into points,
though eventually creating singularities, and the same in higher
dimensions), one can conclude that the dimensions of the differences
$X_k\ssm X_{k+1}$ are non-increasing:
\begin{equation}\label{dimen-nonincrease}
  \dim (X_k\ssm X_{k+1})\le \dim (X_{k-1}\ssm X_k),\qquad k=1,2,\dots.
\end{equation}

\subsection{Dynamically generated chains stabilize fast}
Consider a strictly descending chain of varieties \eqref{chainofvarieties}
satisfying the additional condition \eqref{dimen-nonincrease} on the
dimensions. Then the associated words $\nu_k=\nu(X_k)\in\mathbb Z_+^n$, in
addition to the lexicographic decrease \eqref{lexic-chain}, display
stronger monotonicity properties.

Consider the sequence of ``heads''  (as opposed to ``tails''), obtained by
truncating the words $\nu_k$ to their first $s$ symbols,
$[\nu_k]_s=(\nu_k^{n-1},\cdots,\nu_k^{n-s})\in\mathbb Z_+^s$, for all
$k=1,2,\dots$. For any lexicographically non-increasing sequence, the
sequence of heads of any fixed length $s$ will be again non-increasing (by
definition of the lexicographic order). Yet this monotonicity can be
non-strict in general.

However, under the additional assumption \eqref{dimen-nonincrease} any
truncated sequence of heads $\{[\nu_k]_s\}_{k=1}^\infty$ for the sequence
$\nu_k=\nu(X_k)$ will be \emph{strictly decreasing} in the following sense:
\begin{equation}\label{ultra-lex}
 \forall s=1,\dots,n-1\qquad
  [\nu_1]_s\succ[\nu_2]_s\succ\cdots\succ[\nu_{k-1}]_s
  =[\nu_k]_s=\cdots
\end{equation}
In other words, as soon as the first equality between truncated words
occurs, the rest of the chain will have the same heads. Indeed,
$[\nu_{k-1}]_s=[\nu_{k}]_s$ if and only if all irreducible components of
$X_{k-1}$ and $X_k$ of dimensions $n-s$ and higher, are the same, which
means that the difference $X_{k-1}\ssm X_k$ is at most
$(n-s-1)$-dimensional.

Such type of descent ensures much faster convergence. Indeed, the first
letter stabilizes after no more than $\|\nu_1\|$ steps, after which the
problem is reduced to that for the tails, which are words of length $n-1$.
The corresponding inductive inequality is very simple.

\begin{Ex}\label{ex:d+k}
Denote by $g(n,d)$ the maximal length of decreasing sequence of words of
length $n$ growing no faster than linear in the sense of the norm,
$\|\nu_k\|\le d+k$, with $d$ being an integer parameter as in \eqref{d+k},
under the additional assumption \eqref{ultra-lex}. Then the above
observation implies that the length until stabilization of the first
letter is at most $d$, whereas the sequence of tails of length $n-1$
starts from the word of norm no greater than $d+d=2d$ and hence its length
by the induction assumption is no greater than $g(n-1,2d)$, so that finally
\begin{equation}\label{inductive-bound-d+k}
  g(n,d)\le d+g(n-1,2d),\qquad g(1,d)\le d.
\end{equation}
This gives an upper bound for the length,
\begin{equation*}
  g(n,d)\le (2^{n}-1)d,
\end{equation*}
which is much better than in the general case considered in
Example~\ref{ex:appearance-Ack}.
\end{Ex}

\begin{Ex}\label{ex:d-to-k}
In a similar way the exponential bound \eqref{d-to-k} can be treated. In
this case the inductive inequality analogous to \eqref{inductive-bound-d+k}
takes the form
\begin{equation}\label{inductive-bound-d-to-k}
  g(n,d)\le d+g(n-1,d^d),
\end{equation}
which gives $g(n,d)$ as a tower function of height $n$. Though not
elementary, this is obviously a primitive recursive function of both
arguments.
\end{Ex}

These simple examples illustrate the algorithmic complexity of the problem
on lengths of lexicographically descending chains. For descending chains of
algebraic varieties the answer follows in fact from
Example~\ref{ex:d-to-k}, since the degrees of equations defining the
algebraic varieties $X_k$ grow exponentially: $\deg F^k=d^k$, where
$d=\deg F$ (without loss of generality we may assume that $\deg X=d$ with
the same $d$). The bounds for the norms $\|\nu(X_k)\|$ follow from a
version of B\'ezout theorem due to J.~Heintz \cite{joos:bezout}.

The most technically difficult case is that of ascending chains of
polynomial ideals, mainly since there is no uniqueness in the primary
decomposition of such ideals, hence one cannot associate a word $\nu(I)$
with an ideal $I$, counting the number of primary components of various
dimensions. However, the components of the maximal dimensions are
correctly defined and may be counted, which allows to implement a similar
inductive proof. Details can be found in \cite{annalif-99}.

\section{Restricted quasialgebraicity of Picard--Vessiot fields}

Quasialgebraicity of a function field was defined as a property allowing
for counting zeros globally. However, for technical reasons we need a
weaker notion of \emph{restricted quasialgebraicity}. Assuming $U\subset\C$
being a simply connected (usually polygonal or circular) domain containing
no singular points on its boundary (but eventually some singularities
\emph{inside} $U$), we can restrict functions from $\C(X)$ on $U$. The
result, denoted by $\C_U(X)$, is another functional field, again consisting
of multivalued functions.

\begin{Def}
We say that $\C(X)$ is \emph{restricted quasialgebraic in $U$}, if the
upper bounds on the number of isolated zeros in any triangle $T\subset
U\ssm\S$ can be given in the same terms as before, plus eventually some
geometric parameters describing the relative position of $\S$ and $U$,
most often the distance between $\S$ and the boundary $\partial U$.
\end{Def}

\subsection{Bounds in the disk for polynomial systems}
The above discussion allows to analyze completely the polynomial case.

\begin{Thm}\label{thm:meandering}
If the coefficients matrix $A(t)$ of a linear system \eqref{linsys} is
polynomial as in \eqref{matr-poly} with bounded matrix coefficients of
known total norm $R$, then the corresponding Picard--Vessiot extension
$\C(X)$ will be quasialgebraic after restriction on any disk $D_r\subset\C$
of radius $r$.

In other words, the bounds on the number of zeros of functions from
$\C_{D_r}(X)$ will be explicit but depending on $r$. As $r\to+\infty$, the
bounds explode.
\end{Thm}

Explosion of the bounds occurs for a very simple reason: the singular point
at infinity is in general an irregular (and certainly non-Fuchsian)
singularity which may be an accumulation point for isolated roots of
polynomial combinations (see \secref{sec:insearch}).

\subsection{Fuchsian system: reduction to the polynomial case}
To apply the previous technique to a Fuchsian system \eqref{fuchsian2}, it
is sufficient to introduce the new independent (complex time) variable.

The matrix function $A(t)$ from \eqref{fuchsian2} can be reduced to the
common denominator,
\begin{equation*}
  A(t)=\tfrac1{\chi(t)}P(t),\qquad P=\sum_{i=0}^{d-1} P_i t^i,
  \qquad \chi(t)=\prod_{i=1}^d (t-t_i).
\end{equation*}
The corresponding linear system is \emph{orbitally} equivalent to the
\emph{polynomial} system
\begin{equation}\label{fuch->poly}
  \left\{
  \begin{aligned}
   \dot X&=P(t)X,
   \\
   \dot t&=\chi(t),
  \end{aligned}\right.
  \qquad \cdot=\frac{d}{d\tau},
\end{equation}
where $\tau\in\C$ is the new time variable. The map
\begin{equation*}
  t\mapsto\tau(t)=\int_0^t\frac{dz}{\chi(z)}
\end{equation*}
is defined on the universal covering surface of $\C\ssm\S$ and takes
explicitly bounded values away from $\S$. The inverse map $\tau\mapsto
t(\tau)$ covers the complement $\C\ssm\S$ so that for any triangle
$T\Subset\C\ssm\S$ on distance $\e>0$ from $\S$ one can find $\rho$ such
that $T$ is covered by the image of the disk $\{|\tau|<\rho\}\subset\C$.

Application of Theorem~\ref{thm:meandering} to the system
\eqref{fuch->poly} implies the following corollary (after some preliminary
work).

\begin{Cor}\label{cor:quasialg-nonsing}
The field $\C(X)$ constructed from the Fuchsian system \eqref{fuchsian2},
is quasialgebraic in any domain containing no singular points  from $\S$.
The bounds would depend on $r=r(U)=\dist(\partial U,\S)=\dist(\overline
U,\S)$ and explode as $r\to 0^+$.
\end{Cor}

\subsection{Rational systems with apparent singularities}\label{sec:apparent}
 Actually, not all singular points are dangerous. As was already remarked,
Theorem~\ref{thm:lin-eq} is not about the number of zeros of solutions
inside a triangle $T$, but rather about the variation of argument of these
solutions along the boundary of $T$.

This means that the above theorems on restricted quasialgebraicity of
Fuchsian systems away from their poles, can be reformulated as assertions
on computability of upper bounds for the variation of argument along
polygonal paths distant from the singular locus $\S$.

If the domain $U$ has only \emph{apparent} singularities inside, then all
solutions of the system \eqref{linsys} are meromorphic in $U$ and  the
order of their poles can be explicitly bounded in terms of the
corresponding residue norms $\|A_j\|$. By the argument principle, this
together with the bounds on the index along the boundary, implies
restricted quasialgebraicity of $\C_U(X)$.

This remark shows that it is the multivaluedness of solutions rather than
any other circumstance, that is an obstruction for the global
quasialgebraicity. In the next sections we show how multivalued functions
can be treated, first locally and then globally.

\chapter{Isomonodromic reduction principle and Riemann--Hilbert problem}

\section{Isomonodromic fields: local theory}

\subsection{Euler field: an example}
The Euler field \eqref{euler-field}, obtained by adjoining  to $\C$ all
entries of the (multivalued) matrix function $t^A=\exp (A\ln t)$ solving
the system \eqref{euler}, already occurred as a first example when
quasialgebraicity could be expected.

The necessary condition for quasialgebraicity is given in terms of the
spectrum of the constant matrix $A$ (the only residue of the system). We
have already seen that if some of the eigenvalues are non-real, then
infinitely many real zeros may easily accumulate to the origin. It can be
relatively easily seen that \emph{under} this assumption on the spectrum
of $A$, such accumulation is impossible not only for real, but also for
complex zeros. Actually, a stronger assertion holds: the \emph{number} of
isolated roots can be explicitly majorized.

The following result (together with explicit bounds that in this case are
rather accurate) can be found in \cite{jdcs-96a}. Consider a finite
\emph{set of exponents} $S\subset\C$ of diameter $\diam S=\max_{\l,\l'\in
S}|\l-\l'|$ whose points may have non-trivial multiplicities
$\nu(\l)\in\mathbb N$, so that by definition $\sum_{\l\in S}\nu(\l)=\#S$.
\begin{Thm}\label{thm:complex-rolle}
If $S\subset\R$, then the number of isolated roots of any finite sum
\begin{equation}\label{quasipolynomial}
  \sum_{\l,k}c_{k\l}\,t^\l\ln^{k-1}t,\qquad \l\in S,\ k\le\nu(\l),
  \ c_{k\l}\in\C,
\end{equation}
in any triangle $T\subset\C\ssm\{0\}$ is no greater than $\#S-1+2\diam
S$.\qed
\end{Thm}

This theorem implies quasialgebraicity of the Euler field, since any
element $f$ from $\C[t^A]$ can be represented in the form
\eqref{quasipolynomial} with explicit control over $\diam S$ and $\#S$
expressed in terms of $\|A\|$, $\dim A$ and $\deg f$.

The assumption of this theorem can be formulated in terms of the spectrum
of the (only) monodromy operator.

\begin{Cor}\label{cor:euler}
If the spectrum of the monodromy operator $M$ corresponding to a small
loop around the origin, belongs to the unit circle \(i.e., all eigenvalues
have modulus $1$\), then the Euler field $\C(t^A)$ is globally
quasialgebraic, and the number of isolated zeros can be bounded in terms
of the norm $\|A\|$.
\end{Cor}

\subsection{Restricted quasialgebraicity near a Fuchsian point}
Consider a Fuchsian system \eqref{fuchsian2} with a singular point $t_1=0$
at the origin and all residues of norm $\le R$. Let $\rho$ be a
(sufficiently small) positive number such that all other singularities are
at least $2\rho$-distant and at most $1/\rho$-distant from the origin.

Let $M$ be a monodromy matrix corresponding to a small positively oriented
loop around the origin (this matrix is defined up to a conjugacy) and
$S\subset\C$ the spectrum of $M$ (i.e., the collection of eigenvalues which
is independent of anything but the system and the singular point).

\begin{Thm}\label{thm:quasialg-local}
If $S\subset\{|\l|=1\}$, then the field $\C(X)$ is quasialgebraic in the
disk $D_\rho=\{|t|<\rho\}$. The bound on the number of isolated zeros can
be given in terms of $R$ and $1/\rho$.
\end{Thm}

The outline of the proof of this theorem occupies the rest of this section.
A similar result was proved in \cite{mrl-96} for \emph{linear equations} of
order $n$ near a Fuchsian singular point.

\subsection{Joint fields}
Consider a Fuchsian singular point $t_1=0$ of the system \eqref{fuchsian2}
with the corresponding residue matrix $A_1=B$. Then the fundamental matrix
solution $X(t)$ in any disk $D$ containing no other singularities, can be
represented as
\begin{equation}\label{local-repr}
  X(t)=Y(t)\,t^B, \qquad t\in D,
\end{equation}
where $Y$ is a meromorphic (single-valued) function of $t$ in $D$.
Replacing if necessary $B$ by $B+r E$ with an appropriate $r\in\mathbb Z$,
one can always consider $Y(t)$ as a holomorphic function at $t=0$.
Expanding all matrix products, we conclude that the field $\C(X)$ belongs
to a bigger field $\C(Z,Y)$ spanned jointly by entries of the two matrix
functions, $Z(t)=t^B$ and $Y(t)$. Whereas the field $\C(Z)$ is a
Picard--Vessiot extension for the Euler system, this is not immediately
obvious concerning the extension field $\C(Y)$. However, we claim that the
field $\C(Y)$ is a subfield of a bigger Picard--Vessiot field for some
other Fuchsian system having $t=0$ as an apparent singularity. Indeed,
\begin{equation}\label{sys-for-H}
  \dot Y=\dot X\,t^{-B}-t^{-1}Xt^{-B}\,B=A(t)Y-t^{-1}YB.
\end{equation}
This is \emph{not} a linear system for $Y$ as a matrix function to satisfy
(since it involves the matrix multiplication from both sides). But if all
entries of $Y$ are arranged as one column vector, then \eqref{sys-for-H}
becomes a system of $n^2$ linear ordinary differential equations with
rational coefficients, exhibiting a Fuchsian singularity at $t=0$
(actually, at all other poles of $A(t)$ as well). The residues of this
larger-size system can be explicitly constructed and their norms bounded
from above.

The joint field $\C(Z,Y)$ possesses a property that was already used for
``single'' extension fields.

\begin{Lem}\label{lem:joint-argument}
Variation of argument of any element from the joint field $\C(X,Y)$ along a
polygonal path distant from the union of singular loci, is a computable
function. The same is true also for sufficiently small circular arcs
around each singular point.
\end{Lem}

Note that $t=0$ by construction is an apparent singularity for this
system. Thus the Picard--Vessiot field $\C(Y)$ is quasialgebraic restricted
on $D$, as follows from \secref{sec:apparent}.

The problem that we face now, is to prove the \emph{restricted}
quasialgebraicity in $D$ of the joint extension field $\C_D(Z,Y)$
containing $\C_D(X)$ by virtue of \eqref{local-repr}, having already
established restricted quasialgebraicity of each of the fields $\C_D(Z)$
and $\C_D(Y)$ separately and the fact that $Y$ has only one apparent
singularity in $D$. This will be done by reducing the question on
(restricted) quasialgebraicity of $\C_D(Z,Y)$ to that for $\C(Z)$, using
the fact that $Y$ has a trivial monodromy in $D$. This reduction can be
extracted from the papers by Petrov
\cite{petrov:cubic-nonosc,petrov:trick}.

First we need some real analysis.

\subsection{Real closedness}
Real part of an analytic function $f$ on an open domain $U$ is not analytic
(unless it is a constant). Yet the real part of the restriction $f|_\gamma$
on an analytic curve $\gamma\subset U$  can be extended to a neighborhood
of $\gamma$ and sometimes to the whole of $U$ as an analytic function (of
course, taking non-real values outside of $\gamma$).

If $U=D$ is a disk centered at the origin and $\gamma$ the real line $\R$,
then for any function $f$ meromorphic in $D$ one can take
\begin{equation}\label{re-im}
  \Re f(t)=\tfrac12(f(t)+\overline{f(\bar t)}),
  \qquad
  \Im f(t)=\tfrac1{2i}(f(t)-\overline{f(\bar t)}),
\end{equation}
and then $\Re f$ and $\Im f$ will again be meromorphic in $D$ and equal to
the real (resp., imaginary) part of $f$ on $\R$.

This observation allows to assume without loss of generality that all
generators of the field $\C(Y)$ are real on the real axis. Otherwise, one
should take their real and imaginary parts $\Re Y$ and $\Im Y$ as above
(both are meromorphic matrix functions real on $\R\cap D$) and consider
the field $\C(\Re Y,\Im Y)$ obviously containing $\C(Y)$. All properties of
the field $\C(Y)$ are inherited by $\C(\Re Y,\Im Y)$.

For a multivalued function $f$ ramified over the origin $t=0\in D$, taking
its real or imaginary part on the whole real axis is an ambiguous
operation because of the ramification: the real part on the positive
semiaxis $\R_+$ can cease to be real after continuation on the negative
semiaxis $\R_-$. Yet one can often guarantee that the field as a whole is
closed by taking real/imaginary parts along any segment of the real axis.

\begin{Lem}\label{lem:real-closed}
Let $A\in\Mat_{n\times n}(\R)$ be a real constant matrix and $\C(Z)$,
$Z(t)=t^A$, the corresponding Euler field. Then for any function
$f\in\C(Z)$ its real or imaginary part can be extended from the positive
semiaxis $\R_+$ or the negative semiaxis $\R_-$ to functions $g_\pm$ again
belonging to $\C(Z)$:
\begin{equation*}
  \forall f\in \C(Z)\ \exists g_\pm\in\C(Z)\:\quad
  \Re f|_{\R_\pm}=g_\pm|_{\R_\pm},
\end{equation*}
and the same for the imaginary parts $\Im f|_{\R_\pm}$.
\end{Lem}

We will denote the functions $g_\pm=\Re f_{\R_\pm}$ and their imaginary
counterparts by $\Re_\pm f$ and $\Im_\pm f$ respectively.

\begin{proof}
The field $\C(Z)$ is independent of the choice of the fundamental
solution, hence for each semiaxis $\R_\pm$ one can choose a solution
$Z_\pm$ of the Euler system $t\dot Z=AZ$ that is real on that semiaxis
(recall that the matrix $A$ is real so $A(t)=t^{-1}A$ is real-valued on
$\R$). Then it remains only to define $g_\pm$ using the identities
\begin{equation*}
  \Re\bigl(\sum c_\alpha Z_\pm^\alpha\bigr)=\sum_\alpha (\Re
  c_\alpha)Z_\pm^\alpha\in\C[Z],\qquad c_\alpha\in\C,
\end{equation*}
and similarly for the imaginary part (with obvious modifications for the
field of fractions $\C(Z)$).
\end{proof}

\begin{Rem}
The fact that the matrix $A$ is real, is not a restriction: one can always
consider a bigger Euler field of dimension $2n$ with the real
block-diagonal matrix $\operatorname{diag}(\Re A,\Im A)$.
\end{Rem}

\subsection{Variation of argument and zeros of imaginary part}
The following elementary statement will play the central role in the
constructions below.

\begin{Lem}\label{lem:im-arg}
If $f\:\R\supset [a,b]\to\C$ is a complex-valued function having no zeros
on the real interval $[a,b]$, then variation of argument of $f$ along this
interval is no greater than $\pi(\#\{\Im f=0\}+1)$.
\end{Lem}

Note that here no extension from $[a,b]$ is required, hence $\Im f$ stands
for the usual imaginary part of the restriction.

\begin{proof}
If $|\Arg f(t_1)-\Arg f(t_2)|>\pi$, then the imaginary part $\Im f$ must
vanish somewhere between $t_1$ and $t_2$ on $[a,b]$.
\end{proof}

Actually, this result is true for any linear combination $\alpha\Re
f+\beta\Im f$. The function $f$ may be allowed to have isolated zeros (but
certainly not vanishing identically), provided that zeros of $\Im f$ are
counted with multiplicities.

\subsection{Isomonodromic reduction: the local case}
Now everything is ready to prove that quasialgebraicity of $\C(Z)$ implies
the restricted quasialgebraicity of $\C(Z,Y)$ in $D=\{|t|<1\}$, provided
that $Y(t)$ is single-valued in $D$ and satisfies a Fuchsian system with
bounded residues and all singularities $\e$-distant from $\partial D$.

Any function $f\in\C(Z,Y)$ can be written in $D$ as
\begin{equation}\label{sum1N}
  f(t)=\sum_1^N z_i(t)y_i(t),\qquad z_i\in\C(Z),\ y_i\in\C_D(Y),\
  i=1,\dots,N,
\end{equation}
where $y_i$ are \emph{monomials} from the field $\C(Y)$ and hence real on
$\R$, and the number of terms $N$ as well as the degrees $\deg z_i,\deg
y_i$ (with respect to the corresponding fields $\C(Z)$ and $\C_D(Y)$) are
explicitly computable in terms of $\deg f$ with respect to
$\C(X)\subset\C(Z,Y)$. Indeed, for each generator $X_{ij}$ of degree $1$
in $\C(X)$ such representation can be chosen with no more than $n$ terms,
each being bilinear in $Y$ and $Z$.

Assume that the field $\C(Z)$ is  quasialgebraic. Then division by $z_1$
changes the number of zeros by no more than some known number. Hence when
counting the number of zeros of $f$ in any triangle $T\subset D\ssm\{0\}$,
one can assume that $z_1(t)\equiv1$, that is, $f=y_1+\sum_2^N z_i y_i$.

Also without loss of generality we may assume that $T$ does not intersect
the real axis (this can be also always achieved by subdividing it into
smaller triangles if necessary), and belongs to the upper half-plane. Then
the number of zeros of $f$ in $T$ is no greater than the sum of four
terms, variation of argument of $f$ along the sufficiently small arc
$\{|t|=\e_0\ll 1,\ \Im t>0\}$ around the origin, and similar contributions
from the semicircle $\{|t|=1,\ \Im t>0\}$, and  two rectilinear intervals
$[-1,-\e_0)\subset\R_-$ and $(\e_0,1]\subset\R_+$ respectively.

The contributions from the two circular arcs are computable functions by
Lemma~\ref{lem:joint-argument}. As for the two rectilinear segments, the
contribution of each of them is bounded by the number of zeros of $\Im f$
on them. By the assumption on the matrix function $Y$, $\Im y_i\equiv0$,
hence
\begin{equation*}
  \Im_\pm f=\Im y_1+\sum_2^N y_i\Im z_i=\sum_2^N y_i\Im_\pm z_i.
\end{equation*}
By Lemma~\ref{lem:real-closed}, the field $\C(Z)$ is closed by taking
imaginary parts, hence for all $i$ the functions $\Im_\pm z_i$ are again in
$\C(Z)$. Thus the question on the number of zeros of $f$ represented by
\ref{sum1N} in $T$ is reduced to that for the number of zeros of two
functions $\Im_+f,\Im_-f\in\C(Z,Y)$ (one for each semiaxis), each of them
involving less terms (at most $N-1$). This allows to continue the process
inductively, reducing the problem for zeros of $f$ to that for some $2^N$
functions from $\C(Z)$. On the last step the number of isolated zeros of a
product $z_N y_N$ is majorized by the sum of zeros of each term (known
since quasialgebraicity of $\C(Z)$ and $\C(Y)$ is already established).

\begin{Rem}
The above described construction reducing quasialgebraicity of the joint
field $\C_D(Z,Y)$ to that of $\C_D(X)$ if $Y$ has trivial monodromy in $D$,
is very similar to the standard differentiation-division scheme (based on
the Rolle lemma) used to obtain bounds on the number of isolated zeros of
real functions. Here the role of the differentiation is played by the
operators $\Im_\pm$, whereas single-valued functions real on $\R$ play the
role of constants killed by differentiation.
\end{Rem}

\subsection{Non-uniform quasialgebraicity}
From Theorem~\ref{thm:quasialg-local} one can immediately derive
non-uniform quasialgebraicity of Fuchsian systems. Consider such a system
with $d$ singular points.

Notice that the residue matrices $A_j$ are invariant by conformal
automorphisms of the independent variable, hence the residual norm of the
corresponding matrix function $R=\max_{i=1,\dots,d}\|A_i\|$ is also
invariant.

Using such conformal automorphisms, one can always place any three poles
of $A(t)$ at any three points, say, $0,1$ and $\infty$, but starting from
the fourth pole, one has a nontrivial parameter characterizing the spread
of singular points on the sphere. Let $\rho$ be a small positive number
such that:
\begin{enumerate}
  \item $|t_i-t_j|\ge 2\rho$, $i,j=1,\dots,d$, $i\ne j$;
  \item $|t_i|\le 1/\rho$, $i=1,\dots,d$.
\end{enumerate}

Consider the monodromy matrices $M_j$ corresponding to small loops going
around $t_j$. As before, their spectra $S_j$ are uniquely defined.

\begin{Thm}\label{thm:quasialg-restr}
If all spectra $S_j$ belong to the unit circle, then the field $\C(X)$ is
quasialgebraic. The bound for the number of zeros can be given in terms of
$R$ and $\rho$.
\end{Thm}

\begin{proof}
Draw disjoint circles of radius $\rho$ around each singularity and of
radius $1/\rho$ around the origin (this circle bounds a neighborhood of
infinity on $\C P^1$). Restricted on each circle, the field is
quasialgebraic by Theorem~\ref{thm:quasialg-local}. On the complement
there are no singularities, so after triangulation of this multiply
connected domain one can apply Corollary~\ref{cor:quasialg-nonsing}.
\end{proof}

\section{Uniform quasialgebraicity of Fuchsian systems}

In this section we briefly explain additional work to be done in order to
obtain the bounds for quasialgebraicity, that would be independent on the
relative position of singular points. In other words, we look for bounds
that would remain explicit and uniform over $\rho$ as the latter tends to
zero.

\subsection{Isomonodromic reduction: the general case}
The constructions of the previous section can be easily modified for a more
general situation. Assume that two Fuchsian systems with fundamental
solutions $X(t)$ and $Z(t)$ are \emph{isomonodromic relative to a domain
$U$}, that is, they have the same singular locus $\S$ in $U$ and the
monodromy matrices $M_\gamma$ for all loops entirely belonging to $U$, are
the same for the two. Next, assume that the residues (all of them,
including those at singular points outside $U$) are all explicitly
bounded. Finally, assume that $U$ is a polygonal domain (say, a triangle)
and the boundary $\partial U$ is away from all singularities.

Then one can claim that the two fields $\C(X)$ and $\C(Z)$ when restricted
on $U$ are both quasialgebraic or not quasialgebraic simultaneously.

The proof in the case when all singularities in $U$ fall on one straight
line, is very similar to the local case. Namely, consider the matrix
fraction $Y(t)=X(t)Z^{-1}(t)$ possessing trivial monodromy in $U$ and
embed $\C(X)$ into the joint field $\C(Z,Y)$. Writing elements of this
joint field as $\sum y_i z_i$ and applying the above algorithm of
alternating division and taking the imaginary parts, the question on the
number of zeros of $f\in\C(X)$ can be reduced to that for several
auxiliary functions from $\C(Z)$.

\subsection{Inductive strategy}\label{sec:induction}
The isomonodromic reduction principle as described above, would allow for
an inductive proof of the uniform quasialgebraicity of Fuchsian systems if
one could always construct a Fuchsian system that would be isomonodromic
to a given one in a specified simply connected domain, while having no
other singularities outside this domain.

The inductive proof may look as follows. For Fuchsian system with only two
singularities (Euler systems) the quasialgebraicity is known. Assume that
it is already established for all systems with less than $d$ finite
singularities, and consider a system with $d$ finite singular points
forming the locus $\S\subset\C$.

As was already noticed, it is impossible to make a conformal
transformation placing all $d$ points of $\S$ well apart from each other.
However, one can always achieve a situation when all finite singular points
form a set of diameter \emph{exactly} $1$ inside the disk of radius $1$
centered at the origin. In this case one can draw a line that is at least
$1/2d$-distant from all points of $\S$ and such that to each side of this
line lies at least one (hence at most $d-1$) point(s) of $\S$.

Now one can easily construct two polygonal domains each containing no more
than $d-1$ points of $\S$, together covering the whole of $\S$ and with
boundaries distant from both $\S$ and the infinity.

Assume that for each such domain $U$ a Fuchsian system can be found so that
it will be isomonodromic with the given one in $U$, while still having the
residual norm bounded in terms of the residual norm of the initial system.

Then application of the isomonodromic reduction principle would allow to
reduce the question on quasialgebraicity of the initial system in $U$ to
that for a Fuchsian system with $\le d-1$ finite singular points. By the
inductive assumption, the latter question can be explicitly answered.

The only assumption to monitor along this inductive process, is that on
eigenvalues of the monodromy operators. Clearly, one should assume that
the spectral condition (on unit absolute values of eigenvalues) holds for
all small loops around singular points. However, it is not sufficient,
since this condition does not survive the above surgery (cutting out part
of the singularities inside $U$ and pasting out the rest). Indeed, after
replacing all singularities outside $U$ by one singular point at infinity,
we create a point whose local monodromy coincides with that of the
boundary $\partial U$. Thus one must additionally assume that the
monodromy along the boundary of $U$ must also satisfy the spectral
condition.

Unfortunately, there is no way to predict how the partition into distant
``clusters'' of singular points will proceed when carried out inductively.
Instead it is sufficient to assume that the spectral condition is
satisfied for any simple loops (geometrically non-selfintersecting closed
Jordan curves).

This conditional construction would prove the following theorem.

\begin{Thm}\label{thm:quasialg-uniform}
If the monodromy operators along all simple loops have only eigenvalues of
modulus $1$, then the Fuchsian system is uniformly quasialgebraic. Upper
bounds for the number of isolated zeros can be given in terms of the
number of singular points and the residual norm of the system, uniformly
over all configurations of the singularities.
\end{Thm}

However, the proof above gives only a general idea of how the actual proof
is organized. The difficulties are of two kinds, technical and fundamental.

An example of the technical problem is the isomonodromic reduction
principle: it was formulated (and is actually proved) for a particular
case when all singular points lie on just one straight line (which can
then be identified with the real axis). However, it is sufficient for the
purposes of the proof after suitable preparation of the initial Fuchsian
system.

The fundamental problem concerns the (im)possibility of constructing a
Fuchsian system with the prescribed monodromy. It was recently discovered
that there exist obstructions to solvability of this problem, that have to
be somehow circumvented. In addition, one has  to redress the proof of the
corresponding positive results so that they would become
\emph{constructive}, yielding bounds for the residual norm of the
constructed systems.

\section{Quantitative Riemann--Hilbert problem of matrix factorization}

\subsection{Riemann--Hilbert problem: background}
Given $d$ distinct points $t_1,\dots,t_d$ on the Riemann sphere $\C P^1$
and $d$ invertible matrices $M_1,\dots,M_d$ satisfying the identity
$M_1\cdots M_d=E$ (the identity matrix), construct a Fuchsian system
having singular points at $t_1,\dots,t_d$ and only there, for which the
matrices $M_j$ would be monodromy factors for some fundamental solution.

This is the strongest form of the problem known as Hilbert 21st problem or
the \emph{Riemann--Hilbert problem}. Here some of many known results
concerning its solvability.

\begin{description}
  \item[Plemelj theorem] For any collection of points
  and any monodromy matrices, one can construct a linear system with all
  but one singularities Fuchsian; the last singular point is regular and
  can be made Fuchsian if the corresponding monodromy matrix is
  diagonalizable. As a corollary, one can always construct a Fuchsian
  system with one extra singular point that would be an apparent
  singularity \cite{plemelj,ai:ode}.
  \item[Bolibruch--Kostov theorem] If the monodromy group generated by the
  matrices $M_j$ is irreducible, then the Riemann--Hilbert problem is
  solvable \cite{bolibr:irred-monodr,kostov:irred-monodr}.
  \item[Bolibruch counterexamples] There exists a reducible group generated
  by three $4\times 4$-matrices, which cannot be realized as a monodromy
  group of a Fuchsian system \cite{bolibr:umn}.
\end{description}
What is missing in all these formulations, is the possibility of placing
any upper bound on the norms of the residue matrices of the corresponding
Fuchsian system (provided that the problem is solvable).

\subsection{Surgery}\label{sec:surgery}
In this section we formulate an analytic problem that is sufficient to
solve in order to construct a system isomonodromic to the given one in a
domain with all singularities away from the boundary.

Recall that we have a Fuchsian system that already possesses the poles
inside $U$ at specified points and the monodromy around these poles is as
required. What is necessary to do is to remove all finite singularities
from $\C\ssm U$, leaving only one singular point at infinity.

For simplicity we assume that $U$ is the unit disk centered at the origin,
and the annulus $K=\{\tfrac12<|t|<2\}$ is free from singular points (i.e.,
$\S\cap K=\varnothing$). Denote by $X(t)$ the (multivalued) fundamental
matrix solution. In general, $X$ is multivalued even in $K$, however, for
an appropriate constant matrix $B$ the product $W(t)=Xt^{-B}$ is a
single-valued matrix function that is holomorphic and holomorphically
invertible in $K$.

\emph{Suppose} that the matrix function $W=Xt^{-B}$, holomorphic and
invertible in $K$ can be represented as the matrix ratio of two other
matrix functions, $H_0(t)$ and $H_\infty(t)$, so that:
\begin{enumerate}
  \item $H_0$ is holomorphic and holomorphically invertible in the disk
  $D_0=\{|t|<2\}$;
  \item $H_\infty(t)$ is holomorphic and holomorphically invertible in the
  disk $D_\infty=\{|t|>\tfrac12\}$, including the point $t=\infty$;
  \item on the intersection $K=D_0\cap D_\infty$,
\begin{equation}\label{factorization}
  X(t)\,t^{-B}=H_0^{-1}(t)H_\infty(t),\qquad t\in D_0\cap D_\infty.
\end{equation}
\end{enumerate}

Then the two expressions, $X_0(t)=H_0(t) X(t)$ defined in $D_0$ and
$X_\infty(t)=H_\infty(t)\,t^B$ defined on $D_\infty$, agree on the
intersection and hence their ``logarithmic derivatives'' $B_0(t)=\dot
X_0(t) X_0^{-1}(t)$ and $B_\infty(t)=\dot X_\infty(t)X^{-1}_\infty(t)$,
coincide on $K$ and together define a single-valued meromorphic function
$B(t)$ on the entire sphere $\C P^1$.

The poles of $B(t)$ can be easily described: due to the holomorphic
invertibility of $H_{0,\infty}(t)$ they can occur either at the poles of
$A(t)$ that are inside $D_0$, or at $t=\infty$. In both cases the residues
can be easily computed: for any $t_j\in\S\cap D_0$,
\begin{equation*}
  B(t)=H_0(t_j)\cdot\frac{A_j}{t-t_j}\cdot H_0^{-1}(t_j)+\cdots,
\end{equation*}
where the dots stand for terms holomorphic at $t_j$. In a similar way, at
infinity we have
\begin{equation*}
  B(t)=H_\infty(\infty)\cdot\frac{B}{t}\cdot H_\infty^{-1}(\infty)+O(1/t^2),
\end{equation*}
where $O(1/t^2)$ stands for a holomorphic matrix function of the
corresponding growth. This means that in fact
\begin{equation}\label{global-B}
  B(t)=\sum_{t_j\in D_0}\frac{B_j}{t-t_j},
  \qquad B_j=H_0(t_j)\,{A_j}\,H_0^{-1}(t_j),
\end{equation}
and $\sum_j B_j=H_\infty(\infty)\,{B}\, H_\infty^{-1}(\infty)$.

The Fuchsian system with the matrix of coefficients $B(t)$ given by
\eqref{global-B}, would serve our purposes, since its solution $H_0 X$ has
the same monodromy in $D_0$ as the initial matrix solution $X$.

In order to complete the proof, one would have to obtain upper bounds on
the matrices $H_j=H_0(t_j)$ and their inverses $H_j^{-1}$ for all singular
points $t_j\in\S\cap D_0$. This would imply an upper bound on the residual
norm of the matrix function $B(t)$ and finally would allow for the
inductive proof as explained in \secref{sec:induction}.

\subsection{Birkhoff--Grothendieck factorization}
Unfortunately, finding factorization \eqref{factorization} satisfying all
properties above, is impossible.

\begin{Ex}\label{ex:fact-obstr}
One obstruction to holomorphic factorization can be immediately seen.
Consider the determinant $w(t)=\det W(t)$: this is a holomorphic
invertible function in the annulus $K$, and variation of argument of this
function along, say, the middle circle of the annulus is an integer number
$\nu$. If this number is nonzero, then either $\det H_0$ or $\det
H_\infty$ must have zeros and/or poles in the respective domains $D_0$
resp., $D_\infty$.

If $n=1$, i.e., all matrices are of size $1\times1$, then one can always
achieve factorization of the form
\begin{equation}\label{fact-scalar}
  W(t)=H_0^{-1}(t)\,t^\nu\,H_\infty(t)
\end{equation}
with $1\times1$-matrices $H_0,H_\infty$ holomorphically invertible in the
respective domains.
\end{Ex}

Both the positive and the negative assertions present in the above
example, admit generalization for the general $n$-dimensional matrix case.

\begin{Thm}[see \cite{krein-gohberg-rus,krein-gohberg}]\label{thm:bir-gro}
A matrix function $W(t)$ holomorphic and invertible in the annulus $K$,
can be factorized as follows,
\begin{equation}\label{fact-matr}
  W(t)=H_0^{-1}\,t^G\,H_\infty(t),\qquad
  G=\operatorname{diag}(\nu_1,\dots,\nu_n),\ \nu_i\in\mathbb Z,
\end{equation}
with the matrix factors $H_0,H_\infty$ holomorphic and invertible in $D_0$
resp., $D_\infty$.

The integer numbers $\nu_i$, called partial indices, are uniquely
determined \(the same for all representations with the above properties\).
\end{Thm}

The decomposition \eqref{fact-matr} can be used for surgery, if we
incorporate the term $t^G$ into $H_\infty(t)$. Then the matrix function
$Y(t)$ defined as $H_0 X$ in $D_0$ and as $t^G H_\infty \,t^B$ in
$D_\infty$, satisfies a linear system of ordinary differential equations
with only Fuchsian singular points in $D_0$ and a regular non-Fuchsian
point $t=\infty$. The principal Laurent part of the matrix $B(t)=\dot Y
Y^{-1}$ at $y=\infty$ can be easily described: the multiplicity of the pole
at infinity is bounded in terms of $\nu=\|G\|=\max_j|\nu_j|$, and the
magnitude (norm) of the Laurent coefficients of $B$ at $t=\infty$ is
bounded in terms of the Taylor coefficients of order $\le \nu$ of the
matrix $H_\infty$ and its inverse.

Actually, the form of the matrix $t^GH_\infty(t)$ is not important: in
order to have a system with a regular singularity at infinity, it is
sufficient to have a matrix factorization
\begin{equation}\label{fact-merom}
  W(t)=H_0^{-1}H_\infty,\qquad
  H_\infty(t)^{\pm 1}=C_{\pm,0}(t)+\sum_{j=1}^\nu C_{\pm,j}t^j,
\end{equation}
with holomorphic invertible matrix $H_0$ and meromorphic invertible matrix
function $H_\infty(t)$ having (together with its inverse $H_\infty^{-1}$)
the only pole at $t=\infty$ of order $\le \nu$ with the Laurent (matrix)
coefficients $C_{\pm,j}$.

\begin{Rem}
The assertion of Theorem~\ref{thm:bir-gro} is certainly not the strongest
known. The matrix function $W(t)$ can be defined only on the middle circle
$\{|t|=1\}$ of the annulus $K$ and be rather weakly regular on it, still
the factorization will be possible then, with the terms $H_0,H_\infty$
holomorphic invertible inside (resp., outside) the circle. Some
non-circular contours can be also allowed.
\end{Rem}

\subsection{Modified surgery}
One can easily modify the structure of the inductive construction above to
cover the case of systems having only Fuchsian finite singular points and
a regular eventually non-Fuchsian singularity at $t=\infty$. Such system
has the matrix of coefficients that can be written always as
\begin{equation}\label{class}
  A(t)=\sum_{j=1}^d \frac{A_j}{t-t_j}+\sum_{i=0}^\nu B_i t^i,
\end{equation}
and the residual norm for such systems should be defined as
\begin{equation}\label{class-norm}
  \|A(\cdot)\|=\sum_{j}\|A_j\|+\sum_i\|B_i\|.
\end{equation}

The surgery described in \secref{sec:surgery}, using the factorization
\eqref{fact-merom}, allows to pass from one system from such class to
another system from the same class, having no finite singularities outside
$D_0$. In order to carry out inductively the bounds for zeros, one has to
majorize the magnitude of all Laurent coefficients of the new system in
terms of the norm $R=\|A(\cdot)\|$ of the initial system.

\subsection{Bounds}
Suppose that a system from the class \eqref{class} of explicitly bound\-ed
norm $R=\|A(\cdot)\|$ has no singularities in the annulus
$K=\{\tfrac12<|t|<2\}$. Then the following bounds can be explicitly
computed.
\begin{enumerate}
  \item the norm of the monodromy $M$ corresponding to the circle
  $\{|t|=1\}\subset K$ and its matrix logarithm $B$;
  \item the pointwise upper bound for norm of the fundamental
  solution $X(t)$ with $X(1)=E$ and its
  inverse $X^{-1}(t)$ in any  smaller annulus $K'$, say,
  $\{\frac23<|t|<\frac32\}$;
  \item the pointwise upper bound for $\|W(t)\|+\|W^{-1}(t)\|$ in the
  smaller annulus $K'$.
\end{enumerate}

In order to estimate the Laurent coefficients at all singular points after
the surgery, it would be sufficient to find factorization
\eqref{fact-merom} and supply the following bounds,
\begin{enumerate}
  \item $\max_{t\in D_0}\|H_0(t)\|+\|H_0^{-1}(t)\|$ (this would allow to
  estimate the norms of residues at all finite singularities),
  \item the bound $\nu$ for the order of the pole of $H^{\pm 1}_\infty$;
  \item $\max_{t\in D_\infty}\|C_{+,0}(t)\|+\|C_{-,0}(t)\|$ together with
  \item $\sum_{i=1}^\nu\|C_{\pm,i}\|$ to majorize the norms of all Laurent
  coefficients at infinity.
\end{enumerate}

\subsection{Quantitative matrix factorization}
Unfortunately for our pur\-pos\-es, the known methods of constructing the
Birkhoff--Grothendieck factorization \eqref{fact-matr}, based on index
theory for integral operators, do not allow for quantitative conclusions.
Moreover, in some sense the problem admits no solution. The reason for
this is the known \emph{instability} of the partial indices
$\nu_1,\dots,\nu_n$: an arbitrarily small variation (in the uniform norm)
of the function $W$ can result in a jump of the partial indices. This is
clearly incompatible with existence of any bounds that would be continuous
in the $C^0$-norm.

However, if we give up with the uniquely defined Birkhoff--Grothendieck
decomposition, then one can satisfy all the above conditions.

\begin{Thm}[Novikov and Yakovenko \cite{factor}]\label{thm:eff-factorization}
A $n\times n$-matrix function $W(t)$ holomorphic and holomorphically
invertible in the annulus
\begin{equation}\label{annulus-e}
  K=\{(1+2\e)^{-1}<|t|<(1+2\e)\},\qquad \e>0,
\end{equation}
and bounded together with its inverse there,
\begin{equation}\label{bounds-annulus}
  \|W(t)\|+\|W^{-1}(t)\|<R<+\infty,\qquad t\in K,
\end{equation}
can be factorized as $W(t)=H_0^{-1}H_\infty(t)$ with the matrix functions
$H_0,H_\infty$ satisfying the following conditions.
\begin{enumerate}
  \item $H_0(t)$ is holomorphic invertible in the disk
  $D_0'=\{|t|<(1+\e)\}$ and satisfies the inequality
  $\|H_0(t)\|+\|H_0^{-1}(t)\|\le R'$ in this disk,
  \item $H_\infty(t)$ is holomorphic and holomorphically invertible in the
  complement $D_\infty'=\{(1+\e)^{-1}<|t|<+\infty\}$ and both $H_\infty$
  and $H_\infty^{-1}$ have at most a pole of order $\nu$ at $t=\infty$;
  \item the coefficients $C_{\pm,i}$ of the Laurent expansions
\begin{equation*}
  P_\pm(t)=H_\infty^{\pm 1}(t)=\sum_{i=0}^\nu C_{\pm,i}t^i
\end{equation*}
  are all bounded by $R'$ in the sense of the matrix norm;
  \item the ``regular parts'' are bounded uniformly in
  $D_\infty'$ so that $\|H_\infty^+-P_+(t)\|+\|H_\infty^--P_-(t)\|\le R'$
  there.
\end{enumerate}
The integer parameter $\nu$ and real parameter $R'$ can be expressed as
explicit elementary functions of $n,R$ and $\e$, the ``width'' of the
annulus $K$.
\end{Thm}

\subsection{Conclusion}
Theorem~\ref{thm:eff-factorization} provides the last tool necessary to
run the inductive proof and construct an explicit primitive recursive
upper bound for the number of isolated roots of solutions to Fuchsian
systems. The bounds can be explicitly written down in the sense that all
primitive recursions describing them, can be extracted from the
constructions above. However, this has not been done, among other things,
because of the very excessive bounds that appear on this way. Besides
nested inductive constructions, this occurs because of very excessive
bounds on the lengths of ascending chains of polynomial ideals (see
\secref{sec:meandering}). However, some very recent works of A.~Grigoriev
suggest that this crucial step can be considerably improved and instead of
a tower of four stories, something like a double exponential bound can be
achieved, at least for linear systems with rational coefficients. This
improvement may affect all other parts of the global construction as well,
resulting in a considerably more realistic quasialgebraicity-type
statements.

\subsection*{Acknowledgements}
The main theme of these lecture notes was surveying several recent results
obtained jointly with Dmitry Novikov, who discussed with me the text on
numerous occasions and greatly contributed to putting it to the final form.

I am grateful to Dana Schlomiuk who invited me to deliver these lectures
on the Workshop organized in Universit\'e de Montreal in June--July 2000.
Part of these notes was written during my sabbatical stay in University of
Toronto. During these periods I tried various ways of exposing parts of
these notes and collected remarks and advises from many listeners, among
them
     A.~Bolibruch,
    Yu.~Ilyashenko,
     M.~Jakobson,
     V.~Katsnelson,
     A.~Khovanskii,
     V.~Matsaev,
    Ch.~Miller,
     P.~Milman,
     B.~Schapiro,
     P.~Speisegger,
     Y.~Yomdin.

Finally, I am very much indebted to the colleagues who read preliminary
versions of these notes and pointed to numerous errors, especially to
Christiane Rousseau and Iliya Iliev. Of course, the responsibility for the
remaining bugs remains entirely with the author.

%\bibliographystyle{amsalpha}
%\bibliography{tangent16}

\def\BbbR{$\mathbf R$}\def\BbbC{$\mathbf
  C$}\providecommand\cprime{$'$}\providecommand\mhy{--}\font\cyr=wncyr9
\providecommand{\bysame}{\leavevmode\hbox
to3em{\hrulefill}\thinspace}
\providecommand{\MR}{\relax\ifhmode\unskip\space\fi MR }
\providecommand{\MRhref}[2]{%
  \href{http://www.ams.org/mathscinet-getitem?mr=#1}{#2}
} \providecommand{\href}[2]{#2}

\end{document}